\documentclass[nonblindrev]{informs3-vvm}

\SingleSpacedXI

\usepackage{endnotes}
\let\footnote=\endnote

\usepackage{url}
\usepackage{amsmath}
\usepackage{amssymb}
\usepackage{natbib}
\bibpunct[, ]{(}{)}{,}{a}{}{,}%
\def\bibfont{\small}%
\usepackage{tikz}
\usepackage{pgfplots}
\usetikzlibrary{arrows}
\usetikzlibrary{calc,arrows.meta,positioning}
\usepackage{ctable}
\usepackage{booktabs}
\usepackage{multirow}
\usepackage{algorithm}
\usepackage[noend]{algorithmic}
\usepackage{color}
\usepackage{graphicx}
\usepackage{enumitem}

\usepackage{mathtools}
\usepackage[hidelinks]{hyperref}

\usepackage{xspace}

\usepackage{pgfplotstable}
\usetikzlibrary{arrows.meta}
\usepgfplotslibrary{groupplots}

\usepackage[title]{appendix}
\usepackage{titletoc}
\newcommand\DoToC{%
	\startcontents
	\printcontents{1}{1}{\centering \textbf{Table of Contents}\vskip3pt\hrule\vskip5pt}
	\vskip3pt\hrule\vskip5pt
}

\def\Acal{\mathcal{A}}

\def\Dcal{\mathcal{D}}
\def\Fcal{\mathcal{F}}
\def\Gcal{\mathcal{G}}

\def\Scal{\mathcal{S}}
\def\Wcal{\mathcal{W}}

\def\GMGaussianNum{2\xspace}
\def\GMUniformSphereNum{2-U\xspace}
\def\GMDirichletNum{1\xspace}
\def\GMBernoulliCoveringNum{3\xspace}

\def\Ibb{\mathbb{I}}

\def\zerob{\mathbf{0}}
\def\oneb{\mathbf{1}}

\def\port{\text{portfolio}}

\def\distr{\text{distr}}

\def\Incremental{I}
\def\MNL{\text{MNL}}
\def\BiasedIncremental{BI}
\def\UniformRS{U}

\def\CR{\text{CR}}

\def\Dirichlet{\text{Dirichlet}}
\def\Normal{\text{Normal}}
\def\Bernoulli{\text{Bernoulli}}
\def\Binomial{\text{Binomial}}
\def\Uniform{\text{Uniform}}

\def\rank{\mathrm{rank}}

\def\feas{\text{feas}}
\def\covering{\text{covering}}
\def\packing{\text{packing}}

\def\Nbb{\mathbb{N}}
\def\Rbb{\mathbb{R}}
\def\Ebb{\mathbb{E}}

\def\Pr{\textbf{Pr}}

\def\Ab{\mathbf{A}}

\def\Eb{\mathbf{E}}
\def\Pb{\mathbf{P}}
\def\Ib{\mathbf{I}}
\def\ab{\mathbf{a}}
\def\bb{\mathbf{b}}
\def\cb{\mathbf{c}}
\def\db{\mathbf{d}}
\def\eb{\mathbf{e}}

\def\pb{\mathbf{p}}
\def\Pb{\mathbf{P}}
\def\rb{\mathbf{r}}
\def\sb{\mathbf{s}}
\def\Tb{\mathbf{T}}
\def\ub{\mathbf{u}}
\def\vb{\mathbf{v}}
\def\wb{\mathbf{w}}

\def\xb{\mathbf{x}}
\def\tildexb{\mathbf{\tilde{x}}}
\def\yb{\mathbf{y}}
\def\zb{\mathbf{z}}
\def\Zb{\mathbf{Z}}

\def\alphab{\boldsymbol \alpha}

\def\thetab{\boldsymbol \theta}
\def\epsilonb{\boldsymbol \epsilon}

\def\xib{\boldsymbol \xi}

\def\Sigmab{\boldsymbol \Sigma}
\def\lambdab{\boldsymbol \lambda}

\def\OLP{\mathrm{OLP}}
\def\OLA{\mathrm{OLA}}
\def\worep{\text{wo-rep}}

\def\Misic{Mi\v{s}i\'{c} }

\def\Halmos{$\square$}

\usepackage{natbib}
 \bibpunct[, ]{(}{)}{,}{a}{}{,}%
 \def\bibfont{\small}%

\TheoremsNumberedThrough     %
\ECRepeatTheorems

\EquationsNumberedThrough    %

\begin{document}

\RUNAUTHOR{Akchen and \Misic}

\RUNTITLE{Column-Randomized Linear Programs}

\TITLE{Column-Randomized Linear Programs: Performance Guarantees and Applications}

\ARTICLEAUTHORS{%
	\AUTHOR{Yi-Chun Akchen}
	\AFF{School of Management, University College London, London E14 5AB, United Kingdom, \EMAIL{ \tt yi-chun.akchen@ucl.ac.uk}}
	\AUTHOR{Velibor V. Mi\v{s}i\'{c}}
	\AFF{UCLA Anderson School of Management, University of California, Los Angeles, California 90095, United States, \EMAIL{\tt velibor.misic@anderson.ucla.edu}} %
}
\ABSTRACT{We propose a randomized method for solving linear programs with a large number of columns but a relatively small number of constraints. Since enumerating all the columns is usually unrealistic, such linear programs are commonly solved by column generation, which is often still computationally challenging due to the intractability of the subproblem in many applications. Instead of iteratively introducing one column at a time as in column generation, our proposed method involves sampling a collection of columns according to a user-specified randomization scheme and solving the linear program consisting of the sampled columns. While similar methods for solving large-scale linear programs by sampling columns (or, equivalently, sampling constraints in the dual) have been proposed in the literature, in this paper we derive an upper bound on the optimality gap that holds with high probability.  This bound converges at a rate $1 / \sqrt{K}$, where $K$ is the number of sampled columns, to the optimality gap of a linear program related to the sampling distribution. We analyze the gap of this latter linear program, which we dub the distributional counterpart, and derive conditions under which this gap will be small. Finally, we numerically demonstrate the effectiveness of the proposed method in the cutting-stock problem and in nonparametric choice model estimation. 
}%


 \HISTORY{First version: July 20, 2020. Second version: July 28, 2022. Third version: June 21, 2023. Forthcoming in {\it Operations Research}.}

\maketitle

\section{Introduction}
We consider solving a linear program (LP) in standard form:
\begin{subequations}
	\label{problem:LP_standard}
	\begin{alignat}{2}
	& \underset{\xb \in \Rbb^n}{\text{minimize}}  \quad & & \cb^T\xb \\   
	& \text{such that}  \quad & & \Ab \xb = \bb, \label{eq:LP_standard_Axeqb} \\  
	& & & \xb \geq \zerob, \label{eq:LP_standard_xgeq0}
	\end{alignat}
\end{subequations}
where $\xb \in \Rbb^n$, $\cb \in \Rbb^n$, $\Ab \in \Rbb^{m \times n}$, and $\bb \in \Rbb^m$. In various applications of linear programming, such as the cutting-stock problem \citep{gilmore1961linear}, the vehicle routing problem \citep{dumas1991pickup}, and the choice-based network revenue management \citep{bront2009column}, it is often the case that the number of variables $n$ is significantly larger than the number of constraints $m$. For example, in the choice-based network revenue management problem, $n = O(2^N)$ and $m = O(N)$, where $N$ is the number of products in a market and each column corresponds to a subset of these $N$ products. In the cutting stock problem, the number of columns $n$ represents the number of feasible cutting patterns and it grows exponentially with respect to $m$. In both cases, one can easily observe instances of LP \eqref{problem:LP_standard} such that $m$ is on the scale of a few hundred while $n$ is more than a billion. Moreover, since the constraint matrix $\Ab$ is too large in these large-scale LPs, one usually cannot explicitly write it down but only specify it as a matrix that consists of all columns satisfying a certain property.

Given that there are many more columns than constraints and enumerating all of the columns, i.e., obtaining the full constraint matrix $\Ab$, is impossible in most cases, a standard solution method is column generation (CG), which works as follows: (i) start with an initial set of columns from $\Ab$; (ii) solve the corresponding restricted linear program to optimality; (iii) solve a subproblem to find the column with the lowest reduced cost; (iv) add the new column to the current set of columns; (v) go back to step (i) until problem~\eqref{problem:LP_standard} is solved to optimality (i.e., the minimum reduced cost in step (iii) is nonnegative). The subproblem that ones solves to introduce a new column is often computationally challenging. For example, in the cutting-stock problem, which is a well-known large-scale LP that is typically solved using column generation, the subproblem is a knapsack problem that is known to be NP-hard \citep{garey1979computers}. In practice, the subproblem is often formulated as an integer program, and can be difficult to solve at a large scale.  %
In addition, CG is a \emph{sequential} method, that is, the subproblem that one solves to introduce the $i$th column depends on the computational results of the previous $i-1$ iterations. Such a structure prohibits one from applying parallel computing techniques to implement the column generation method.

Instead of searching for columns by a subproblem that is potentially NP-hard, we propose a randomized method, called \emph{column randomization}. In this method, one first samples a collection of columns according to a user-specified randomization scheme, and then solves the corresponding restricted linear program. We refer to this restricted linear program that consists of sampled columns as the {\it column-randomized linear program}. This approach is attractive because computationally, it is often significantly easier to randomly sample columns than it is to optimize over columns (as is the case in CG). In addition, while CG operates sequentially, the sampling step in column randomization is well-suited to parallelization. 

We note that similar sampling-based methods for large-scale LPs have been previously considered in the operations research literature. In particular, there is a significant literature on solving problems with large numbers of constraints by randomly sampling constraints \citep{de2004constraint,calafiore2005uncertain}. By strong duality of linear programs, sampling the columns of problem~\eqref{problem:LP_standard} is equivalent to sampling the constraints of its dual problem. However, the behavior of the sampled LP in terms of its optimality gap -- the difference in objective value between the sampled problem and the complete problem -- has received scarce attention in the literature. In this paper, our main goal is to answer the following question: {\it Given a user-specified randomization scheme for sampling columns from a linear program, is it possible to probabilistically bound the optimality gap of the column-randomized linear program?}

We provide theoretical results to answer this question and demonstrate how these results can be applied to common applications of large-scale linear programming. We make the following specific contributions:
\begin{enumerate}%
	\item {\bf Theoretical Guarantees}. We show that with high probability over the sample of columns, the optimality gap of the column-randomized linear program is bounded by the sum of two terms: the optimality gap of a linear program related to the sampling distribution and a term that is of order $1/\sqrt{K}$, where $K$ is the number of sampled columns. To best of our knowledge, this is the first simple theoretical result that addresses the behavior of the optimality gap of the column sampling technique for general linear programs using only elementary arguments (in particular, LP sensitivity analysis and McDiarmid's inequality).

	\item {\bf Analysis of the Distributional Counterpart}. A key component of our bound is the optimality gap of an LP related to the sampling distribution that we refer to as the \emph{distributional counterpart}. We undertake a detailed analysis of this quantity. We show theoretically that this gap will be small when there exist many diverse near-optimal basic feasible solutions, where diversity is measured by how infrequently a column appears in the bases. We also study this gap in a probabilistic setting, where we assume that the LP~\eqref{problem:LP_standard} is generated according to a random generative model. We show that under three different generative models, the distributional counterpart gap scales like $O(1/\sqrt{n})$ or $O(\log n / \sqrt{n})$ with high probability, where $n$ is the number of columns in the complete LP.

	\item {\bf Extensions}. We extend our main performance guarantee in two ways. First, we apply the proposed method to several applications of large-scale linear programming and derive problem-specific upper bounds for the optimality gap. The problems include LPs with totally unimodular constraints, Markov decision processes (MDP), covering problems and packing problems. We also extend our approach to the portfolio optimization problem, in which the objective function is only assumed to be Lipschitz continuous (and is not necessarily linear or convex). Second, we generalize our column randomization approach to the case where the sampled columns are no longer i.i.d. and may be statistically dependent. In particular, we develop a theoretical guarantee for when the dependency of the sampled columns is described by a dependency graph, and an alternate guarantee for the case when columns are sampled uniformly without replacement.
	
	\item {\bf Numerical Results}. We numerically demonstrate the effectiveness of the proposed method on two optimization problems that are commonly solved by CG: the cutting-stock problem, which is a classical application of linear programming; and the nonparametric choice model estimation problem, which is a modern application of linear programming. We compare the performance of the column randomization method to that of the CG method and show that for a fixed positive optimality gap, the column randomization method can attain the same optimality gap within a fraction of the time required by CG. Thus, for some problems, the column randomization method can be a viable alternative to CG or can otherwise be used to provide a good warm start solution for CG. 
	\end{enumerate}

We organize the paper as follows. In Section \ref{sec:literature_review}, we review the related literature and highlight our contribution.  In Section \ref{sec:main_results}, we state our theoretical results and discuss their implications. In Section~\ref{sec:P_distr}, we present our detailed analysis of the distributional counterpart. Due to space constraints, our extensions to special problem structures/applications of large-scale LP and to non-i.i.d. column sampling are presented in Sections~\ref{sec:special_structures} and \ref{sec:dependent_columns} of the ecompanion. Section~\ref{sec:numerics_CS} presents our numerical experiments with the cutting stock problem and Section~\ref{sec:numerics_NCME} presents our numerical experiments with the nonparametric choice model estimation problem. We conclude in Section~\ref{sec:conclusion}. Omitted proofs are provided in the electronic companion.

\section{Literature Review}
\label{sec:literature_review}

In this section, we review four streams of literature. %

{\bf Column Generation.} CG has been widely used to solve optimization problems that have a huge number of columns compared to the number of constraints \citep{ford1958suggested,dantzig1960decomposition,du1999stabilized}. Applications include vehicle routing \citep{dumas1991pickup,feillet2010tutorial}, facility location problems \citep{klose2005lower}, and choice model estimation \citep{van2015market,mivsic2016data}; we refer readers to \cite{desrosiers2005primer} for a comprehensive review. By strong duality of linear programs, CG is equivalent to constraint generation that solves linear programs with a large number of constraints \citep{bertsimas1997introduction}. A key component of both methods is the subproblem that one solves to iteratively introduce columns or constraints. Usually, this subproblem is computationally challenging and is often solved by integer programming. For example, in the cutting-stock problem, the CG subproblem is a knapsack problem, which is NP-hard \citep{gilmore1961linear,garey1979computers}. 

{\bf Sampling Columns/Constraints.} Another approach to solving LPs with huge numbers of columns (or equivalently, with huge numbers of constraints), is by sampling \citep{de2004constraint,calafiore2005uncertain,calafiore2006scenario,campi2008exact,campi2018wait}. Specifically, one first samples a set of columns (or constraints) according to a given distribution then solves a linear program that consists of the sampled columns (or constraints). The seminal paper of \cite{de2004constraint} proposed the constraint sampling method for linear programs that arise in approximate dynamic programming (ADP). Given a distribution for sampling the constraints, the paper showed that with high probability over the sampled set of constraints, any feasible solution of the sampled problem is nearly feasible for the complete problem (that is, there is a high probability of satisfying a new random constraint, sampled according to the same distribution). Under the additional assumption that the constraint sampling distribution is a Lyapunov function, the paper also develops a specific guarantee on the error between the optimal value function and the approximate value function that is obtained by solving the sampled problem, but does not provide a bound on the gap between the objective values of the sampled and complete linear programs. 
In contrast, the results of our paper pertain specifically to the objective value of the sampled problem, are free from any assumptions on the sampling distribution and are applicable to general linear programs beyond those arising in ADP. Around the same period, \cite{calafiore2005uncertain,calafiore2006scenario} pioneered the sampling approach to robust convex optimization. With a different perspective from \cite{de2004constraint}, \cite{calafiore2005uncertain,calafiore2006scenario} also characterized the sample complexity needed for the optimal solution (as opposed to an arbitrary feasible solution) of the sampled problem to be nearly feasible for the original problem. However, the performance of the sampled problem in terms of the objective value, and its dependence on the number of samples, was not addressed.

Since the works of \cite{calafiore2005uncertain} and \cite{de2004constraint}, there has been some work that has quantified the dependence of the objective value on the number of sampled constraints. In particular, the paper of \cite{esfahani2014performance} considers a convex program where the decision variable $\xb$ satisfies a family of convex constraints, which are later sampled, and is also constrained to lie in an ambient set $\mathbb{X}$. The paper develops a probabilistic bound on the difference in objective value between the complete problem and its sampled counterpart in terms of a uniform level-set bound (ULB), which is a quantile function of the worst-case probability over all feasible solutions in set $\mathbb{X}$. Our work differs significantly from \cite{esfahani2014performance} in two aspects. First, in terms of the problem setting, \cite{esfahani2014performance} assumes that even before any constraints are sampled, the decision variable is already constrained in the convex compact (and thus bounded) set $\mathbb{X}$, and the associated performance guarantees also rely on properties of $\mathbb{X}$. In our setting, this corresponds to the dual solutions of problem~\eqref{problem:LP_standard} being bounded, which need not be the case in general. Consequently, the result of \cite{esfahani2014performance} is not directly applicable to the research question discussed in this paper. Second, as noted earlier, the performance bound in \cite{esfahani2014performance} relies on the ULB function of the sampling distribution. While sufficient conditions for the existence of a ULB are provided in the paper, in general a ULB cannot be represented explicitly and thus the resulting performance guarantee is less interpretable. In contrast, our theoretical results do not require a ULB or other related functions, and have a more interpretable dependence on the sampling distribution (via the distributional counterpart; see problem~\eqref{problem:LP_restricted} in Theorem~\ref{thm:main_largest_abs_of_all_dual_BFS}). In addition, we also believe our results are more straightforward technically: one only needs McDiarmid's inequality and standard linear programming results to prove them. As we will show in Section~\ref{sec:special_structures}, our theoretical results and proof technique can be applied to many common types of LPs to derive application-specific guarantees.

{\bf Randomized Projection, Stochastic Optimization and Online Linear Programming.} Besides column/constraint sampling, many other randomized methods have been proposed to solve large-scale optimization problems, including methods based on random walks \citep{bertsimas2004solving} and random projection \citep{pilanci2015randomized, vu2018random}. Specifically, the random projection method of \cite{vu2018random} involves selecting a matrix $\Tb \in \Rbb^{k \times m}$, with $k < m$ and then left-multiplying the constraint matrix and the right-hand side vector. This transforms the original LP $\min \{ \cb^T \xb \mid \Ab \xb = \bb, \xb \geq \zerob\}$ into the problem $\min\{ \cb^T \xb \mid \Tb \Ab \xb = \Tb \bb, \xb \geq \zerob\}$, which is a problem with fewer rows. Although there are a number of important differences between this approach and ours, the most significant is the philosophical difference in the intended use case for each approach. In order to apply the random projection approach, one needs to be able to form the full constraint matrix $\Ab$ and the projection matrix $\Tb$ in order to carry out the multiplication $\Tb \Ab$. Consequently, the random projection approach applies to LPs where $\Ab$ is large, but not so large that it cannot be formed and stored in computer memory. In contrast, our approach can be used for truly large-scale LPs where the matrix $\Ab$ is defined implicitly, as a matrix whose columns obey some property, and cannot be formed explicitly because the number of columns $n$ is astronomically large (e.g., the set of patterns for the cutting stock problem, as in our experiments in Section~\ref{sec:numerics_CS}). For this type of large-scale LP that is ordinarily solved via column generation and that can be solved by our approach, the random projection cannot be applied, because it is computationally infeasible to form the matrix $\Ab$. We discuss this difference, and other differences, in more detail in Section~\ref{subsec:comparison_vu2018random} of the ecompanion.

In addition to these randomized methods, there is also a separate literature on optimization problems where stochasticity is part of the problem definition; some examples include stochastic programming \citep{birge2011introduction,shapiro2014lectures}, contextual optimization \citep{elmachtoub2017smart}, and online optimization \citep{shalev2012online}. Within this literature, the problem setting of online linear programming, where columns of a linear program are revealed sequentially to a decision maker, bears a resemblance to ours; some examples of papers in this area include \cite{agrawal2014dynamic,eghbali2018competitive,li2019online}. Despite this similarity, this problem setting differs significantly from ours in that a decision maker is making irrevocable decisions in an online fashion: the decision maker must decide how much to use of a variable/column at the time that it is revealed, and cannot revise this decision in the future. With regard to \cite{agrawal2014dynamic} specifically, we note that this paper comments on the possibility of applying the proposed procedure (the one-time learning algorithm, or OLA) in an offline manner. Deploying OLA in such a manner would involve iterating through all $n$ columns of the LP; although this could be accomplished when $n$ is of a medium scale, it is untenable when $n$ is combinatorially large (e.g., the number of patterns in a cutting stock problem, as in our experiments with the cutting stock problem in Section~\ref{sec:numerics_CS}, or the number of rankings in the nonparametric choice estimation problem in Section~\ref{sec:numerics_NCME}, which is $(N+1)!$ for a set of $N$ products). For a more detailed comparison of the differences between our work and \cite{agrawal2014dynamic}, we refer readers to Section~\ref{subsec:comparison_agrawal2014dynamic} of the ecompanion.

{\bf Other Related Literature.} Our proof technique is inspired by the literature on random feature selection in machine learning \citep{moosmann2007fast,rahimi2008random,rahimi2009weighted}. In particular, our paper generalizes the result of \cite{rahimi2009weighted}, which considers the problem of learning a predictive model that is a weighted sum of random feature functions, to the problem of solving linear programs that consist of random columns. The major difference between our setup and that of \cite{rahimi2009weighted} is that the decision variables in a linear program must satisfy constraints (i.e., constraints \eqref{eq:LP_standard_Axeqb} and \eqref{eq:LP_standard_xgeq0}), while the weights of random feature functions in the setup of \cite{rahimi2009weighted} are not constrained in any way. Because of this difference, the results of \cite{rahimi2009weighted} cannot directly be applied to our problem setting. To overcome this feasibility issue, we utilize classical LP sensitivity analysis and relate a possibly infeasible solution constructed using the random sample of columns to a feasible solution of the sampled LP (see Section~\ref{subsec:main_theorem_proof}).

\section{Column-Randomization Method}
\label{sec:main_results}

In this section, we first describe the basic notations and definitions that will be used throughout the paper (Section~\ref{subsec:notation}). Then we formally define the column randomization method and investigate its theoretical properties (Section~\ref{subsec:sampling_columns_code_and_theorems}). We end this section by discussing implications and interpretations of the theoretical results (Section~\ref{subsec:discussion_theorems}). Proofs of the results are relegated to Section~\ref{sec:proofs}.

\subsection{Notation and Definitions}
\label{subsec:notation}
For any positive integer $n$, let $[n] \equiv \{ 1,2,\ldots,n \}$. Let $\eb_i$ be the $i$th standard basis vector for $\Rbb^n$; that is, $\eb_i = (e_{ij})$ where $e_{i,j} = 1$ if $j = i$ and $e_{i,j} = 0$ if $j \neq i$. Thus, for any $\xb \in \Rbb^n$, we can represent it as $\xb= \sum_{i \in [n]} x_i \eb_i$. We consider a linear program in standard form:
\begin{align}
\label{problem:LP_standard_2}
P: \quad  \min \{  \cb^T \xb \mid \Ab \xb = \bb, \,\,\, \xb \geq \zerob  \},
\end{align}
where $\Ab$ is an $m \times n$ matrix and $\cb \in \Rbb^n$. We will refer to the problem $P$ as the {\it complete problem} throughout the paper, as it contains all of the columns of $\Ab$. 

We make two assumptions on problem $P$. %
First, we assume that problem $P$ is feasible and bounded; this assumption is not too restrictive, since the cases where the complete problem $P$ is either unbounded or infeasible are not interesting to consider. The second assumption we make is that $\rank(\Ab) = m$, i.e., the rows of $\Ab$ are linearly independent. This is also not too restrictive, as one can remove any rows of $\Ab$ that are linear combinations of the other rows without changing the problem. 

We define the dual problem $D$ of problem~\eqref{problem:LP_standard_2} as $\max \{  \pb^T \bb \mid \pb^T \Ab \leq  \cb^T  \}$. For any optimization problem $P'$, we denote its optimal objective value by $v(P'$) and its feasible region by $\Fcal(P')$. By LP strong duality and the assumption that $P$ is feasible and bounded, we have $v(P) = v(D)$. Furthermore, for any optimization problem $P''$ that shares the same objective function as the complete problem $P$ and satisfies $\Fcal(P'') \subseteq \Fcal(P)$, we define $\Delta v (P'') \equiv v(P'') - v(P)$, which is nonnegative and can be interpreted as the optimality gap of solving $P''$ instead of $P$.

For each $i \in [m]$ and $j \in [n]$, we use $\Ab^i$ and $\Ab_j$ to denote the $i$th row and $j$th column of matrix $\Ab$, respectively. For any collection of indices $J \subseteq [n]$, we let $\Ab_J$ represent the submatrix of $\Ab$ that consists of columns whose indices belong to $J$.  In this paper, instead of solving either the complete problem $P$ or its dual $D$, we consider solving a linear program whose columns are randomly selected. We call such a linear program a \emph{column-randomized linear program}, which we formally define below.

\begin{definition}{\it (Column-Randomized Linear Program)
	Let $J$ be a finite collection of random indices, i.e., $J \equiv \{ j_1,j_2,\ldots,j_K  \}$ for an integer $K$, where $j_k \in [n]$ is a random variable for $k = 1,2,\ldots,K$. Then the problem
	\begin{align}
	\label{problem:LP_sampled_2}
	P_J: \quad \min \{  \cb^T_J \tildexb \mid \Ab_J \tildexb = \bb, \,\,\, \tildexb \geq \zerob  \}
	\end{align}
	is called a column-randomized linear program.}
	\end{definition}

Clearly, $P_J$ is equivalent to $\min \left\lbrace \cb^T \xb \mid \Ab \xb = \bb, \,\,\, \xb \geq \zerob,\,\,\, x_j = 0 \,\,\, \forall j \notin J. \right\rbrace$. With this reformulation, any feasible solution of $P_J$ can be represented as an element in $\Fcal(P)$. We can thus define $\Delta v(P_J)$ for the column-randomized LP $P_J$. We sample random indices in $J$ by a randomization scheme $\rho$, which is a computational procedure that randomly selects indices from $[n]$, or equivalently, randomly generates columns from $\Ab$. Let $\xib$ be the probability distribution over $[n]$ that corresponds to $\rho$; that is, the $j$th component of $\xib$, denoted by $\xi_j$, is the probability that index $j$ is selected by $\rho$. Throughout this section, we assume $\rho$ samples each index independently and identically according to $\xib$. We will relax this assumption in Section~\ref{sec:dependent_columns}. We denote the dual problem $D_J$ of $P_J$ as $\max \{  \pb^T \bb  \mid \pb^T \Ab_J \leq  \cb^T_J  \}$.

We will also require the notions of a basis, basic solutions and reduced costs in our theoretical results. A collection of indices $B \subseteq [n]$ of size $m$ is called a \emph{basis} if the matrix $\Ab_B$ is nonsingular, i.e., the collection of $m$ columns $\{ \Ab_j \}_{j \in B}$ is linearly independent. A \emph{basic solution} $\xb$ of the primal problem $P$ corresponding to the basis $B$ is the solution $\xb$ obtained by setting $\xb_B = \Ab^{-1}_B \bb$, where $\xb_B$ is the subvector corresponding to the columns in $B$, and $\xb_N = \zerob$, where $\xb_N$ is the subvector corresponding to the columns in $[n] \setminus B$. A solution $\xb$ is called a \emph{basic feasible solution} of $P$ if it is a basic solution for some basis $B$ and satisfies $\xb \geq \zerob$. For the dual problem, a basic solution $\pb$ corresponding to the basis $B$ is the solution $\pb$ defined by setting $\pb^T = \cb_B^T \Ab^{-1}_B$; if it additionally satisfies $\pb^T \Ab \leq \cb^T$, then it is also a basic feasible solution. Given a basis $B$, we define the reduced cost vector $\bar{\cb}$ for that basis as $\bar{\cb} \equiv \cb^T - \cb^T_B \Ab^{-1}_B \Ab$. 

Finally, we use $\| \cdot \|$ to denote norms. For a vector $\vb \in \Rbb^n$, we let $\| \vb \|_1 = \sum_{j=1}^n |v_j|$ be its $\ell_1$ norm, $\| \vb \|_2 = \sqrt{\sum_{j=1}^n v_j^2}$ be its Euclidean or $\ell_2$ norm, and $\| \vb \|_\infty = \max_{j=1,\ldots,n}|v_j|$ be its $\ell_\infty$ norm. For a matrix $\Ab$, we let $\| \Ab \|_{\max} = \max_{i,j} | A_{i,j} |$. %
Without loss of generality, we assume that the cost vector $\cb$ has unit Euclidean norm, i.e., $\| \cb \|_2 = 1$. 
This is not a restrictive assumption, because by normalizing the cost vector $\cb$ to have unit Euclidean norm, the objectives of the complete problem $P$ and the column-randomized problem $P_J$ are both scaled by $1/\| \cb \|_2$. Thus, the relative performance of problem $P_J$ to the complete problem $P$, which is the main focus of our paper, remains the same. %

\subsection{Performance Guarantees}
\label{subsec:sampling_columns_code_and_theorems}

We propose the \emph{column randomization method} in Algorithm~\ref{alg:main}. We first sample $K$ indices, $j_1,j_2,\ldots,j_K$, by a randomization scheme $\rho$ and let $J = \{ j_1, \ldots,j_K\}$. We then collect the corresponding columns of $\Ab$ as matrix $\Ab_J$ and the corresponding components of $\cb$ as vector $\cb_J$. After forming $\Ab_J$ and $\cb_J$, we solve the LP~\eqref{problem:LP_sampled} and return its optimal value $v(P_J)$ and optimal solution.

\begin{algorithm}
\SingleSpacedXI
	\caption{The Column Randomization Method}
	\label{alg:main}
	\begin{algorithmic}[1]
		\STATE Sample $K$ indices as $J \equiv \{ j_1,\ldots,j_K \}$ by a randomization scheme $\rho$.
		\STATE Define $\Ab_J = [A_{j_1},\ldots,A_{j_K}]$ and $\cb_J = [c_{j_1},\ldots,c_{j_K}]$.
		\STATE Solve the column-randomized linear program, which only has $K$ columns:
		\begin{align}
		\label{problem:LP_sampled}
		P_J:  \quad \min \left\lbrace \cb^T_J \tildexb \mid \Ab_J \tildexb = \bb, \,\,\, \tildexb \geq \zerob \right\rbrace. 
		\end{align}
		\RETURN optimal objective value $v(P_J)$ and an optimal solution $\tildexb^*$.
	\end{algorithmic}
\end{algorithm}

Notice that an optimal solution $\tilde{\xb}^*$ of problem $P_J$ can be immediately converted to a feasible solution for the complete problem $P$ by enlarging $\tilde{\xb}^*$ to length $n$ and setting $\tilde{x}^*_j = 0$ for $j \in [n] \setminus J$. 

We now present two theorems that bound the optimality gap $\Delta v(P_J) \equiv v(P_J) - v(P)$ of problem $P_J$; we defer our discussion of these two theorems to Section~\ref{subsec:discussion_theorems}. Since several preliminary results are needed before we prove the theorems, we also relegate the proofs of the theorems to Section~\ref{sec:proofs}.

\begin{theorem}
	\label{thm:main_largest_abs_of_all_dual_BFS}
	Let $C$ be a positive constant and define the linear program $P_\distr$ as
	\begin{subequations}
		\label{problem:LP_restricted}
		\begin{alignat}{3}
		P_\distr: & \quad & \underset{\xb \in \Rbb^n}{\text{minimize}}  \quad & \cb^T\xb \\ & & \text{such that}  \quad & \Ab \xb = \bb,
		\\ & &  & \zerob \leq \xb \leq C \cdot \xib.
		\end{alignat}
	\end{subequations}
	Let $P_J$ be the column-randomized LP solved by Algorithm~\ref{alg:main}, and $\Ab_J$ be the corresponding constraint matrix. 
	For any $\delta \in (0,1)$, with probability at least $1 - \delta$ over the sample $J$, the following holds: if $P_J$ is feasible and $\rank(\Ab_J) = m$, then
	\begin{align}
	\label{eq:thm_convergence_with_abs_of_BFS}
	\Delta v(P_J)   \leq \Delta  v(P_\distr) +  \frac{C \left( 1 + m \gamma \| \Ab \|_{\max} \right)}{\sqrt{K}} \left(  1 + \sqrt{2 \log \frac{2}{\delta}}  \right),
	\end{align}
	where $\gamma$ is an upper bound on $\|\pb \|_{\infty}$ for every basic solution $\pb$ of the dual problem $D$ and $\| \Ab \|_{\max} = \max_{ij} | A_{ij}|$.
\end{theorem}

Theorem \ref{thm:main_largest_abs_of_all_dual_BFS} shows that, with probability at least $1-\delta$, the optimality gap $\Delta v(P_J)$ of the column-randomized LP $P_J$ is upper bounded by the sum of two terms. The first term is the optimality gap $\Delta v(P_\distr)$ of the problem $P_\distr$, which we refer to as the \emph{distributional counterpart}. The second term involves $\| \Ab \|_{\max}$, the largest absolute value of elements in the constraint matrix; $\gamma$, the upper bound of the $\ell_\infty$ norm of any basic solution of the dual problem; $\delta$, the confidence parameter; and $K$, the number of sampled columns. Most importantly, the second term converges to zero with a rate $1 / \sqrt{K}$. In Section~\ref{sec:special_structures}, we will see how $\gamma$ and $\| \Ab \|_{\max}$ can be further bounded for certain special cases.

We now present our second theorem, which relates the optimality gap to the reduced costs of the complete problem. 

\begin{theorem}
\label{thm:main_reduced_cost}
Define $C$, $P_\distr$, $P_J$ and $\Ab_J$ as in Theorem~\ref{thm:main_largest_abs_of_all_dual_BFS}. For any $\delta \in (0,1)$, with probability at least $1 - \delta$ over the sample $J$, the following holds: if $P_J$ is feasible and $\rank(\Ab_J) = m$, then
\begin{align}
\Delta v(P_J) \leq \Delta v(P_{\distr}) + \frac{C}{\sqrt{K}} \cdot \chi \cdot \left(1 + \sqrt{2 \log \frac{1}{\delta} } \right)
\end{align}
where $\chi$ is an upper bound on $\| \bar{\cb} \|_2$ for every basic solution of the complete problem $P$. 
\end{theorem}

Theorem~\ref{thm:main_reduced_cost} has a similar structure to Theorem~\ref{thm:main_largest_abs_of_all_dual_BFS}. Compared to Theorem~\ref{thm:main_largest_abs_of_all_dual_BFS}, the upper bound in Theorem~\ref{thm:main_reduced_cost} does not involve $\gamma$ and $\| \Ab \|_{\max}$, but instead requires a bound on the norm of the reduced cost vector for all the bases of $P$.

\subsection{Discussion}
\label{subsec:discussion_theorems}

Both Theorem~\ref{thm:main_largest_abs_of_all_dual_BFS} and \ref{thm:main_reduced_cost} provide bounds on the optimality gap $\Delta v(P_J)$ of the following form:
\begin{align}
\label{eq:bound_general_form}
\Delta v(P_J)   \leq \Delta  v(P_\distr) + \frac{C \cdot C_{P} \cdot C_{\delta}}{\sqrt{K}},
\end{align}
where $C_{P}$ only depends on properties of the complete problem $P$ and $C_{\delta}$ only depends on the confidence parameter $\delta$. In Theorem~\ref{thm:main_largest_abs_of_all_dual_BFS}, $C_{P} = 1 + m \gamma \| \Ab \|_{\max}$ and $C_\delta = 1 + \sqrt{2 \log ({2}/{\delta})}$; in Theorem~\ref{thm:main_reduced_cost}, $C_{P} = \chi$ and $C_\delta = 1 + \sqrt{2 \log ({1}/{\delta})}$. In the following discussion, we first focus on the general structure of the upper bounds given in \eqref{eq:bound_general_form}, and subsequently we address the differences between Theorem~\ref{thm:main_largest_abs_of_all_dual_BFS} and Theorem~\ref{thm:main_reduced_cost}.

\subsubsection*{Role of Problem $P_\distr$:}

The distributional counterpart $P_\distr$ is the restricted version of the complete problem $P$, which includes the additional constraint $\xb \leq C \xib$. Ignoring the value of the constant $C$, which we will discuss in more detail below, $P_{\distr}$ can be interpreted as a modification of $P$ where the most we can use each variable $j$ is proportional to the probability of that variable being sampled, $\xi_j$. Therefore, in a certain sense, $P_{\distr}$ measures how well the induced distribution $\xib$ is aligned with columns that are used in good or optimal solutions of $P$. As a crude example, for a fixed $C$, if $\xib$ is such that $\xi_j$ is large for every column $j$ that appears in an optimal basis and small for any other $j$, then we should expect $P_{\distr}$ to be small. On the other hand, for a fixed $C$, if $\xib$ is such that $\xi_j$ is small for columns $j$ that appear in optimal or near optimal bases, and large for $j$'s that appear in highly suboptimal bases, then we should expect $\Delta v(P_{\distr})$ to be large. However, this simple discussion highlights just one aspect of how $P_{\distr}$ behaves, and in general $P_{\distr}$ has a complex dependence on the structure of $P$ and the distribution $\xib$. When $\xib$ corresponds to the uniform distribution over $[n]$, $\Delta v(P_{\distr})$ can be viewed as measuring how well solutions that are non-sparse and have low $\ell_{\infty}$ norm -- i.e., solutions that will satisfy $\xb \leq C \xib = C / n$ -- perform in the problem $P$. We discuss this perspective on $P_{\distr}$ in much greater detail in Section~\ref{sec:P_distr}.

\subsubsection*{Role of Constant $C$:}

Given a randomization scheme $\rho$ and its corresponding distribution $\xib$, as the constant $C$ increases, the optimality gap $\Delta v(P_\distr)$ of problem $P_\distr$ decreases since its feasible set $\Fcal(P_\distr)$ is enlarged. On the other hand, the second term on the RHS of bound \eqref{eq:bound_general_form} increases since it is proportional to $C$. To interpret this phenomenon, we can view bound~\eqref{eq:bound_general_form} as a type of bias-complexity/bias-variance tradeoff, which is common in statistical learning theory \citep{shalev2014understanding}:
\begin{align}
\label{eq:bound_bias_complexity}
\Delta v(P_J)   \leq \underbrace{ \vphantom{  \frac{C \cdot C_{P} \cdot C_\delta }{\sqrt{K}}     }   \Delta  v(P_\distr)}_{\text{Approximation Error}} + \underbrace{\frac{C \cdot C_{P} \cdot C_{\delta}  }{\sqrt{K}}}_{\text{Sampling Error}}.
\end{align}
When the constant $C$ increases, the feasible set $\Fcal(P_\distr)$ gradually becomes a better approximation of the feasible set $\Fcal(P)$, as more feasible solutions in $\Fcal(P)$ are included in $\Fcal(P_\distr)$. The optimality gap $\Delta v(P_\distr)$, which can be viewed as the approximation error, is thus narrowed. On the other hand, as the set $\Fcal(P_\distr)$ expands, one needs more samples to ensure that the sampled feasible set $\Fcal(P_J)$ can approximate $\Fcal(P_\distr)$. In that sense, as we increase $C$, the second term of the right-hand side of \eqref{eq:bound_bias_complexity} also increases.

Since the constant $C$ can be arbitrary in Theorem~\ref{thm:main_largest_abs_of_all_dual_BFS}, we can in theory minimize the right-hand side of the inequality \eqref{eq:bound_bias_complexity} (or more precisely, inequality \eqref{eq:thm_convergence_with_abs_of_BFS}) to obtain a tighter bound, which results in the following corollary.

\begin{corollary}
	Define $P_\distr$, $P_J$, $\Ab_J$, $\gamma$, and $\| \Ab \|_{\max}$ as in Theorem~\ref{thm:main_largest_abs_of_all_dual_BFS}. Define a function $\Phi:(0,1) \times \Nbb \rightarrow \mathbb{R}$ as
	\begin{align*}
	\Phi(\delta,K) \equiv \inf_{C \geq 0 } \left\lbrace  \Delta  v(P_\distr) +  \frac{C \left( 1 + m \gamma \| \Ab \|_{\max} \right)}{\sqrt{K}} \left(  1 + \sqrt{2 \log \frac{2}{\delta}}  \right) \right\rbrace.
	\end{align*}
	For any $\delta \in (0,1)$, with probability at least $1- \delta$ over the sample $J$, the following statement holds: if $P_J$ is feasible and $\rank(\Ab_J) = m$, then $\Delta v(P_J) \leq \Phi(\delta,K)$.
	\end{corollary}

Unfortunately, this result is in general difficult to apply, because the function $\Phi$ is difficult to obtain in closed form. However, in Section~\ref{sec:P_distr}, we will later see how under certain conditions, a choice of $C$ that is small will be sufficient to ensure that $\Delta v(P_{\distr})$ is small, leading to an overall small bound on $\Delta v(P_J)$.

\subsubsection*{Computational Strengths and Weaknesses:} We compare the column randomization method to the CG method from a computational viewpoint. An obvious characteristic of the CG method is that it is a serial algorithm: to introduce a new column, one needs the dual solution of the restricted problem that consists of columns generated in previous iterations. This sequential nature unfortunately prevents the CG method from being parallelized. In contrast, the column randomization method is amenable to parallelization. Given a collection of processors, each processor can be used to sample a column and compute the constraint and objective coefficients in parallel, until $K$ columns in total are sampled across all processors. This can be especially advantageous in cases where the objective or constraint coefficients require significant effort compute, such as solving a dynamic program or integer program. For example, \cite{bertsimas2019airlift} considers a set partitioning model of a pickup and delivery problem arising in airlift operations, where each decision variable $x_{v,S}$ corresponds to an aircraft $v$ being assigned to a collection of shipments $S$ and the cost coefficient $c_{v,S}$ is the optimal value of a scheduling problem that determines the sequence of pickups and dropoffs of the shipments in $S$.

An obvious disadvantage of the column randomization method is that it does not guarantee optimality. Even if there exists an optimal solution of the complete problem $P$ that belongs to the feasible set $\Fcal(P_\distr)$ of problem $P_\distr$, the optimality gap still converges with rate $1 / \sqrt{K}$, which implies that the ``last-mile'' shrinkage of the optimality gap requires an increasing number of additional sampled columns. If optimality is a concern, instead of solely using the column randomization method, one could use it as a warm-start for the CG method. Specifically, let $J_{\text{nz}} = \{  j \mid  \tilde{x}^*_j > 0  \}$, where $\tilde{\xb}^*$ is the solution returned by Algorithm~\ref{alg:main}. Then, the set of variables $(x_j)_{J \in J_{\text{nz}}}$ and the columns $\Ab_{J_{\text{nz}}}$ can be used as the initial solution for the CG method. We test such a hybrid ``column-randomization-then-column-generation'' method in both of our numerical case studies (see Sections~\ref{subsec:numerics_CS_CRthenCG} and \ref{subsec:numerics_NCME_CRthenCG} in the ecompanion) and show that this hybrid method reaches provably optimal solutions in significantly less time than ordinary CG.

\subsubsection*{Additional Comments:} Considering space constraints, we relegate the following discussions to Section~\ref{sec:additional_discussion_on_thm1_thm2} of the e-companion: the lower bound on $v(P_J)$, the feasibility of $P_J$, interpretation of $\gamma$ and $\xi$, the comparison of Theorems~\ref{thm:main_largest_abs_of_all_dual_BFS} and \ref{thm:main_reduced_cost}, and the design of the randomization scheme $\rho$. Furthermore, in Section~\ref{sec:dependent_columns}, we explore the extension of our results to non-IID sampling of columns and sampling without replacement.

\section{Analysis of the distributional counterpart}
\label{sec:P_distr}

A key component of the theoretical guarantees presented in the previous section is the term $\Delta v(P_{\distr})$, which measures the gap between the distributional counterpart $P_{\distr}$ and the complete problem $P$. In this section, we provide more insight on the behavior of this term. In Section~\ref{subsec:P_distr_example}, we consider a toy example to provide some intuition for how this term behaves. Armed with this insight, in Section~\ref{subsec:P_distr_analysis} we formalize a result, Theorem~\ref{theorem:many_BFS_P_distr_bound}, which relates the gap $\Delta v(P_{\distr})$ to the abundance of nearly optimal, ``diverse'' basic feasible solutions of $P$. Lastly, in Section~\ref{subsec:P_distr_GM}, we consider three different random generative models for the complete problem $P$ and show that with high probability, a choice of $C$ that scales gracefully in $n$ (either $O(1)$ or $O(\log n)$) is sufficient to ensure that $\Delta v(P_{\distr})$ is small in terms of $n$ (either $O(1/\sqrt{n})$ or $O(\log n / \sqrt{n})$).

\subsection{A simple example}
\label{subsec:P_distr_example}

Consider the following full LP $P$ and its distributional counterpart:
\begin{align}
 P& : \ \min\{ \cb^T \xb \mid \oneb^T \xb = 1, \xb \geq \zerob\}, \label{prob:LP_easy} \\
P_{\distr}& : \ \min\{ \cb^T \xb \mid \oneb^T \xb = 1, \xb \leq C \xib, \xb \geq \zerob\}. \label{prob:LP_easy_distr}
\end{align}
Suppose also that we set the probability distribution $\xib$ to be the uniform distribution on $[n]$, that is, we set $\xi_n = 1/n$, so that all columns have the same probability to be chosen by the randomization scheme. Note that although the constraint coefficients and the right-hand side of the only equality constraint in $P$ are all equal to 1, it is possible to transform many LPs with a single constraint to this form. Specifically, consider the problem 
\begin{equation*}
P': \ \min \{ \cb'^T \xb' \mid \ab^T \xb' = b, \xb' \geq \zerob\}
\end{equation*}
where $\ab \geq \zerob$ and $b > 0$. We first divide both sides of the constraint $\ab^T \xb' = b$ by $b$; we then divide each variable's constraint coefficient and objective coefficient by $a_i / b$; and finally, we normalize the objective coefficient vector to have unit norm. This results in a problem of the form~\eqref{prob:LP_easy}, where each $c_i$ is defined as $c_i = v_i / \| \vb \|$ and $\vb = (c'_1 \cdot a_1 / b, \dots, c'_n \cdot a_n / b)$. The two problems $P'$ and $P$ are then equivalent, in the following way: $\xb = (x_1,\dots, x_n)$ is an optimal solution of $P$ if and only if $\xb' = ( (b/a_1) \cdot x_1, \dots, (b/a_n) \cdot x_n)$ is an optimal solution of $P'$. 

With this simple LP defined, we now wish to understand how we should set $C$ so that $\Delta v(P_{\distr}) = v(P_{\distr}) - v(P)$ is equal to zero. Let us consider two extreme cases.

\begin{enumerate}
\item When $\cb = -\eb_1$, where $\eb_1 = (1, 0, \dots, 0)$, then $v(P) = -1$. This optimal value can be only achieved by a single optimal solution, $\xb^* = \eb_1$, which is a basic feasible solution. Therefore, to ensure that $\Delta v(P_{distr}) = 0$, we must set $C = n$. Otherwise, if $C < n$, then the feasible region of $P_{\distr}$ will not contain $\xb^*$, and the optimality gap $\Delta v(P_{\distr})$ will not be zero. %

\item On the other hand, when $\cb = (-1 / \sqrt{n}) \oneb$, then we have $v(P) = - 1/ \sqrt{n}$. As in the previous example, $\xb^* = (1,0,0,\ldots,0)$ is an optimal basic feasible solution, and we can again set $C = n$ so that the set $\{ \xb \mid 0 \leq \xb \leq C \cdot \xib \}$ contains this optimal solution, resulting in $\Delta v(P_\distr) = 0$. However, upon closer inspection, one can see that there are actually multiple optimal non-basic solutions to the problem (in fact, every feasible solution is optimal). For example, another optimal solution is $\xb^* = (0.5, 0.5, 0, \dots, 0)$, which would imply that we can use the smaller value $C = 0.5 n$ to guarantee that $\Delta v(P_{\distr})$ to be zero. 

Note that although the non-sparse nature of this solution helps us, the infinity norm of the solution is also important. For example, if we consider the optimal solution $\xb^* = (0.7, 0.3, 0, \dots, 0)$, then we would need to set $C = 0.7 n$. In general, we need to set $C = n \| \xb^* \|_{\infty}$ to ensure that $\xb^* \in \{ \xb' \in \Rbb^n \mid \zerob \leq \xb' \leq C \cdot \xib \}$.

Following this logic, it turns out that the optimal solution with the lowest infinity norm is $\xb^* = (1/n,1/n,\ldots,1/n)$. For this solution, setting $C = 1$ is sufficient to ensure that the set $\{ \xb \mid \xb \leq C \cdot \xib \}$ can still include an optimal solution, resulting in $\Delta v(P_{\distr}) = 0$.
\end{enumerate}

By comparing these two cases, we can see that in the former case that we must have $C = O(n)$ to ensure $\Delta v(P_{\distr}) = 0$, whereas in the latter case, having $C = O(1)$ is sufficient to ensure $\Delta v(P_{\distr})$. What helps to ensure that $C$ can be small in the latter case is the existence of optimal solutions that are non-sparse and in particular, have low infinity norm. Thus, we should intuitively expect that a small gap $\Delta v(P_{\distr})$ can be achieved with a low value of $C$ when there exist optimal or nearly-optimal solutions with low infinity norm. In the following two sections, we build on this intuition to provide two different types of guarantees. First, in Section~\ref{subsec:P_distr_analysis}, we show that such nearly-optimal solutions with low infinity norm exist when there exist many near-optimal basic feasible solutions with low overlap in their bases, guaranteeing that $\Delta v(P_{\distr})$ will be small for a particular choice of $C$. Second, in Section~\ref{subsec:P_distr_GM}, we propose three different random generative models for the complete LP $P$, and show that with high probability, a small value of $C$ (either $O(1)$ or $O(\log n)$) is sufficient to ensure $\Delta v(P_{\distr})$ will be small (either $O(1/ \sqrt{n})$ or $O(\log n / \sqrt{n})$); in all three models, the key will be to show the existence of solutions to $P$ with small infinity norm.

\subsection{Bounding the distributional counterpart gap for a fixed $P$}
\label{subsec:P_distr_analysis}

Building on the intuition obtained in the prior section, we now present our first theoretical result on $\Delta v(P_{\distr})$, which depends on the structure of $P$ in terms of the behavior of nearly-optimal basic feasible solutions to $P$. 

\begin{theorem} \label{theorem:many_BFS_P_distr_bound}
Suppose that:
\begin{itemize}
\item $\xib$ is the uniform distribution on $[n]$, i.e., $\xi_j = 1 / n$ for all $j$; 
\item There exist $M$ basic feasible solutions, $\xb^{1}, \dots, \xb^M$, of $P$ that are within $\epsilon \geq 0$ of the optimal objective value, i.e., $\cb^T \xb^{i} - v(P) \leq \epsilon$;
\item Each variable $j \in [n]$ appears in at most $R$ of the corresponding bases $B^1, \dots, B^M$; and 
\item There exists a value $x_{\max}$ such that $\| \xb^i \|_{\infty} \leq x_{\max}$ for all $M$ BFSs.
\end{itemize}
 Then for $C = (nR / M) \cdot x_{\max}$, we have $\Delta v(P_{\distr}) \leq \epsilon$.
\end{theorem}

The proof of this result follows by showing that the average of the $M$ BFSs that are $\epsilon$-optimal is also an $\epsilon$-optimal feasible solution, and that the infinity norm of this solution is at most $R x_{\max} / M$. 

This result formalizes some of the insight from the stylized single-constraint example in Section~\ref{subsec:P_distr_example}. In particular, the choice of $C$ that achieves the gap of $\epsilon$ is $O(nR/ M)$. Thus, the larger the number $M$ of $\epsilon$-optimal BFSs, the smaller the sampling effort $K$ needs to be to guarantee that the gap of the column-randomized problem $\Delta v(P_J)$ will be within $\epsilon$. This, however, is modulated by $R$, which measures the diversity of the BFSs. When $R$ is small, it implies that the BFSs are different, in that the same column only appears in a small number of BFSs, and that the $M$ bases actually span a large set of columns. When this is the case, it makes sense that the sampling effort should be small, as there are many columns that we could use to form one of the BFSs or a convex combination of the BFSs. On the other hand, when $R$ is large, this implies that there are one or more columns that are common across many of the BFSs. In this case, it is reasonable that the sampling effort should be large, as we would need to sample these specific columns in order to be able to form one of the BFSs or a convex combination of them.

We note that a limitation of this result is that the existence of nearly optimal BFSs is taken as an assumption; our result does not provide conditions on $P$ which would ensure that this is the case. In general, it seems that for many types of large-scale LPs, it should be the case that there are many nearly-optimal BFSs. In Section~\ref{subsec:numerics_CS_small_instance} of the ecompanion, we show empirically that this is indeed the case for the cutting stock problem, which is a classical example of a large-scale LP that is usually solved via column generation. 
In particular, we demonstrate that $M$ can be much larger than $R$, indicating that the factor $nR/M$ in $C$ from Theorem~\ref{theorem:many_BFS_P_distr_bound} could exhibit sublinear growth in $n$. Additionally, note that the following holds:
\begin{align*}
Mm = \sum_{t=1}^M m = \sum_{t=1}^M \sum_{j=1}^n \mathbb{I}\{ j \in B^t \} = \sum_{j=1}^n \sum_{t=1}^M \mathbb{I}\{ j \in B^t \} \leq \sum_{j=1}^n R = nR,
\end{align*}
where the second equality follows the fact that each basis $B^t$ consists of $m$ columns and the inequality follows the definition of $R$. Along with the fact $R \leq M$, we have $m \leq nR/M \leq n$. Thus, the smallest possible value that $C$ can take in Theorem~\ref{theorem:many_BFS_P_distr_bound} is $mx_{\max}$, indicating that a $C$ that is constant or sublinear in $n$ would lead to a small $\Delta v(P_\distr)$ for large-scale LPs. On the other hand, the largest possible value of $C$ is $n x_{\max}$, in which case a large sampling effort $K$ would be required to achieve a small overall gap. We note that this latter case occurs when $R$ is close to or equal to $M$, which corresponds to the case where there is at least one column that appears in most or all near-optimal bases. While such a scenario can occur in some specially constructed LPs, we believe that this will not be the case for a large class of LPs that arise in practice. As noted above, this is not the case in the cutting stock problem (see the aforementioned Section~\ref{subsec:numerics_CS_small_instance} of the ecompanion), and is also generally not the case in the nonparametric choice model estimation problem (see our discussion on multiplicity of optimal solutions at the end of Section~\ref{sec:numerics_NCME}). Beyond these two crude bounds, it is challenging to give a more precise bound on $C$ without any further assumptions on the LP structure. For this reason, in Section~\ref{subsec:P_distr_GM} we will establish that, by assuming that the LP instances are randomly generated by three reasonable models, a carefully selected $C$ of $O(1)$ or $O(\log n)$ can indeed result in a small distributional counterpart gap of $O(1/\sqrt{n})$ or $O(\log n / \sqrt{n})$.

Lastly, we note here that the assumption of $\xib$ being the uniform distribution, i.e., $\xi_j = 1/n$ for all $j \in [n]$, is not a restrictive assumption in analyzing the distributional counterpart, since the result can be used to bound the performance of the column randomization method under a general class of distributions. We demonstrate this as follows. Let us first define $P_\distr^{\text{unif}}$ as
\begin{align*}
P_\distr^{\text{unif}}: \qquad \min \left\lbrace  \cb^T \xb \mid \Ab \xb = \bb, \quad \zerob \leq \xb \leq  (C/n) \cdot \oneb \right\rbrace.
\end{align*}
Now we consider a class of distributions with respect to a constant $\alpha > 0$:
\begin{align*}
\Xi_\alpha = \left\{  \xib \in \Rbb^n  \mid \oneb^T \xib = 1,\  \xi_j \geq \alpha / n, \quad \forall j \in [n] \right\rbrace.
\end{align*}
The set $\Xi_\alpha$ is a general class of distributions. Note that for every distribution $\xib$ in $\Xi_{\alpha}$, the probability $\xi_j$ is bounded away from zero, which means that every column has a positive probability of being sampled. Additionally, compared to the uniform distribution, it allows each column to be sampled with a different probability.

Suppose $P_J$ is a column-randomized LP that is obtained by sampling $K$ columns under a distribution $\xib$ from $\Xi_\alpha$. Following the notation in expression~\eqref{eq:bound_bias_complexity} and letting $C' = C / \alpha$, with probability at least $1-\delta$, we have the following statement: if $P_J$ is feasible and $\rank(\Ab_J) = m$,
\begin{align*}
v(P_J) & \leq  \min \left\lbrace \cb^T \xb \mid \Ab \xb = \bb, \zerob \leq \xb \leq C' \xib \right\rbrace + \frac{C' \cdot C_P \cdot C_\delta}{\sqrt{K}} \\
& \leq  \min \left\lbrace \cb^T \xb \mid \Ab \xb \leq \bb, \zerob \leq \xb \leq C' \cdot (\alpha / n) \cdot \oneb \right\rbrace + \frac{C' \cdot C_P \cdot C_\delta}{\sqrt{K}} \\
& =  v \left(P_\distr^{\text{unif}} \right) + \frac{1}{\alpha} \cdot \left(  \frac{C \cdot C_P \cdot C_\delta}{\sqrt{K}}  \right),
\end{align*}
which implies that the gap of $P_J$ satisfies $\Delta v (P_J) \leq \Delta v \left(P_\distr^{\text{unif}} \right)+ (1/\alpha) \cdot  \left(  {C \cdot C_P \cdot C_\delta}/{\sqrt{K}}  \right)$.

In other words, the second term $O(1/\sqrt{K})$ is scaled by a factor of $1/\alpha$. Therefore, the analysis of $P_{\distr}$ when $\xib$ is the uniform distribution can be used to provide performance guarantees for the column-randomized LP sampled by any $\xib \in \Xi_\alpha$. We will thus continue to use this uniform sampling assumption in our analyses in Section~\ref{subsec:P_distr_GM}.

\subsection{Analysis of the distributional counterpart under random generative models}
\label{subsec:P_distr_GM}

In this section, we present three different random generative models for large-scale linear programs, and investigate the behavior of the distributional counterpart gap $\Delta v(P_{\distr})$ under these three models. The idea is to assume that the complete LP $P$ is generated randomly according to a certain procedure, and to then develop a high probability bound for $\Delta v(P_{\distr})$ for a particular choice of the constant $C$. We will show that for all three models, the corresponding choice of $C$ is either constant (does not have an explicit dependence on $n$) or logarithmic in $n$, while $\Delta v(P_{\distr})$ is correspondingly either $O(1/ \sqrt{n})$ or $O(\log n / \sqrt{n})$ with high probability. 

In each generative model that we present, the final output is the triple $(\Ab, \bb, \cb)$ which fully defines the complete problem $P$. Some of the steps will involve randomly generating some of these objects, while others may involve choosing these objects in any arbitrary way that satisfies certain conditions; for any such case of the latter, the manner in which the object is chosen is not important, as the ensuing analysis of $\Delta v(P_{\distr})$ will not depend on how that object is chosen. It will, of course, depend on the probabilistic behavior of the objects chosen randomly.

\subsubsection{Generative model \GMDirichletNum}

The first generative model that we will consider is generative model \GMDirichletNum. This procedure is formalized as Algorithm~\ref{algorithm:GMDirichlet}. The idea in this procedure is that we start from some arbitrarily chosen set of columns $\Ab_1,\dots, \Ab_n$ and a scaling factor $\eta$. We then set the right-hand side vector $\bb$ as $\bb = \eta \sum_{j=1}^n \theta_j \Ab_j$, where $(\theta_1,\dots, \theta_n)$ is drawn uniformly from the $(n-1)$-dimensional unit simplex or equivalently, drawn from a $\Dirichlet(\alpha_1,\dots,\alpha_n)$ distribution where $\alpha_1 = \dots = \alpha_n = 1$. We then choose $\cb$ as any arbitrary unit norm vector that ensures that the optimal value $v(P)$ of the complete problem is nonnegative. 

\begin{algorithm}
\caption{Generative Model \GMDirichletNum}
\begin{algorithmic}[1]
	\STATE Fix any nonnegative constant $\eta \geq 0$.
	\STATE Fix any matrix $\Ab = [ \Ab_1 \cdots \Ab_n]$ of columns.
	\STATE Generate a random vector $\thetab = (\theta_1,\dots, \theta_n) \sim \Dirichlet(\alpha_1,\dots,\alpha_n)$, where $\alpha_1 = \dots = \alpha_n = 1$. 
	\STATE Set $\bb = \sum_{j=1}^n \eta \theta_j \Ab_j = \eta \Ab \thetab$. 
	\STATE Fix any $\cb \in \{ \vb \in \Rbb^n \mid \| \vb \|_2 = 1 \}$ such that $v(P) \geq 0$.
	\RETURN $(\Ab, \bb, \cb)$.
\end{algorithmic}
\label{algorithm:GMDirichlet}
\end{algorithm}

Before presenting our theoretical result on generative model \GMDirichletNum, we pause to make three important comments about the generative model. First, the $P$ generated in this way is always feasible by construction ($\xb = \eta \thetab$ is a feasible solution). Second, with regard to an interpretation of generative model \GMDirichletNum, note that the right-hand side vector $\bb$ is synthesized as a (scaled) convex combination of the columns $\Ab_1,\dots,\Ab_n$. %
Thus, we can think of the complete optimization problem $P$ as making the decision $\xb$ so that it has the same resource requirements as some reference or status quo decision given by $\eta \thetab$, that is, it satisfies $\Ab \xb = \Ab (\eta \thetab)$, while minimizing the objective function $\cb^T \xb$. By assuming that $\thetab \sim \Dirichlet(1,\dots,1)$, we are making the assumption that all decisions in the set $\{ \yb \in \Rbb^n \mid \sum_{j=1}^n y_j = \eta, \yb \geq \zerob \}$ are equally likely to be the status quo decision. For a particular choice of $\Ab$, we can also interpret the resulting $P$ as an estimation problem over the space of discrete probability distributions, which relates to one of the numerical experiments we consider (on nonparametric choice model estimation; see Section~\ref{sec:numerics_NCME}). Furthermore, it is worth noting that Step 5 in Algorithm~\ref{algorithm:GMDirichlet} is always achievable, and even when $P$ is such that $v(P) < 0$, it is possible to transform $P$ into an equivalent problem for which $v(P) \geq 0$; hence the requirement that $v(P)$ is nonnegative comes without loss of generality. Additional details regarding these two observations can be found in Section~\ref{subsec:additional_comment_on_generative_model_Dirichlet} of the e-companion.

Under generative model \GMDirichletNum, we have the following result, which bounds $\Delta v(P_{\distr})$ with high probability.

\begin{theorem} \label{theorem:GMDirichlet_gap_bound_whp}
Suppose that $P$ is generated according to generative model \GMDirichletNum. Assume that $\xib$ is the uniform distribution over $[n]$, that is, $\xi_j = 1 /n$ for all $j \in [n]$. Let $t \geq 1$.  Suppose that $C$ is set as %
\begin{equation*}
C = t \cdot \eta \cdot (\log n+1).
\end{equation*}
Then with probability at least $1 - 1/t$, $P_{\distr}$ is feasible and 
\begin{equation*}
\Delta v(P_{\distr}) \leq \frac{ t \eta (1 + \log n)}{\sqrt{n}}.
\end{equation*}
\end{theorem}
In words, Theorem~\ref{theorem:GMDirichlet_gap_bound_whp} states that for most problems $P$, setting $C$ to be logarithmic in $n$ is sufficient to ensure $\Delta v(P_{\distr}) = O( \log n / \sqrt{n})$. This result is particularly attractive because while $n$ may be unmanageably large, $\log n$ can be much smaller.

The proof of Theorem~\ref{theorem:GMDirichlet_gap_bound_whp} relies on an alternate characterization of the $\Dirichlet(1,\dots,1)$ distribution as the distribution of uniform spacings, and then using results on ordered uniform spacings, which are the order statistics of uniform spacings, to obtain a high probability bound on the random variable $\max_{j \in [n]} \theta_j$. This, together with a result that relates $\Delta v(P_\distr)$ to a bound $\beta$ on the minimum infinity norm of any feasible solution of $P$, yields the result. 

\subsubsection{Generative model \GMGaussianNum}

We now consider our second generative model. In generative model \GMGaussianNum, we assume that $\bb$ is chosen arbitrarily, and then we generate the columns of the $\Ab$ matrix. In particular, each of the $n$ columns, $\Ab_1, \dots, \Ab_n$, are drawn independently from a standard multivariate normal distribution. Upon selecting $\bb$ and drawing the columns $\Ab_1,\dots, \Ab_n$, we finally choose the objective coefficient vector $\cb$. We assume that $\cb$ is chosen as any vector with unit norm and that ensures that $v(P) \geq 0$, if $P$ is feasible; if $P$ is not feasible, we then simply select any $\cb$ with unit norm. This procedure is formalized below as Algorithm~\ref{algorithm:GMGaussian}. 

\begin{algorithm}
\caption{Generative Model \GMGaussianNum}
\begin{algorithmic}[1]
	\STATE Fix any $\bb \in \Rbb^m$.
	\STATE Generate $n$ i.i.d. random vectors $\Ab_1,\dots,\Ab_n \sim \Normal( \zerob, \Ib)$, where $\Normal(\zerob, \Ib)$ denotes a standard multivariate normal distribution with $\zerob \in \Rbb^m$ as the mean vector and the $m$-by-$m$ identity matrix $\Ib$ as the covariance matrix.
	\STATE Set $\Ab = [\Ab_1 \ \cdots \ \Ab_n]$. 
	\STATE If $\{\xb \mid \Ab \xb = \bb, \xb \geq \zerob\}$ is non-empty, fix any $\cb \in \{ \vb \in \Rbb^n \mid \| \vb \|_2 = 1 \}$ such that $v(P) \geq 0$; otherwise, fix any $\cb \in \{ \vb \in \Rbb^n \mid \| \vb \|_2 = 1 \}$. 
	\RETURN $(\Ab, \bb, \cb)$.
\end{algorithmic}
\label{algorithm:GMGaussian}
\end{algorithm}

For this generative model, we have the following guarantee on the distributional counterpart gap. Note that unlike generative model \GMDirichletNum, the problem $P$ generated by generative model \GMGaussianNum need not be feasible. However, this guarantee also ensures that both $P$ and $P_{\distr}$ are feasible.

\begin{theorem} \label{theorem:GMGaussian_gap_bound_whp}
Suppose that $P$ is generated according to generative model \GMGaussianNum. Assume that $\xib$ is the uniform distribution over $[n]$, that is, $\xi_j = 1 /n$ for all $j \in [n]$. Let $t \geq 1$, and suppose that $n > 4 \pi t^2 m$. Suppose that $C$ is set as 
\begin{equation*}
C = \| \bb \|_2  \cdot \frac{1}{ \frac{\sqrt{2}}{2\sqrt{\pi}} - \frac{ t \sqrt{2m}}{\sqrt{n}}}.
\end{equation*}
Then, with probability at least $1 - 1/t$, both $P$ and $P_{\distr}$ are feasible, and
\begin{equation*}
\Delta v(P_{\distr}) \leq \frac{ \| \bb \|_2 }{ \sqrt{n} } \cdot \frac{1}{ \frac{\sqrt{2}}{2\sqrt{\pi}} - \frac{ t \sqrt{2m} }{\sqrt{n}}}.
\end{equation*}
\end{theorem}
The proof of Theorem~\ref{theorem:GMGaussian_gap_bound_whp} follows by bounding the minimum infinity norm problem $\min\{ \| \xb \|_{\infty} \mid \Ab \xb = \bb, \xb \geq \zerob \}$. To obtain a bound, it turns out that the dual problem can be written as a maximization problem over a variable $\pb$ subject to a constraint that can be written as a sample average of a certain function of the columns that depends on the dual variable $\pb$. The true expectation of this function of a random column at a given $\pb$ can be found in closed form and turns out to be $\| \pb \|_2$ multiplied by a scaling constant, which gives rise to a dual problem that is essentially the optimization of a linear function $\pb^T \bb$ subject to a constraint that looks like $C' \| \pb \|_2 \leq 1/n $, where $C'$ is a constant. This is what gives rise to the $\| \bb \|_2$ part of the definition of $C$, and the $\| \bb \|_2 / \sqrt{n}$ part of the bound on $\Delta v(P_{\distr})$. To ensure that the sample average is close to this true expectation, we consider the Rademacher complexity of a certain function class, which allows us to bound with high probability the difference between the aforementioned sample average and its true expectation using a term of the form $C'' \| \pb \|_2$ where $C''$ is a constant that depends on $\sqrt{n}$. The two constants $C'$ and $C''$ are what gives rise to the factor $[ \sqrt{2} / (2\sqrt{\pi}) - t \sqrt{2m} / \sqrt{n}]^{-1}$ in the bound.

As with our previous generative model result, Theorem~\ref{theorem:GMGaussian_gap_bound_whp} states that when the columns of $\Ab$ are drawn i.i.d. from a standard multivariate normal distribution, then a choice of $C$ that is $O(1)$ will result in a distributional counterpart gap that is $O(1/\sqrt{n})$. Since $C$ directly translates into the sampling effort $K$, this implies that a constant sampling effort should be sufficient to ensure that the gap of the column-randomized LP $\Delta v(P_J)$ is $O(1/\sqrt{n})$. %

With regard to the generality of generative model \GMGaussianNum, we make the following two remarks. First, generative model \GMGaussianNum is quite general, in the following informal sense. For any linear program $P \equiv \min \{ \cb^T \xb \mid \Ab \xb = \bb, \xb \geq \zerob\}$ where the columns $\Ab_1,\dots, \Ab_n$ have arbitrary non-zero Euclidean norm, we can transform the problem into an equivalent problem with columns that are normalized to have unit norm. In particular, let $\Ab' = [\Ab'_1 \ \cdots \ \Ab'_n]$, where $\Ab'_j = \Ab_j / \| \Ab_j \|_2$, define $\cb' = (c'_1, \dots, c'_n)$ as $c'_j = c_j / \| \Ab_j \|_2$, and define the new problem $P'$ as $P' \equiv \min \{ \cb'^T \xb' \mid \Ab' \xb' = \bb, \xb' \geq \zerob \}$. Observe that for every feasible solution $\xb$ of $P$, the solution $\xb' = (\| \Ab_1 \|_2 \cdot x_1, \dots, \| \Ab_n \|_2 \cdot x_n)$ is a feasible solution whose objective in $P'$ is the same as the objective of $\xb$ in $P$, so by solving $P'$ we can solve $P$, and vice versa. Now, recall that when $\Ab_j$ follows the standard multivariate normal distribution, then $\Ab_j / \| \Ab_j \|_2$ is uniformly distributed on the $m$-dimensional unit sphere $S^m = \{ \vb \in \Rbb^m \mid \| \vb \|_2 = 1\}$. Thus, by assuming that the columns of $\Ab$ are drawn from the standard multivariate normal distribution, we ensure that each column of the transformed matrix $\Ab'$ is uniformly distributed on the unit sphere $S^m$. 

Second, building on the intuition in the prior remark, we can consider a modification of generative model \GMGaussianNum, where instead of sampling columns $\Ab_1,\dots,\Ab_n$ independently from a standard multivariate normal distribution, we sample them from the uniform distribution on $S^m$. The resulting generative model is described in Section~\ref{subsec:proofs_GMUniformSphere} of the ecompanion. The same proof machinery used for Theorem~\ref{theorem:GMGaussian_gap_bound_whp}, with a few careful (albeit tedious) modifications, goes through for this new generative model, resulting in a similar guarantee for this new model (Theorem~\ref{theorem:GMUniformSphere_gap_bound_whp} in Section~\ref{subsec:proofs_GMUniformSphere}). The main difference in the new guarantee is that the factor $[ \sqrt{2} / (2\sqrt{\pi}) - t \sqrt{2m} / \sqrt{n})]^{-1}$ that appears in Theorem~\ref{theorem:GMGaussian_gap_bound_whp} is replaced by the factor $[ \sqrt{2} / (2\sqrt{\pi} \mu_m) - t \sqrt{2} / \sqrt{n})]^{-1}$, where $\mu_m = \sqrt{2} \Gamma((m+1)/2)/ \Gamma(m/2)$ is the mean of a chi distributed random variable with $m$ degrees of freedom. 

\subsubsection{Generative model \GMBernoulliCoveringNum}

Lastly, we turn our attention to our final generative model, generative model~\GMBernoulliCoveringNum. For this generative model, we deviate slightly from the previous two models by considering a covering LP, as opposed to a standard form LP. The covering LP is defined as
\begin{equation*}
P^{\covering}: \qquad \min \{ \cb^T \xb \mid \Ab \xb \geq \bb, \xb \geq \zerob\},
\end{equation*}
where each entry of $\Ab$ is nonnegative, $\bb$ is a nonnegative $m$-dimensional vector and $\cb$ is a nonnegative $n$-dimensional vector, which we again assume to be normalized to have unit norm, i.e., $\| \cb \|_2 = 1$. The distributional counterpart of this problem, $P^{\covering}_{\distr}$ is defined as
\begin{equation*}
P^{\covering}_{\distr}: \qquad \min \{ \cb^T \xb \mid \Ab \xb \geq \bb, \xb \leq C \xib, \xb \geq \zerob \},
\end{equation*}
and the distributional counterpart gap can be defined as $\Delta v(P^{\covering}_{\distr}) = v(P^{\covering}_{\distr}) - v(P^{\covering})$. It can be shown that a modified version of Theorem~\ref{thm:main_largest_abs_of_all_dual_BFS}, which bounds the gap of the column-randomized LP in terms of the gap of the distributional counterpart plus a $O(1/\sqrt{K})$ term that depends on the maximum infinity norm of any dual basic solution, holds for $P^{\covering}$ (see Section~\ref{subsec:special_structures_covering}). 

Our final generative model, generative model~\GMBernoulliCoveringNum, is defined below as Algorithm~\ref{algorithm:GMBernoulliCovering}. In this model, each entry of $A_{i,j}$ is 0 or 1, generated as an independent Bernoulli random variable with a row-dependent probability $q_i$. 

\begin{algorithm}
\caption{Generative Model \GMBernoulliCoveringNum}
\begin{algorithmic}[1]
	\STATE Fix $m$ probabilities, $q_1,\dots, q_m \in (0,1)$. 
	\STATE Fix any $\bb \in \Rbb^m$ such that $\bb \geq \zerob$. 
	\STATE Generate $n$ i.i.d. random vectors $\Ab_1,\dots,\Ab_n$, where each $A_{i,j} \sim \Bernoulli(q_i)$.
	\STATE Set $\Ab = [\Ab_1 \ \cdots \ \Ab_n]$. 
	\STATE If $\{\xb \mid \Ab \xb \geq \bb, \xb \geq \zerob\}$ is non-empty, fix any $\cb \in \{ \vb \in \Rbb^n \mid \| \vb \|_2 = 1 \}$ such that $v(P) \geq 0$; otherwise, fix any $\cb \in \{ \vb \in \Rbb^n \mid \| \vb \|_2 = 1 \}$. 
	\RETURN $(\Ab, \bb, \cb)$.
\end{algorithmic}
\label{algorithm:GMBernoulliCovering}
\end{algorithm}

For this generative model, we have the following result which bounds $\Delta v(P^{\covering}_{\distr})$ with high probability. 

\begin{theorem} \label{theorem:GMBernoulliCovering_gap_bound_whp}
Suppose that $P$ is generated according to generative model \GMBernoulliCoveringNum. Assume that $\xib$ is the uniform distribution over $[n]$, that is, $\xi_j = 1 /n$ for all $j \in [n]$. Let $\delta \in (0,1)$ and that $n > \log(m/\delta) / [ 2 (\min_{i \in [m]} q_i)^2 ]$. Suppose that $C$ is set as 
\begin{equation*}
C = \max_{i \in [m]} b_i  \cdot \frac{1}{ \min_{i' \in [m]} q_{i'} - \sqrt{ \frac{1}{2n} \log \frac{m}{\delta}} }
\end{equation*}
then with probability at least $1 - \delta$, both $P$ and $P_{\distr}$ are feasible and 
\begin{equation*}
\Delta v(P^{\covering}_{\distr}) \leq \frac{ \max_{i \in [m]} b_i }{\sqrt{n}} \cdot \frac{1}{ \min_{i' \in [m]} q_{i'} - \sqrt{ \frac{1}{2n} \log \frac{m}{\delta}} }.
\end{equation*}
\end{theorem}
Similarly to our previous results, we establish this result by bounding (the dual of) the infinity norm problem $\min\{ \| \xb \|_{\infty} \mid \Ab \xb \geq \bb, \xb \geq \zerob \}$. Although the dual is generally challenging to analyze due to the presence of the $(\cdot)_+ = \max\{0, \cdot\}$ function, we can leverage the fact that the dual variable is nonnegative (due to the primal covering constraint $\Ab \xb \geq \bb$, as opposed to the equality constraint $\Ab \xb = \bb$ in the standard form LP) and that the columns of $\Ab$ are nonnegative. This allows us to bound the dual optimal objective in closed form in terms of the row sums of $\Ab$. An application of Hoeffding's inequality and the union bound allows us to then bound the deviation of the row sums of $\Ab$ with high probability, leading to the above result. 

An attractive aspect of generative model~\GMBernoulliCoveringNum is that the constraint matrix that one obtains is a sparse 0-1 matrix; in contrast, under generative model \GMGaussianNum, the matrix $\Ab$ is almost surely not sparse. A limitation of generative model~\GMBernoulliCoveringNum and Theorem~\ref{theorem:GMBernoulliCovering_gap_bound_whp} is that it only applies to the covering LP $P^{\covering}$, as opposed to the general standard form LP $P$. Unfortunately, it seems difficult to generalize the proof approach to the case where $\Ab$ is used in the standard form LP $P$, as the dual variable of the constraint $\Ab \xb = \bb$ in the infinity norm problem $\min\{ \| \xb \|_{\infty} \mid \Ab \xb = \bb, \xb \geq \zerob \}$ is no longer forced to be nonnegative.

\subsubsection{Concluding remarks}

Overall, the main takeaway from this section is that under three different and reasonably broad generative models for the complete LP $P$, there exists a choice of $C$ that (1) scales gracefully in $n$ and (2) ensures that the distributional counterpart gap, $\Delta v(P_{\distr})$, will be small in terms of $n$ with high probability. A limitation of these results is that many type of LPs in practice are highly structured, and may not look like LPs that would be produced by our generative models. For example, as discussed above, under generative model \GMGaussianNum the matrix $\Ab$ is almost surely not sparse, whereas this is the case for many LPs with special structure, such as network flow problems. Similarly, under generative model \GMDirichletNum, the right hand side vector $\bb$ arises as a scaled random convex combination of the columns, but of course $\bb$ could be chosen in a way that is unlikely to come about from such a combination (e.g., a scaled convex combination of a small set of columns). Nevertheless, we believe that our results are useful in providing intuition for how $\Delta v(P_{\distr})$ will behave under three stylized models for how the complete LP is formed.

\section{Numerical experiments with the cutting stock problem}
\label{sec:numerics_CS}

In this section, we apply the column randomization method to the cutting-stock problem, a well-known large-scale linear program that is commonly solved by CG. We follow the notation in \cite{bertsimas1997introduction} and briefly review the problem for completeness.

A paper company needs to satisfy a demand of $b_i$ rolls of paper of width $w_i$, for each $i \in [m]$. The company has supply of large rolls of paper of width $W$ such that $W \geq w_i$ for $i \in [m]$. %
To meet the demand, the company slices the large rolls into smaller rolls according to \emph{patterns}. A pattern is a vector of nonnegative integers $ (a_1,a_2,\ldots,a_m)$ that satisfies $\sum_{i=1}^m a_i w_i \leq W$, where each $a_i$ is the number of rolls of width $w_i$ to cut from the large roll. Let $n$ be the number of all feasible patterns and let $(a_{1j},a_{2j},\ldots,a_{mj})$ be the $j$th pattern for $j \in [n]$. Let $\Ab$ be the matrix such that $A_{ij} = a_{ij}$ for $i \in [m]$ and $j \in [n]$. The \emph{cutting-stock problem} is to minimize the number of large rolls of papers used while satisfying the demand, which can be formulated as the following covering LP:
\begin{align}
	\label{problem:LP_cutting_stock}
	P^{\text{CS}}: \quad \underset{\xb \in \Rbb^n, \xb \geq \zerob}{\text{minimize}}  \left\lbrace \,\, \sum_{j=1}^n x_j \,\, \bigg| \,\, \sum_{j=1}^n a_{ij} x_j \geq b_i ,  \,\, \forall i \in [m] \right\rbrace.
\end{align}
Explicitly representing the constraint matrix $\Ab$ in full is usually impossible: the number of feasible patterns $n$ can be huge even if the number of demanded widths $m$ is small. A typical solution method is column generation, in which each iteration proceeds as follows. Given a set of patterns $J = \{j_1,j_2,\ldots,j_K\}$, solve the restricted problem $P^{\text{CS}}(J): \quad \underset{\tilde{\xb} \in \Rbb^{K}}{\text{minimize}}  \left\lbrace \sum_{k=1}^K \tilde{x}_{k}  \mid \sum_{k=1}^K \Ab_{j_k} \tilde{x}_k \geq \bb , \tilde{\xb} \geq \zerob \right\rbrace$ and let $\pb$ be the optimal dual solution. Then find a new pattern $j_{K+1}$ such that the corresponding new column has the most negative reduced cost $1 - \pb^T \Ab_{j_{K+1}}$. If the reduced cost is nonnegative, the current solution is optimal and the procedure terminates; otherwise, we add $j_{K+1}$ to the collection $J$ and repeat the procedure. The problem of finding the column with the most negative reduced cost is equivalent to solving the following subproblem:
\begin{align}
	\label{problem:sub_cutting_stock}
	P^{\text{CS-sub}}: \quad \underset{\ab \in \Nbb_+^m}{\text{maximize}}  \left\lbrace \,\, \sum_{i=1}^m p^*_i a_i \,\, \bigg| \,\, \sum_{i=1}^m w_i a_i \leq W \right\rbrace,
	\end{align}
where $\Nbb_+$ is the set of nonnegative integers; if the optimal value $v(P^{\text{CS-sub}})$ is smaller than $1$, then we terminate the column generation procedure; otherwise, we let pattern $j_{K+1}$ correspond to the optimal solution of $P^{\text{CS-sub}}$ and add it to $J$.

Instead of column generation, we can consider solving the cutting-stock problem by the column randomization method. In our implementation of the column randomization method, we consider the randomization scheme described in Algorithm~\ref{alg:sampling_cutting_stock}. The randomization scheme essentially starts with an empty pattern, i.e., $(a_1,\dots, a_m) = (0, \dots, 0)$ and at each iteration, it increments $a_i$ for a randomly chosen $i$, while ensuring that it does not exceed the available width $W$. We refer to this randomization scheme as the incremental randomization scheme and denote it by $\rho_{\Incremental}$. We note that Algorithm~\ref{alg:sampling_cutting_stock} is not the only way to sample columns of $\Ab$, and one can consider other randomization schemes that would lead to potentially better performance of the column randomization method. In Sections~\ref{subsec:numerics_CS_uniform_vs_incremental} and \ref{subsec:numerics_CS_biased_vs_incremental}, we will see two other randomization schemes for this problem.

\begin{algorithm}
\SingleSpacedXI
	\caption{Incremental randomization scheme $\rho_{\Incremental}$ for the cutting-stock problem.}
	\label{alg:sampling_cutting_stock}
	\begin{algorithmic}[1]
		\STATE Column $\ab$ is a zero vector of length $m$ and $\zeta \leftarrow W$.
		\WHILE {$\zeta > 0$}
			\STATE $I \leftarrow \{ i \mid w_i \leq \zeta \} $.
			\IF {$| I | \geq 1$}
				\STATE Sample an index $i$ uniformly at random from $I$.
				\STATE Update $a_i \leftarrow a_i + 1$ and $\zeta \leftarrow \zeta - w_i$.
			\ELSE
				\STATE Break the while loop
			\ENDIF
		\ENDWHILE
		\RETURN Column $\ab$.
	\end{algorithmic}
\end{algorithm}

\begin{figure}[h!]
	\centering
	\includegraphics{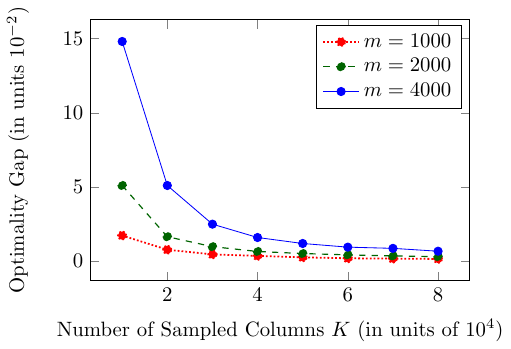}
	\caption{Performance of the column randomization method on the cutting-stock problem with respect to number of columns $K$ and number of required widths $m$.}
	\label{fig:cutting_stock_varying_K_and_m}
\end{figure}

In Figure~\ref{fig:cutting_stock_varying_K_and_m}, we illustrate the performance of column-randomized linear programs for the cutting-stock problem with respect to number of columns $K \in \{ 2 \times 10^4, 4 \times 10^4, 6 \times 10^4, 8 \times 10^4 \}$ and number of required widths $m \in \{ 1000,2000,4000\}$. We note that the value of $m$ significantly affects size and complexity of the problem: as $m$ increases, there are more possible patterns and thus $n$ increases as well. For the CG approach, $m$ defines the number of integer variables in the subproblem~\eqref{problem:sub_cutting_stock}; as it increases, the subproblem becomes more challenging. We set $W = 10^5$; we draw each $w_i$ uniformly at random from $\{W/10, W/10 + 1, \ldots, W/4-1, W/4\}$ without replacement; and we draw each $b_i$ independently uniformly at random from $\{1,\ldots,100\}$. %
We measure the performance of column-randomized linear programs $P^{\text{CS}}_J$, where each column is obtained by Algorithm~\ref{alg:sampling_cutting_stock}, by its relative optimality gap ${\Delta v(P^{\text{CS}}_J)}/{ v(P^{\text{CS}})}$. For each value of $m$ and $K$, we run the column-randomized method 20 times and compute the average optimality gap, which is plotted in Figure~\ref{fig:cutting_stock_varying_K_and_m}. Before continuing, we note here that there are many ways to randomly generate cutting-stock instances. Our goal is not to exhaustively evaluate the numerical performance of the column randomization method on every possible family of instances, but rather to understand its performance on a reasonably general set of instances.

We first observe that the curves in Figure~\ref{fig:cutting_stock_varying_K_and_m} approximately match the convergence rate of $1/\sqrt{K}$ in Theorems~\ref{thm:main_largest_abs_of_all_dual_BFS} and \ref{thm:main_reduced_cost}. In addition, the speed of convergence significantly slows down after the optimality is smaller than $2\%$; see the curve for $m=1000$. %
Second, as the problem size increases, we need more samples to return comparable performance in terms of optimality gap. This is reflected by the fact that for a fixed number of columns $K$, the optimality gap is larger for larger $m$.

We further compare the runtime of the column randomization method to that of the CG method in Table~\ref{table:cutting_stock_CR_vs_CG}. The first column of the table indicates the value of $m$, which quantifies the problem size and subproblem complexity. The second column indicates the number of sampled columns $K$ in the column-randomized linear program. The third and fourth columns indicate relative optimality gap ${\Delta v(P^{\text{CS}}(J))}/{ v(P^{\text{CS}})}$ and runtime of the column randomization method, respectively; for both of these metrics, we report the average over 20 runs of the column-randomized method. The fifth column shows the time required by the CG method to reach the same (average) relative optimality gap. We also list the total duration for CG (i.e., the time required for CG to reach a 0\% optimality gap) in the fifth column, and denote it by ``(total)''. %

Table~\ref{table:cutting_stock_CR_vs_CG} shows that, when the problem is small ($m=1000$), the column randomization method returns a high-quality solution with an optimality gap below $1\%$, within $30$ seconds and with $2 \times 10^4$ sampled columns. Doubling or tripling the number of sampled columns does not significantly improve the performance, as the optimality gap is already small. Meanwhile, CG also works well when $m=1000$, obtaining the optimal solution in a reasonable time (within fifteen minutes). On the other hand, when the problem is large ($m=4000$), the runtime of CG dramatically increases, as it needs almost 5000 seconds (just under 1.5 hours) to reach a $5\%$ optimality gap. The computational limiting factor comes from solving the subproblem, which becomes more difficult as $m$ increases. On the other hand, the column randomization method only needs ten minutes to reach a $1 \%$ optimality gap. This demonstrates the value of solving linear programs by the column randomization method in lieu of CG when the subproblem is intractable.

\begin{table}[]
\SingleSpacedXI
	\centering
	\small
	\begin{tabular}{crrrr}
		\toprule
		\multicolumn{1}{c}{{Demand Types $(m)$}}    &  \multicolumn{1}{c}{{Columns $(K)$}}     & \multicolumn{1}{c}{{Optimality Gap $(\%)$}}          & \multicolumn{1}{c}{{Runtime (s)}}     & \multicolumn{1}{c}{{CG Runtime (s)}}    \\
		\midrule
		$1000$ & $2 \times 10^4$  & 0.78 & 28.4   & 365.5    \\
		 & $4\times 10^4$  & 0.36 & 56.4   & 411.7   \\
		 & $6\times 10^4$  & 0.20 & 89.3  & 456.4    \\
		 & $8\times 10^4$ & 0.16  & 122.5  & 475.1    \\ & & & & (total) 775.4 
		 \\ \midrule
		2000 & $2\times 10^4$  & 1.65 & 58.9   & 1330.6   \\
		 & $4\times 10^4$  & 0.65 & 120.1  & 1622.8   \\
		 & $6\times 10^4$  & 0.43 & 197.9  & 1732.2   \\
		 & $8\times 10^4$ & 0.31  & 287.6  & 1805.0
		 \\ & & & & (total) 2932.92    \\ \midrule
		4000 & $2\times 10^4$  & 5.10 & 139.4  & 4979.8   \\
		 & $4\times 10^4$  & 1.59  & 314.2  & 7175.2   \\
		 & $6\times 10^4$  & 0.95 & 527.1  & 7670.1  \\
		 & $8\times 10^4$ & 0.68 & 768.6 & 7940.0
		 \\ & & & & (total) 13336.1
		 \\ \bottomrule
	\end{tabular}
	\caption{Performance of the column randomization method on the cutting stock problem for different problem sizes and numbers of sampled columns. 	\label{table:cutting_stock_CR_vs_CG}}
\end{table}

Finally, we have included additional numerical experiments in Section~\ref{sec:cutting_stock_continued} of the e-companion. Specifically, we compare the column randomization method under different sampling schemes in Sections~\ref{subsec:numerics_CS_uniform_vs_incremental} and \ref{subsec:numerics_CS_biased_vs_incremental}. Furthermore, we investigate the combined effectiveness of the column randomization method with CG in Section~\ref{subsec:numerics_CS_CRthenCG}. In Section~\ref{subsec:numerics_CS_small_instance}, we offer further insights into the strong performance of the column randomization method in the cutting stock problem, drawing connections to Theorem~\ref{theorem:many_BFS_P_distr_bound}.

\section{Numerical experiments with nonparametric choice model estimation}

\label{sec:numerics_NCME}

The second problem we consider is nonparametric choice model estimation, which is a modern application of large-scale linear programming and CG. In particular, we consider estimating the ranking-based choice model from data \citep{farias2013nonparametric,van2015market,mivsic2016data}. For completeness, we first briefly review the problem.

In the ranking-based nonparametric choice model, we assume that a retailer offers $N$ different products, indexed from $1$ to $N$. We use the index $0$ to represent the no-purchase alternative, which is always available to customer. Together, we refer to the set $[N]^+ \equiv\{ 0,1,\ldots,N \}$ as the set of purchase options. %
A ranking-based choice model $(\Sigmab,\lambdab)$ consists of two components. The first component $\Sigmab$ is a collection of rankings over options $[N]^+$, in which each ranking represents a customer type. We use $\sigma(i)$ to indicate the rank of option $i$, where $\sigma(i) < \sigma(j)$ implies that $i$ is more preferred to $j$ under the ranking $\sigma$.  When a set of products $S \subseteq [N]$ is offered, a customer of type $\sigma$ selects option $i$ from the set $S \cup \{ 0 \}$ with the lowest rank, i.e., the option $\arg\min_{i  \in S \cup \{ 0 \}} \sigma(i)$. The second component $\lambdab$ is a probability distribution over rankings in the set $\Sigmab$; the element $\lambda_\sigma$ can be interpreted as the probability that a random customer would make decisions according to ranking $\sigma$.

To estimate a ranking-based model, we utilize data in the form of past sales rate information. Here we consider the type of data described in \cite{farias2013nonparametric}; we refer readers to that paper for more details. Assume that the retailer has provided $M$ assortments $\Scal = \{S_1,S_2,\ldots,S_M \}$ in the past, where each $S_m \subseteq [N]$. For each assortment $S_m$, the retailer observes the choice probability $v_{i,m}$ for assortment $S_m$ and option $i$, which is the fraction of past transactions in which a customer chose $i$, given that assortment $S_m$ was offered. We let $v_{(i,m)} \equiv 0$ if $i \notin S \cup \{ 0 \}$. %

The estimation of a ranking-based choice model $(\Sigmab,\lambdab)$ can be formulated in the form of problem $P^\port$
(Section~\ref{subsec:special_structures_portfolio}). We first notice that there are in total $(N+1)!$ rankings over $[N]^+$, which we enumerate as $\sigma_1,\sigma_2,\ldots,\sigma_{(N+1)!}$. We let the $k$th column of the problem correspond to ranking $\sigma_k$, for $k \in [ (N+1)! ]$. We use $\alpha_{(i,m),k}$ to indicate whether a customer following ranking $\sigma_k$ would choose option $k$ when offered assortment $S_m$. The estimation problem can then be written as
\begin{subequations}
	\begin{alignat}{2}
	P^{\text{EST}}: \quad\underset{\lambdab \zerob, \hat{\vb}}{\text{minimize}}  \quad & \Dcal(\hat{\vb},\vb) \\   
	\text{such that}  \quad & \sum_{k=1}^{ (N+1)!} \alpha_{(i,m),k} \cdot \lambda_k = \hat{v}_{(i,m)}, \quad \forall m \in [M], i \in [N]^+, \label{problem:ncm_est_Alambdaeqv}\\ 
	& \sum_{k=1}^{(N+1)!} \lambda_k = 1, \label{problem:ncm_est_unitsum}\\ 
	& \lambdab \geq \zerob, \label{problem:ncm_est_nonnegative}
	\end{alignat}
\end{subequations}
where $\hat{\vb}$ and $\vb$ are vectors of $\hat{v}_{(i,m)}$ and $v_{(i,m)}$ values, respectively, for $ i \in [N]^+$ and $m \in [M]$. The function $\Dcal$ measures the error between the predicted choice probabilities $\hat{\vb}$ and the actual choice probabilities $\vb$. We follow \cite{mivsic2016data} and set $\Dcal = \| \hat{\vb} - \vb \|_1$, which has Lipschitz constant $\sqrt{M(N+1)}$. %

We notice that even if $N$ is merely $10$, problem $P^{\text{EST}}$ has nearly $4 \times 10^7$ columns. Given that problem $P^{\text{EST}}$ may have an intractable number of columns, \cite{van2015market} and \cite{mivsic2016data} applied CG to solve the problem. 
Alternatively, we can apply the column randomization method. We consider the randomization scheme described in Algorithm~\ref{alg:sampling_ranking}, where we first randomly generate a ranking (line 2) and then map its decision under each assortment to form a column (lines 3-5). We refer to this randomization scheme as the \emph{uniform randomization scheme} and denote it by $\rho_{\Uniform}$. 

Before continuing, we pause to make three important remarks. First, we note that sampling a ranking uniformly at random (line 2) requires minimal computational effort, and can be done by a single function call in most programming languages. %
Second, we also note that while in Algorithm~\ref{alg:sampling_cutting_stock} we directly sample the coefficients of a column, in Algorithm~\ref{alg:sampling_ranking} we instead first sample the underlying ``structure'' of the column (a ranking) then obtain the corresponding coefficients; this illustrates the problem-specific nature of the randomization scheme. Lastly, we note that the paper of \cite{farias2013nonparametric} considered a linear program for computing the worst-case revenue of an assortment, which is effectively the minimization of a linear function of $\lambdab$ subject to constraints~\eqref{problem:ncm_est_Alambdaeqv}--\eqref{problem:ncm_est_nonnegative}. The paper considered a solution method for this problem based on sampling constraints in the dual (which is equivalent to sampling columns in the primal), but did not compare this approach to column generation, which will do shortly.  %

\begin{algorithm}
	\SingleSpacedXI
	\caption{Uniform randomization scheme $\rho_{\Uniform}$ for the nonparametric choice estimation problem.}
	\label{alg:sampling_ranking}
	\begin{algorithmic}[1]
		\STATE Initialize $\alpha_{(i,m)} \leftarrow 0$ for $i \in [N]^+$ and $m \in [M]$.
		\STATE Sample a ranking/permutation $\sigma: [N]^+ \to [N]^+$ uniformly at random.
		\FOR {$m \in [M]$}
		\STATE $i^* \leftarrow \arg \min_{ i \in S_m \cup \{ 0 \} } \sigma(i)$.
		\STATE $\alpha_{(i^*,m)} \leftarrow 1$
		\ENDFOR
		\RETURN Column $\alphab = (\alpha_{(i,m)})_{i \in [N]^+, m \in [M]}$.
	\end{algorithmic}
\end{algorithm}

We compare the performance of the column randomization method to that of CG with the following experiment setup. We assume that customers follow multinomial logit (MNL) model to make decision, that is, the choice probability $v_{i,m}$ follows $v_{i,m} = \exp(u_i) / \left( 1 + \sum_{j \in S_m} \exp(u_j)  \right)$ for a given assortment $S_m$, where each parameter $u_i$ represents the expected utility of product $i$. We choose each $u_i \sim U[0,1]$, i.e., uniformly at random from interval $[0,1]$. We also choose the set of historical assortments $\Scal = \{ S_1,\ldots,S_M \}$ uniformly at randomly from all possible $2^N$ assortments of $N$ products. We examine the performance of the column randomization method under various problem sizes, using different values of $N$ and $M$. For the CG method, we use the method in \cite{mivsic2016data}, and solve the subproblem as an integer program (IP) from \cite{van2015market}.

Table~\ref{table:ranking_estimation_CR_vs_CG} shows the performance of the column randomization method. The first two columns of the table indicate the problem size. The third column shows the number of sampled columns. The fourth and fifth columns display the objective value and the runtime, respectively; for both of these metrics, we report the average value of the metric over 20 runs of the column randomization method. The sixth column denotes the duration of the CG method to reach the same (average) objective value as the column randomization method. We remark that the optimal objective value $v(P^{\text{EST}})$ is always zero, since random utility maximization models such as the MNL model can be represented as ranking-based models \citep{block1959random}. Thus, instead of showing relative optimality gap as in Table~\ref{table:cutting_stock_CR_vs_CG}, we directly show the objective value of the column-randomized linear program in Table~\ref{table:ranking_estimation_CR_vs_CG}.

In all cases listed in Table~\ref{table:ranking_estimation_CR_vs_CG}, the column randomization method outperforms the CG method by a large margin. It only requires a fraction of the runtime of the CG method to reach the same optimality level. In particular, when $(N,M)=(10,150)$, the column randomization method only needs three seconds to reach the optimal objective value, which is zero, while the CG method needs over ten thousand seconds (almost three hours). In real-world applications, the number of products $N$ is usually significantly larger than $10$. In those cases, the advantage of column randomization will be even more pronounced. We note that in the IP formulation of the CG subproblem, the number of binary variables scales as $O(N^2 + NM)$. Thus, as $N$ increases, the subproblem quickly becomes intractable (\cite{van2015market} showed this subproblem to be NP-hard).

Lastly, we comment on why column randomization performs well for the ranking-based choice estimation problem. For the ranking-based choice estimation problem, it is well-known that the problem is extremely underdetermined. Specifically, for a fixed collection of $m$ assortments with choice probabilities generated according to a random utility maximization problem, there can be multiple distributions $\lambdab$ that solve $\Ab \lambdab = \vb$, where $\Ab = [\alphab_1 \ \alphab_2 \ \cdots \alphab_{(N+1)!}]$; see \cite{farias2013nonparametric,van2015market,mivsic2016data,sturt2021value}. This multiplicity arises in two ways. First, we can find two different collections of columns $\alphab_{j_1}, \dots, \alphab_{j_K}$ and $\alphab_{\tilde{j}_1}, \dots, \alphab_{\tilde{j}_K}$, each of which can be used to perfectly fit the vector $\vb$. Second, even for a fixed collection of columns $\alphab_{j_1},\dots, \alphab_{j_K}$, each column $\alphab_j$ could be rationalized by more than one ranking; in other words, there could exist many rankings that give rise to the same column. (As a simple example of this, suppose that $N = 4$, $M = 2$ and $S_1 = \{1,2\}$, $S_2 = \{3,4\}$, and consider the column $\alphab = (1,0,0,1,0,0)$ corresponding to the option-assortment pairs $( (1,1), (2,1), (0,1), (3,2), (4,2), (0,2) )$. This column can correspond to the ranking $1 \prec 3 \prec 2 \prec 4 \prec 0$, which prefers product 1 the most, but can also correspond to the ranking $3 \prec 1 \prec 2 \prec 4 \prec 0$, which prefers product 3 the most. In fact, any ranking that obeys $1 \prec 2$ and $3 \prec 4$ will be consistent with $\alphab$.) Consequently, the ranking-based choice estimation problem will often have an extremely large number of optimal or near-optimal solutions that are diverse in terms of which variables (rankings) those solutions are supported on. Although $P^{\text{EST}}$ is not a standard form LP, this type of structure and the numerical performance exhibited in Table~\ref{table:ranking_estimation_CR_vs_CG} is consistent with Theorem~\ref{theorem:many_BFS_P_distr_bound}, which suggests that column randomization will do well in the presence of many diverse, near-optimal solutions.

\begin{table}[]
	\SingleSpacedXI
	\centering
	\small
	\begin{tabular}{rrrrrr}
		\toprule
		\multicolumn{1}{c}{{$N$}} &\multicolumn{1}{c}{{$M$}}    &  \multicolumn{1}{c}{{Columns $(K)$}}     & \multicolumn{1}{c}{{Objective}}          & \multicolumn{1}{c}{{Runtime (s)}}     & \multicolumn{1}{c}{{CG Runtime (s)}}    \\
		\midrule
		6 &$50$ & $500$  & 0.05 & 0.03   & 20.58    \\
		& & $1000$  & 0.00 & 0.07   & 30.44   \\ \midrule %
		8& 50 & 500  & 0.13 & 0.10  & 52.32    \\
		& & 1000 & 0.00  & 0.12  & 88.25  
		\\ \midrule
		8 & 100 & 500  & 0.92 & 0.21   & 120.14   \\
		& & 1000  & 0.07 & 0.45  & 414.43   \\
		& & 1500 & 0.00 & 0.66  & 632.23   \\
		\midrule
		10 & 50 & 500  & 0.27 & 0.17  & 11.93   \\
		& & 1000  & 0.00  & 0.22  & 282.78   \\ \midrule
		10& 100 & 500  & 1.60 & 0.28  & 240.23  \\
		& & 1000 & 0.40 & 0.53 & 774.66
		\\ & & 1500 & 0.06 & 0.71 & 1423.71 
		\\ & & 2000 & 0.00 & 1.57 & 2234.52 \\ \midrule
		10& 150 & 500  & 2.91 & 0.69  & 507.63  \\
		& & 1000 & 0.98 & 1.07 & 1399.22
		\\ & & 1500 & 0.43 & 1.33 & 2635.36 
		\\ & & 2000 & 0.18 & 2.01 & 4524.72
		\\ & & 2500 & 0.00 & 3.14 & 10143.93
		\\ \bottomrule
	\end{tabular}
	
	\caption{Performance of the column randomization method on the estimation problem $P^{\text{EST}}$ under varying problem sizes and numbers of sampled columns. 	\label{table:ranking_estimation_CR_vs_CG}}
\end{table}

Finally, we have included supplementary numerical experiments in Section~\ref{sec:numerics_NCME_appendix} of the e-companion. In Section~\ref{subsec:numerics_NCME_uniform_vs_MNL}, we compare the method's performance under a different sampling scheme. In Section~\ref{subsec:numerics_NCME_CRthenCG}, we explore the benefits of combining the column randomization method with CG.

\section{Conclusion}
\label{sec:conclusion}

In this paper, we analyzed the column-randomization method for solving large-scale linear programs with an intractably large number of columns, which involves simply randomly sampling a collection of $K$ columns from the constraint matrix, and solving the corresponding problem. We developed performance guarantees for the solution one obtains from this approach. We derived an upper bound on the optimality gap that holds with high probability.  This bound converges at a rate $1 / \sqrt{K}$, where $K$ is the number of sampled columns, to the optimality gap of a linear program that we named as distributional counterpart. We further analyzed the gap of the distributional counterpart and discussed conditions under which this gap will be small. In numerical experiments with the cutting stock problem and the nonparametric choice model estimation problem, we showed that the proposed approach can obtain near-optimal solutions in a fraction of the computational time required by column generation. Given the computational simplicity of randomly sampling columns in many problems, we hope that this paper will spur further research into large-scale optimization that leverages the synergy of randomization and optimization. 

\section*{Acknowledgments}
We sincerely thank the area editor Daniel Kuhn, the associate editor, and the three anonymous referees for their thoughtful comments that helped to strengthen this work. The authors also thank Vishal Gupta for helpful comments on an early version of this work.

\bibliographystyle{plainnat}
\bibliography{sampling_method.bib}

\ECSwitch

\ECDisclaimer

\ECHead{Electronic companion for ``Column-Randomized Linear Programs: Performance Guarantees and Applications'' by Akchen and Mi\v{s}i\'{c}}

\begin{appendices}

\DoToC

\hspace{1cm}\\

The ecompanion is organized as follows. Section~\ref{sec:additional_discussion_on_thm1_thm2} follows Section~\ref{subsec:discussion_theorems} and continues the discussion of Theorems~\ref{thm:main_largest_abs_of_all_dual_BFS} and~\ref{thm:main_reduced_cost}. Section~\ref{sec:proofs} completes all omitted proofs for the theoretical results in Section~\ref{sec:main_results}. Section~\ref{sec:proofs_distributional_counterpart} provides the proofs for all theoretical results in Section~\ref{sec:P_distr} regarding the distributional counterpart and provides additional comments. Section~\ref{sec:special_structures} shows that the parameters $\gamma$ in Theorem~\ref{thm:main_largest_abs_of_all_dual_BFS} can be further obtained for several applications. Section~\ref{sec:dependent_columns} extends the proposed framework by proposing sampling statistically-dependent columns in the column randomization method. Section~\ref{sec:cutting_stock_continued} follows Section~\ref{sec:numerics_CS} and  completes the numerical experiments on the cutting stock problem. Section~\ref{sec:numerics_NCME_appendix} follows Section~\ref{sec:numerics_NCME} and completes the numerical experiments on the nonparametric choice model estimation. Section~\ref{sec:comparison} provides a detailed comparison between the proposed framework and other large-scale LP solvers based on randomized algorithms \citep{agrawal2014dynamic,vu2018random}. All of the code is available at the repository $\texttt{Column-Randomized\_LP}$ at 
\begin{center}
\tt \href{https://github.com/yi-chun-akchen/Column-Randomized_LP}{https://github.com/yi-chun-akchen/Column-Randomized\_LP}.
\end{center}

\hspace{1cm}

\section{Additional Discussion on Results in Section~\ref{sec:main_results}}
\label{sec:additional_discussion_on_thm1_thm2}

This section continues the discussion in Section~\ref{subsec:discussion_theorems}.

\subsubsection*{Lower Bound on $v(P_J)$:}

We note that neither Theorem~\ref{thm:main_largest_abs_of_all_dual_BFS} nor \ref{thm:main_reduced_cost} implies that the optimality gap $\Delta v(P_J)$ of the column-randomized linear program $P_J$ can be arbitrarily small with large $K$. Indeed, if $\xib$ is not ``comprehensive'' enough -- that is, its support is small, and does not include the complete set of columns of any optimal basis for $P$ -- then $\Delta v(P_{\distr}) > 0$ no matter what $C$ is, and one would not expect the column-randomized program $P_J$ to perform closely to the complete problem $P$, even if $K$ is large. We can formalize this intuition in the following proposition, where $I^+$ denotes the support of the distribution $\xib$.
\begin{proposition}
	\label{proposition:P_J_lower_bound} 
	Define $I^+ = \{  j \in [n]  \mid \xi_j > 0 \}$ and let $P^+_{\distr} \equiv \min \{ \cb^T_{I^+} \xb^+ \mid \Ab_{I^+} \xb^+ = \bb, \xb^+ \geq \zerob \}$. Then $v(P_J) \geq v(P^+_{\distr})$ almost surely, and $v(P_J) \to v(P^+_{\distr} )$ almost surely as $K \to \infty$.  
\end{proposition}
The proof is straightforward and omitted for brevity, as any solution of $P_J$ can be reformulated as a feasible solution of $P^+_\distr$, and as $K \to \infty$, every column in $I^+$ is sampled at least once almost surely, ensuring that $v(P_J) = v(P^+_{\distr})$. An obvious consequence of this proposition is that if $I^+ = [n]$, i.e., every column has a positive probability of being sampled, then $v(P_J)$ will converge to $v(P)$ when enough columns are sampled. From this perspective, the value of our bounds in Theorems~\ref{thm:main_largest_abs_of_all_dual_BFS} and \ref{thm:main_reduced_cost} is that they provide finite sample guarantees, for the case where $K \ll n$ and it is impossible that one will have sampled all of the $n$ columns.

\subsubsection*{Feasibility of $P_J$:}

We make several important remarks regarding the feasibility of $P_J$ and how feasibility is incorporated in our guarantee. First, note that in general, the sampled problem $P_J$ need not be feasible. As a simple example, consider the following complete problem:
\begin{equation*}
P = P_{\Ib} \equiv \min \{ \oneb^T \xb \mid \Ib \xb = \oneb, \xb \geq \zerob \},
\end{equation*} 
where $\Ib$ is the $n$-by-$n$ identity matrix and $m = n$. In this problem, the only way that the sampled problem $P_J$ can be feasible is if the collection $j_1,\dots,j_K$ includes every index in $[n]$; if any column $j \in [n]$ is not part of the sample $J$, then the sampled problem $P_J$ is automatically infeasible. Thus, when $K < n$, $P_J$ is infeasible almost surely. When $K \geq n$, it is still possible that $j_1, \dots, j_K$ does not include all indices in $[n]$, and thus $P_J$ is infeasible with positive probability.

For this reason, our guarantee on the optimality gap is stated as a \emph{conditional} guarantee: with high probability over the sample $j_1,\dots,j_K$, the optimality gap of $P_J$ obeys a particular bound \emph{if the column-randomized LP is feasible}. Formally, our two guarantees can be represented as
\begin{equation*}
\Pr \left[ \left\{ P_J\ \text{is feasible} \right\} \Rightarrow \left\{ \Delta v(P_J) \leq \Delta v(P_{\distr}) + \frac{C}{\sqrt{K}} \cdot C_{P} \cdot C_{\delta} \right\} \right] \geq 1 - \delta.
\end{equation*}
Since the implication $A \Rightarrow B$ is logically equivalent to $A^C \cup B$, an alternative equivalent restatement of the general form of our guarantee is 
\begin{equation*}
\Pr \left[ \left\{P_J\ \text{is infeasible} \right\} \ \cup\ \left\{ \Delta v(P_J) \leq \Delta v(P_{\distr}) + \frac{C}{\sqrt{K}} \cdot C_{P} \cdot C_{\delta} \right\} \right] \geq 1 - \delta.
\end{equation*}
We note that this type of guarantee is distinct from probabilistically conditioning on $j_1,\dots,j_K$, i.e., our guarantee is \emph{not} the same as
\begin{equation*}
\Pr\left[ \Delta v(P_J) \leq \Delta v(P_{\distr}) + \frac{C}{\sqrt{K}} \cdot C_{P} \cdot C_{\delta}  \  \vline \  P_J\ \text{is feasible} \right] \geq 1 - \delta,
\end{equation*}
because upon conditioning on the feasibility of $P_J$, the random variables $j_1,\dots,j_K$ are in general no longer an i.i.d. sample. As an example of this, consider again problem $P_{\Ib}$ above, with $K = n$ and a randomization scheme $\rho$ corresponding to the uniform distribution $\xib = (1/n, \dots, 1/n)$ over $[n]$. By conditioning on the event that $P_J$ is feasible, the sample $J = \{ j_1,\dots, j_K\}$ must then be exactly equal to $[n]$, and we obtain that $\Pr[ j_k = t, j_{k'} = t] = 0 \neq \Pr[ j_k = t] \cdot \Pr[ j_{k'} = t]$ for any $k, k' \in [K]$ with $k \neq k'$ and $t \in [n]$. In this example, the indices $j_1,\dots, j_K$ are thus not independent. 

With regard to the feasibility of column-randomized LPs, it appears to be difficult to guarantee feasibility in general. However, one can use similar techniques as in the proofs of our main results to characterize the near-feasibility of a column-randomized LP. Consider the following complete problem, and its sampled and distributional counterparts:
\begin{align*}
P^{\feas} & = \min \{ \| \Ab \xb - \bb \|_1 \mid \xb \geq \zerob \}, \\
P^{\feas}_J & = \min \{ \| \Ab_J \tilde{\xb} - \bb \|_1 \mid \tilde{\xb} \geq \zerob \}, \\
P^{\feas}_{\distr} & = \min \{ \| \Ab \xb - \bb \|_1 \mid \zerob \leq \xb \leq C \xib \}.
\end{align*}
The objective function in each problem measures how close $\Ab \xb$ is to $\bb$ for a given nonnegative solution $\xb$, and the optimal value measures the minimum total infeasibility, as measured by the lowest attainable $\ell_1$ distance between $\Ab \xb$ and $\bb$. Note that an optimal value of zero for a given problem implies that the feasible region contains a solution $\xb$ that satisfies $\Ab \xb = \bb$. With a slight abuse of notation, let us use $v(P^{\feas})$, $v(P^{\feas}_J)$ and $v(P^{\feas}_{\distr})$ to denote the optimal objective value of each problem. We then have the following result.
\begin{proposition}
	\label{prop:near_feasibility_bound}
	Let $C$ be a nonnegative constant. For any $\delta \in (0,1)$, with probability at least $1 - \delta$ over the sample $J$,
	\begin{equation*}
	v(P^{\feas}_J) \leq v(P^{\feas}_{\distr}) + \frac{C}{\sqrt{K}} \cdot m \cdot \|\Ab\|_{\max} \cdot \left(1 + \sqrt{2 \log \frac{1}{\delta}} \right). 
	\end{equation*}
\end{proposition}
The proof of Proposition~\ref{prop:near_feasibility_bound} (see Section~\ref{proof:near_feasibility_bound} of the ecompanion) follows using a similar but simpler procedure than those used in the proofs of Theorems~\ref{thm:main_largest_abs_of_all_dual_BFS} and \ref{thm:main_reduced_cost}. The guarantee in Proposition~\ref{prop:near_feasibility_bound} has a similar interpretation to Theorems~\ref{thm:main_largest_abs_of_all_dual_BFS} and \ref{thm:main_reduced_cost}: the magnitude of the total infeasibility of the columns $J$ is bounded with high probability by the minimum infeasibility of the distributional counterpart $P^{\feas}_{\distr}$ plus a $O(1 / \sqrt{K})$ term.

\subsubsection*{Feasibility-guaranteed column randomization algorithm:} One practical way in which one can modify Algorithm~\ref{alg:main} to ensure that the sampled problem is always feasible is to augment the column set $J$ with a set of columns $J_F$ such that $P_{J \cup J_F}$ is a feasible problem. We define the new procedure, Algorithm~\ref{alg:main_FG}, below.

\begin{algorithm}
	\SingleSpacedXI
	\caption{The Feasibility Guaranteed Column Randomization Method}
	\label{alg:main_FG}
	\begin{algorithmic}[1]
		\STATE Set $J_F \subseteq [n]$ to be a set of columns such that $P_{J_F}$ is feasible and $\rank(\Ab_{J_F}) = m$. 
		\STATE Sample $K$ indices as $J \equiv \{ j_1,\ldots,j_K \}$ by a randomization scheme $\rho$.
		
		\STATE Define $\Ab_{J \cup J_F} = [A_{j}]_{j \in J \cup J_F}$ and $\cb_{J \cup J_F} = [c_j]_{j \in J \cup J_F}$.
		\STATE Solve the column-randomized linear program, which has $|J \cup J_F|$ columns:
		\begin{align}
		\label{problem:LP_sampled_FG}
		P_{J \cup J_F}:  \quad \min \left\lbrace \cb^T_{J \cup J_F} \tildexb \mid \Ab_{J \cup J_F} \tildexb = \bb, \,\,\, \tildexb \geq \zerob \right\rbrace. 
		\end{align}
		\RETURN optimal objective value $v(P_{J \cup J_F})$ and an optimal solution $\tildexb^*$.
	\end{algorithmic}
\end{algorithm}

There are two important aspects of Algorithm~\ref{alg:main_FG} to be cognizant of. First, for this procedure, we can adapt Theorems~\ref{thm:main_largest_abs_of_all_dual_BFS_FG} and Theorem~\ref{thm:main_reduced_cost_FG} so as to obtain guarantees on $\Delta v(P_{J\cup J_F})$. This results in the following two guarantees; importantly, these guarantees are no longer conditional guarantees. 

\begin{theorem}
	\label{thm:main_largest_abs_of_all_dual_BFS_FG}
	Let $J$ and $J_F$ be as defined in Algorithm~\ref{alg:main_FG}. For any $\delta \in (0,1)$, with probability at least $1 - \delta$ over the sample $J$, then
	\begin{align}
	\label{eq:thm_convergence_with_abs_of_BFS_FG}
	\Delta v(P_{J \cup J_F})   \leq \Delta  v(P_\distr) +  \frac{C \left( 1 + m \gamma \| \Ab \|_{\max} \right)}{\sqrt{K}} \left(  1 + \sqrt{2 \log \frac{2}{\delta}}  \right),
	\end{align}
\end{theorem}

\begin{theorem}
	\label{thm:main_reduced_cost_FG}
	Let $J$ and $J_F$ be as defined in Algorithm~\ref{alg:main_FG}. For any $\delta \in (0,1)$, with probability at least $1 - \delta$ over the sample $J$, then
	\begin{align}
	\Delta v(P_{J \cup J_F}) \leq \Delta v(P_{\distr}) + \frac{C}{\sqrt{K}} \cdot \chi \cdot \left(1 + \sqrt{2 \log \frac{1}{\delta} } \right)
	\end{align}
\end{theorem}

Second, in the statement of Algorithm~\ref{alg:main_FG}, the first step is to obtain a set of columns $J_F$ so that $P_{J_F}$ is feasible and $\Ab_{J_F}$ has full row rank. This can be accomplished easily by applying the first phase of the two-phase method, which is a standard method for obtaining an initial feasible solution to a linear program for which an initial basis is not obvious (see Chapter 3, Section 5 of \citealt{bertsimas1997introduction}). In particular, one formulates the following problem:
\begin{subequations}
	\begin{alignat}{2}
	& \underset{\xb, \epsilonb^+, \epsilonb^-}{\text{minimize}} & \quad & \oneb^T \epsilonb^+ + \oneb^T \epsilonb^- \\
	& \text{subject to} & & \Ab \xb + \Ib \epsilonb^+ - \Ib \epsilonb^- = \bb, \\
	& & & \xb \geq \zerob, \\
	& & & \epsilonb^+, \epsilonb^-
	\end{alignat}
	\label{prob:first_phase_two_phase_method}%
\end{subequations}
where $\epsilonb^+, \epsilonb^- \in \mathbb{R}^m$ and $\zerob$ is an appropriately sized vector of zeros in the two nonnegativity constraints. As with the main problem $P$, this problem obviously cannot be formulated explicitly, but can be solved using column generation. At the start of column generation, we do not include any columns from $\Ab$, and the principal constraint is $\Ib \epsilonb^+ - \Ib \epsilonb^- = \bb$, for which an initial basis can be found trivially: for each $i \in [m]$, set $\epsilon^+_i = \max\{b_i, 0\}$ and $\epsilon^-_i = \min\{b_i, 0\}$. As column generation progresses, the objective value will decrease and the $\epsilon^+_i$ and $\epsilon^-_i$ variables will gradually leave the basis. Upon termination, one will obtain a basic feasible solution for which all of the $\epsilon^+_i, \epsilon^-_i$ variables are non-basic. The resulting set of basic columns of $\xb$, $J_F$, is such that $P_{J_F}$ is feasible and $\rank(\Ab_{J_F}) = m$. Observe now that if $P_{J_F}$ is feasible, then $P_{J \cup J_F}$, for \emph{any} set of columns $J \subseteq [n]$, must also be feasible: one can take any solution $\xb_{J_F}$ to $P_{J_F}$ and set $x_j = 0$ for all $j \in J \setminus J_F$, resulting in a solution $\xb_{J \cup J_F}$ that is nonnegative, and satisfies $\Ab_{J \cup J_F} \xb_{J \cup J_F} = \bb$. In addition, the matrix $\Ab_{J \cup J_F}$ must also have rank $m$. Although this approach requires column generation, it is reasonable to expect that column generation applied to problem~\eqref{prob:first_phase_two_phase_method} to find an initial solution should generally be faster than when it is applied to the complete problem $P$. Proofs of Theorem~\ref{thm:main_largest_abs_of_all_dual_BFS_FG} and~\ref{thm:main_reduced_cost_FG} can be found in Section~\ref{subsec:feasbility_guarantee_proof}.

\subsubsection*{Interpretation of $\gamma$ and $\chi$:}

We first note that the technique of bounding the objective value of a linear program using the $\ell_\infty$ norm of basic feasible solutions has been applied previously in the literature \citep{ye2011simplex,kitahara2013bound}. The presence of $\gamma$ and $\chi$ in Theorem~\ref{thm:main_largest_abs_of_all_dual_BFS} and \ref{thm:main_reduced_cost}, respectively, arises due to the use of sensitivity analysis results from linear programming with respect to the right-hand side vector $\bb$. As we discuss in the proof in Section~\ref{sec:proofs}, any optimal solution $\xb^{*0}$ of problem $P_\distr$ has a sparse counterpart $\xb'$ in the space $\Scal_J \equiv \{ \xb \mid x_j = 0 \,\,\, \forall j \notin J \}$ such that $\xb'$ is in the vicinity of $\xb^{*0}$ in terms of Euclidean distance. However, $\xb'$ does not necessarily belong to the feasible set $\Fcal(P_J)$ of the column-randomized linear program $P_J$, since $\Fcal(P_J)$ is a subset of $\Scal_J$. To relate the optimal objective value $v(P_J)$ of problem $P_J$ to $\cb^T \xb'$, which is close to $\cb^T \xb^{*0}$, we use sensitivity analysis arguments which involve either $\gamma$ or $\chi$.

\subsubsection*{Comparison of Theorems~\ref{thm:main_largest_abs_of_all_dual_BFS} and \ref{thm:main_reduced_cost}:}

While both Theorem~\ref{thm:main_largest_abs_of_all_dual_BFS} and \ref{thm:main_reduced_cost} provide valid bounds for the optimality gap $\Delta v(P_J)$, Theorem \ref{thm:main_largest_abs_of_all_dual_BFS} is in general easier to apply; indeed, in Section~\ref{sec:special_structures} we discuss two notable examples where $\gamma$ can be easily computed (specifically, LPs with totally unimodular constraint matrices $\Ab$ and infinite horizon discounted Markov decision processes). For problems that are not standard form LPs, neither guarantee directly applies, but we can obtain specialized guarantees by carefully modifying a result (Proposition~\ref{prop:objective_bound_by_dual_solution} in Section~\ref{subsec:main_theorem_proof}) that leads to Theorem~\ref{thm:main_largest_abs_of_all_dual_BFS} and designing bounds for the $\ell_{\infty}$ norm of feasible or optimal solutions of $D_J$ (as opposed to basic solutions of $D$). We will later showcase two examples of such guarantees, for covering LPs (Section~\ref{subsec:special_structures_covering}) and packing LPs (Section~\ref{subsec:special_structures_packing}). 

With regard to Theorem~\ref{thm:main_reduced_cost}, we expect for most problems that Theorem~\ref{thm:main_reduced_cost} will be difficult to apply, as it requires a universal bound for the norm of the reduced cost vector for every basis, feasible or not, of problem $P$. Nevertheless, Theorem~\ref{thm:main_reduced_cost} is interesting because it involves reduced costs, which are also of importance in column generation. For a basic feasible solution, the reduced cost of a non-basic variable $j$ can be thought of as the rate at which the objective changes as one increases $x_j$ to move from the current basic feasible solution to an adjacent/neighboring basic feasible solution in which $j$ is part of the basis. With this perspective of reduced costs, one can informally interpret the result in the following way: if $\chi$ is small, then the rate at which the objective changes between adjacent basic feasible solutions is small. In such a setting, it is reasonable to expect that there will be many basic feasible solutions that are close to being optimal and that solving the sampled problem $P_J$ should return a solution that performs well. On the other hand, if there exist non-optimal basic feasible solutions where the reduced cost vector has a very large magnitude (which would imply a large $\chi$), then this would suggest that the objective changes by a large amount between certain adjacent basic feasible solutions, and that there are certain ``good'' columns that are more important than others for achieving a low objective value. In this setting, we would expect the sampled problem objective $v(P_J)$ to only be close to $v(P)$ if $J$ includes the ``good'' columns, which would be unlikely to happen in general.

\subsubsection*{Design of Randomization Scheme $\rho$:}

The quantity $\xi_j$, which is the probability that the $j$th column is drawn by the randomization scheme $\rho$, can be interpreted as the relative importance of $x_j$ compared to other components of $\xb \in \mathbb{R}^n$ in the complete problem $P$; indeed, when the corresponding column is randomly chosen, $x_j$ is allowed to be nonzero, and can thus be utilized to solve the optimization problem. 
For example, in a network flow optimization problem, $x_j$ represents the amount of flow over edge $j$; a nonzero $\xi_j$ can thus be interpreted as the belief that edge $j$ should be used for flow. 
As another example, consider the LP formulation of an MDP, where each component of $\xb$ corresponds to a state-action pair $(s,a)$ (i.e., $x_{(s,a)}$ is the expected discounted frequency of the system being in state $s$ and action $a$ being taken). In this setting, a nonzero $\xi_{(s,a)}$ can be interpreted as the relative importance of $(s,a)$ to other state-action pairs. %

One can design the randomization scheme based on prior knowledge of the problem. For example, one could use a heuristic solution to a network flow problem to design a randomization scheme $\rho$ resulting in a distribution $\xib$ that is biased towards this heuristic solution. Similarly, if one has access to a good heuristic policy for an MDP, one can design a distribution $\xib$ that is biased towards state-action pairs $(s,a)$ that occur frequently for this policy. If such prior knowledge is not available, a uniform or nearly-uniform distribution over $[n]$ is adequate.
We provide several concrete examples on how to design randomization schemes in our numerical experiments in Sections~\ref{sec:numerics_CS} and \ref{sec:numerics_NCME}. Finally, we note that the indices in $J$ have been assumed to be i.i.d. In Section~\ref{sec:dependent_columns}, we derive analogous guarantees for the case when the indices are sampled non-independently.

\subsubsection*{Minor Remarks on the Upper Bound:}

We mention two other interesting properties of the bound~\eqref{eq:bound_general_form}. First, the second term in \eqref{eq:bound_general_form} is independent of the distribution $\xib$; no matter how $\xib$ is designed, the optimality gap $\Delta v(P_J)$ is guaranteed to converge with rate $1/\sqrt{K}$. Second, the dependence of the bound on the confidence parameter $\delta$ is via $\sqrt{2 \log (2/ \delta)}$ in Theorem~\ref{thm:main_largest_abs_of_all_dual_BFS} or $\sqrt{2 \log (1/ \delta)}$ in Theorem~\ref{thm:main_reduced_cost}. This implies that very small values of $\delta$ will not significantly increase the upper bound on $\Delta v(P_J)$.

\section{Omitted Proofs of Results in Section~\ref{sec:main_results}}
\label{sec:proofs}

In this section, we prove Theorem~\ref{thm:main_largest_abs_of_all_dual_BFS} and \ref{thm:main_reduced_cost}. We start with some preliminary results (Section~\ref{subsec:preliminary_results}) then prove the main theorems (Section~\ref{subsec:main_theorem_proof}).

\subsection{Preliminary Results and Lemmas}
\label{subsec:preliminary_results}

Lemma~\ref{lemma:averaged_point_and_distance} and \ref{lemma:averaged_point_and_distance_L1} bound the distance between the sample mean and the expected value of a collection of i.i.d. vectors, in terms of $\ell_2$ norm and $\ell_1$ norm, respectively. Lemma~\ref{lemma:averaged_point_and_distance} is Lemma 4 from \cite{rahimi2009weighted}, which utilizes McDiarmid's inequality to show that the scalar function $\|  \bar{\wb} - \Ebb \left[\bar{\wb} \right] \|_2$, where $\bar{\wb}$ is the mean of $K$ i.i.d. vectors $\wb_1,\dots, \wb_K$, concentrates to zero with rate $O\left( 1 / \sqrt{K} \right)$.

\begin{lemma}\citep{rahimi2009weighted}
	\label{lemma:averaged_point_and_distance}
	Let $\wb_1,\wb_2,\ldots,\wb_K$ be i.i.d. random vectors such that $\| \wb_k\|_2 \leq C$ for $k = 1,\ldots,K$. Let $\bar{\wb} = (1/K) \cdot \sum_{k=1}^K \wb_k$. Then for any $\delta \in (0,1)$, we have, with probability at least $1 - \delta$,
	\begin{align*}
	\|  \bar{\wb} - \Ebb \left[\bar{\wb} \right] \|_2 \leq \frac{C}{\sqrt{K}}
	\cdot \left(  1 + \sqrt{2 \log \frac{1}{\delta}}  \right).
	\end{align*}
\end{lemma}

\begin{lemma}
	\label{lemma:averaged_point_and_distance_L1}
	Let $\wb_1,\wb_2,\ldots,\wb_K$ be i.i.d. random vectors of size $m$ such that $\| \wb_k\|_\infty \leq C$ for $k = 1,\ldots,K$. Let $\bar{\wb} = (1/K) \cdot \sum_{k=1}^K \wb_k$. Then for any $\delta \in (0,1)$, we have, with probability at least $1 - \delta$,
	\begin{align*}
	\|  \bar{\wb} - \Ebb \left[ \bar{\wb} \right] \|_1 \leq \frac{mC}{\sqrt{K}}
	\cdot \left(  1 + \sqrt{2 \log \frac{1}{\delta}}  \right).
	\end{align*}
\end{lemma}

{\it Proof:} Since $\| \wb_k \|_2 \leq \sqrt{m} \|  \wb_k\|_\infty \leq \sqrt{m} C$, we apply Lemma~\ref{lemma:averaged_point_and_distance} and obtain that with probability at least $1 - \delta$, $\|  \bar{\wb} - \Ebb \left[ \bar{\wb} \right] \|_2 \leq \sqrt{m} \cdot {C}/{\sqrt{K}}
\cdot \left(  1 + \sqrt{2 \log \frac{1}{\delta}}  \right).$
Combining this with the fact that $\|  \bar{\wb} - \Ebb \left[\bar{\wb} \right] \|_1 \leq \sqrt{m} \cdot \|  \bar{\wb} - \Ebb \left[\bar{\wb} \right] \|_2$, we obtain the desired result. \hfill \Halmos \\

Lemma~\ref{lemma:sensitivity_of_obj} is a standard result of sensitivity analysis of linear programming; see Chapter 5 of \cite{bertsimas1997introduction}. In fact, one can view the optimal objective value of problem $P$ as a convex function in $\bb$ and show that any optimal dual solution $\pb$ is a subgradient at $\bb$.

\begin{lemma}
	\label{lemma:sensitivity_of_obj}
	Let $z(\bb) = \min \left\lbrace \cb^T_0 \yb \mid \Ab_0 \yb = \bb, \yb \geq \zerob \right\rbrace$ and $z(\bb') = \min \left\lbrace \cb^T_0 \yb \mid \Ab_0 \yb = \bb', \yb \geq \zerob \right\rbrace$. Then $z(\bb) - z(\bb') \leq \pb^T (\bb - \bb')$, where $\pb$ is an optimal dual solution of the former problem.
\end{lemma}

\subsection{Proofs of Theorem \ref{thm:main_largest_abs_of_all_dual_BFS} and \ref{thm:main_reduced_cost}}
\label{subsec:main_theorem_proof}

We first establish a useful result.
\begin{proposition}
	\label{prop:objective_bound_by_dual_solution}
	Let $C$ be a nonnegative constant and define the linear program $P_\distr$ as in Theorem~\ref{thm:main_largest_abs_of_all_dual_BFS}, i.e., $P_\distr: \,\,\, \min \left\lbrace \cb^T \xb \mid \Ab \xb = \bb, \zerob \leq \xb \leq C \xib \right\rbrace$. Let $P_J$ be the column-randomized LP solved by Algorithm~\ref{alg:main}. For any $\delta \in (0,1)$, with probability at least $1 - \delta$ over the sample $J$, the following holds: if $P_J$ is feasible, then 
	\begin{equation*}
	\Delta v(P_J) \leq \Delta v(P_{\distr}) + \frac{C}{\sqrt{K}} \cdot (1 + \| \pb \|_{\infty} \cdot m \cdot \| \Ab \|_{\max} ) \cdot \left(1 + \sqrt{2 \log \frac{2}{\delta} } \right)
	\end{equation*}
	for any optimal solution $\pb$ of problem $D_J$ (the dual of problem $P_J$).
	
\end{proposition}

\proof{Proof:} Let $j_1,\dots, j_K$ be the set of indices sampled according to the distribution $\xib$ by the randomization scheme $\rho$. Let $\xb^{*0}$ be an optimal solution of the distributional counterpart problem $P_{\distr}$. Consider the solution $\xb'$ that is defined as
\begin{equation*}
\xb' \equiv \frac{1}{K} \sum_{k=1}^K \frac{x^{*0}_{j_k}}{ \xi_{j_k}} \cdot \eb_{j_k},
\end{equation*}
where we use $\eb_{j}$ to denote the $j$th standard basis vector for $\Rbb^n$. In addition, define the vector $\bb'$ as 
\begin{equation*}
\bb' \equiv \Ab \xb'.
\end{equation*}
To prove our result, we proceed in three steps. In the first step, we show how we can probabilistically bound $\| \xb' - \xb^{*0} \|_2$. In the second step, we show how we can probabilistically bound $\| \bb' - \bb \|_1$. In the last step, we use the results of our first two steps, together with sensitivity results for linear programs, to derive the required bound. In what follows, we use $I_+$ to denote the support of $\xib$, that is, $I_+ = \{ j \in [n] \mid \xi_j > 0\}$. \\

\textbf{Step 1: Bounding $\| \xb' - \xb^{*0} \|_2$.} To show that $\xb'$ will be close to $\xb^{*0}$, let us first define the vector $\wb_k$ as 
\begin{equation*}
\wb_k = \frac{x^{*0}_{j_k}}{\xi_{j_k}} \cdot \eb_{j_k}
\end{equation*}
for each $k \in [K]$.  The vectors $\wb_1,\dots, \wb_K$ constitute an i.i.d. collection of vectors, and possess three special properties. First, observe that $\xb'$ is just the sample mean of $\wb_1, \dots, \wb_K$. Second, observe that the expected value of each $\wb_k$ can be calculated as
\begin{align*}
\Ebb[ \wb ] & = \sum_{j \in I_+} \xi_j \cdot \frac{x^{*0}_j}{\xi_j} \cdot \eb_j = \sum_{j \in I_+} x^{*0}_j \eb_j  = \sum_{j \in [n]} x^{*0}_j \eb_j  = \xb^{*0}
\end{align*}
where we use $\wb$ to denote a random vector following the same distribution as each $\wb_k$. In the above, we note that the third step follows because the distributional counterpart $P_{\distr}$ includes the constraint $\xb \leq C \xib$, so $j \notin I_+$ automatically implies that $x^{*0}_j = 0$. 

Finally, observe that the $\ell_2$ norm of each $\wb_k$ can be bounded as 
\begin{align*}
\| \wb_k \|_2 & = \left| \frac{ x^{*0}_{j_k}}{\xi_{j_k}} \right| \cdot \| \eb_{j_k} \|_2 \leq C \cdot 1 = C,
\end{align*}
where the inequality follows because $\xb^{*0}$ satisfies the constraint $\zerob \leq \xb \leq C \xib$. With these three properties in hand, and recognizing that $\| \xb' - \xb^{*0} \|_2 = \| (1/K) \sum_{k=1}^K \wb_{k} - \Ebb[\wb] \|_2$, we can invoke Lemma~\ref{lemma:averaged_point_and_distance} to assert that, with probability at least $1 - \delta/2$,
\begin{equation}
\| \xb' - \xb^{*0} \|_2 \leq \frac{C}{\sqrt{K}} \cdot \left(  1 + \sqrt{2 \log \frac{2}{\delta}}  \right). \label{eq:bound_sampled_mean_in_x} \\
\end{equation}

\textbf{Step 2: Bounding $\| \bb' - \bb \|_1$.} To show that $\bb'$ will be close $\bb$, we proceed similarly to Step 1. In particular, we define $\bb_k$ for each $k \in [K]$ as 
\begin{align*}
\bb_k \equiv \Ab \wb_{k} = \frac{x^{*0}_{j_k}}{\xi_{j_k}} \cdot \Ab \eb_{j_k} = \frac{x^{*0}_{j_k}}{\xi_{j_k}} \Ab_{j_k}.
\end{align*}
Observe that by definition of $\bb_k$, we have that the sample mean of $\bb_1, \dots, \bb_K$ is equal to $\bb'$: 
\begin{equation}
\frac{1}{K} \sum_{k=1}^K \bb_k = \frac{1}{K} \sum_{k=1}^K \Ab \wb_k = \Ab \left( \frac{1}{K} \sum_{k=1}^K \wb_k \right) = \Ab \xb' \equiv \bb'.
\end{equation}
In addition, the expected value of each $\bb_k$ can be calculated; letting $\tilde{\bb}$ denote a random variable with the same distribution as each $\bb_k$, we have
\begin{align*}
\Ebb[ \tilde{\bb} ] & = \Ab \Ebb[ \wb_k ] = \Ab \xb^{*0} = \bb.
\end{align*}
Lastly, we can bound the $\ell_{\infty}$ norm of each vector $\bb_k$ as 
\begin{equation*}
\| \bb_k \|_{\infty} = \left\| \frac{x^{*0}_{j_k}}{\xi_{j_k}} \Ab_{j_k} \right\|_{\infty} = \left| \frac{x^{*0}_{j_k}}{ \xi_{j_k} } \right| \cdot \| \Ab_{j_k} \|_{\infty} \leq C \| \Ab \|_{\max},
\end{equation*}
where the inequality follows by the definition of $\| \Ab \|_{\max}$ and the fact that $\xb^{*0}$ satisfies $\zerob \leq \xb \leq C \xib$. 

With these observations in hand, we now recognize that $\| \bb' - \bb \|_1 = \| (1/K) \sum_{k=1}^K \bb_k - \Ebb[ \tilde{\bb} ] \|_1$, i.e., $\| \bb' - \bb \|_1$ is just the $\ell_1$ norm of the deviation of a sample mean from its true expectation; we can therefore invoke Lemma~\ref{lemma:averaged_point_and_distance_L1} to assert that, with probability at least $1 - \delta/2$,
\begin{equation}
\| \bb' - \bb   \|_1 \leq \frac{m \cdot C \cdot \| \Ab \|_{\max} }{\sqrt{K}} \cdot \left(  1 + \sqrt{2 \log \frac{2}{\delta}}  \right).\label{eq:bound_sampled_mean_in_b}
\end{equation}

\textbf{Step 3: Completing the proof.} With Steps 1 and 2 complete, we are now ready to bound the optimality gap. For any vector $\bb'' \in \mathbb{R}^m$, we define the linear program $P_J(\bb'')$ as
\begin{align}
\label{problem:LP_sampled_perturb_on_b}
P_J(\bb''): \quad  \min \left\lbrace \cb^T\xb \mid \Ab \xb = \bb'', \xb \geq \zerob,\,\,\, x_j = 0\,\, \forall j \notin J \right\rbrace.
\end{align}
Then $v(P_J(\bb')) \leq \cb^T \xb'$; this follows because $\Ab \xb' = \bb'$ and $\xb' \geq \zerob$, which means that $\xb'$ is a feasible solution to problem~$P_J(\bb')$. In addition, since $ \cb^T \xb^{*0} = v(P_\distr)$, we have
\begin{align}
\label{eq:combined_from_sample_to_restriced_problem}
v(P_J(\bb')) \leq \cb^T \xb' = \cb^T \left( \xb^{*0} + (\xb' - \xb^{*0})  \right) = v(P_\distr) + \cb^T(\xb' - \xb^{*0}).
\end{align}
If the column-randomized problem $P_J$ is feasible, then by letting $\pb$ be any optimal solution of the dual of $P_J$ and applying Lemma~\ref{lemma:sensitivity_of_obj}, we have
\begin{align}
v(P_J) = v(P_J(\bb)) & \leq  v(P_J(\bb'))  + \pb^T (\bb - \bb')  \label{eq:application_of_sensitivity_analysis}\\ 
& \leq v(P_\distr) + \cb^T(\xb' - \xb^{*0})   + \pb^T (\bb - \bb') \label{eq:objective_bound_on_both_x_and_b} \\
& \leq v(P_\distr)   +  \| \cb \|_2 \cdot  \| \xb' - \xb^{*0} \| + \| \pb \|_\infty \cdot \| \bb' - \bb \|_1 \\ 
& =  v(P_\distr)   +  \| \xb' - \xb^{*0} \|_2 +  \| \pb \|_\infty \cdot \| \bb' - \bb \|_1 \label{eq:objective_bound_on_both_x_and_b_abs},
\end{align}
where the first inequality comes from Lemma~\ref{lemma:sensitivity_of_obj}, the second inequality comes from \eqref{eq:combined_from_sample_to_restriced_problem}, the third inequality comes from the Cauchy-Schwarz inequality and H\"{o}lder's inequality, and the last equality comes from the assumption that $\| \cb \|_2 = 1$. 

We now bound expression \eqref{eq:objective_bound_on_both_x_and_b_abs} by applying the inequalities~\eqref{eq:bound_sampled_mean_in_x} and \eqref{eq:bound_sampled_mean_in_b}, each of which hold with probability at least $1 - \delta/2$, and combining them using the union bound. We thus obtain that, with probability at least $1 - \delta$,
\begin{align}
v(P_J) \leq v(P_\distr)  +  \frac{C}{\sqrt{K}} \cdot \left( 1 + \| \pb \|_\infty \cdot m \cdot \Ab_{\max} \right) \cdot \left(  1 + \sqrt{2 \log \frac{2}{\delta}}  \right).
\end{align}
Subtracting $v(P)$ from both sides gives us the required inequality. \hfill \Halmos \\

With Proposition~\ref{prop:objective_bound_by_dual_solution}, we can smoothly prove Theorem~\ref{thm:main_largest_abs_of_all_dual_BFS} as follows. %

\proof{Proof of Theorem~\ref{thm:main_largest_abs_of_all_dual_BFS}:}

By invoking Proposition~\ref{prop:objective_bound_by_dual_solution}, we obtain that with probability at least $1 - \delta$, if $P_J$ is feasible, then 
\begin{align*}
\Delta v(P_J) \leq \Delta v(P_\distr)  +  \frac{C}{\sqrt{K}} \cdot \left( 1 + \| \pb \|_\infty \cdot m \cdot \Ab_{\max} \right)
\cdot \left(  1 + \sqrt{2 \log \frac{2}{\delta}}  \right),
\end{align*}
for any dual optimal solution $\pb$ of $D_J$. To prove the theorem, let us set $\pb$ to an optimal basic feasible solution of the problem $D_J$. Note that such a dual optimal solution is guaranteed to exist by the assumption that $\rank(\Ab_J) = m$. %
Since $\pb$ is a basic feasible solution of $D_J$, it is automatically a basic (but not necessarily feasible) solution of the complete dual problem $D$. By the definition of $\gamma$ in the theorem, we have that $\| \pb \|_{\infty} \leq \gamma$, and the theorem follows. \hfill \Halmos
\endproof

To prove Theorem~\ref{thm:main_reduced_cost}, we prove a complementary result to Proposition~\ref{prop:objective_bound_by_dual_solution}.

\begin{proposition}
	\label{prop:objective_bound_by_dual_slack}
	Let $C$, $P_J$ and $P_{\distr}$ be defined as in the statement of Proposition~\ref{prop:objective_bound_by_dual_solution}. For any $\delta \in (0,1)$, with probability at least $1 - \delta$ over the sample $J$, the following holds: if $P_J$ is feasible, then 
	\begin{equation*}
	\Delta v(P_J) \leq \Delta v(P_{\distr}) + \frac{C}{\sqrt{K}} \cdot  \| \cb^T - \pb^T \Ab \|_2 \cdot \left(1 + \sqrt{2 \log \frac{1}{\delta} } \right)
	\end{equation*}
	for any optimal solution $\pb$ of problem $D_J$ (the dual of problem $P_J$).
\end{proposition}

\proof{Proof:} We follow the proof of Proposition~\ref{prop:objective_bound_by_dual_solution} until inequality~\eqref{eq:objective_bound_on_both_x_and_b} and continue as follows:
\begin{equation}
\label{eq:key_step_of_proof_of_theorem_2}
\begin{aligned}
v(P_J) = v(P_J(\bb)) & \leq  v(P_J(\bb'))  + \pb^T (\bb - \bb')  \\ 
& \leq v(P_\distr) + \cb^T(\xb' - \xb^{*0})   + \pb^T (\bb - \bb') \\ 
& = v(P_\distr) + \cb^T(\xb' - \xb^{*0})   + \pb^T \Ab (\xb^{*0} - \xb') \\ 
& = v(P_\distr) + \left( \cb^T - \pb^T \Ab \right)(\xb' - \xb^{*0}) \\ 
& \leq v(P_\distr) + \| \cb^T - \pb^T \Ab \|_2 \cdot  \| \xb' - \xb^{*0}\|_2,
\end{aligned}
\end{equation}
where the bound holds for any optimal solution $\pb$ of the sampled dual problem $D_J$. By invoking Lemma~\ref{lemma:averaged_point_and_distance} with $\delta$ to bound $\| \xb' - \xb^{*0} \|_2$, and subtracting $v(P)$ from both sides, we obtain the desired result. \hfill \Halmos \\
\endproof

Using Proposition~\ref{prop:objective_bound_by_dual_slack}, we now prove Theorem~\ref{thm:main_reduced_cost}. 

\proof{Proof of Theorem~\ref{thm:main_reduced_cost}:} 
We invoke Proposition~\ref{prop:objective_bound_by_dual_slack} and set $\pb$ to be an optimal basic feasible solution of the sampled dual problem $D_J$; then $\pb^T = \cb^T_B \Ab^{-1}_B$ for some set of basic variables $B \subset [n]$. In this case, we observe that the dual slack vector $\cb^T - \pb^T \Ab$ becomes $\cb^T - \cb_B^T \Ab^{-1}_B \Ab$, which is exactly the reduced cost vector $\bar{\cb}$ associated with the basis $B$ within the full problem $P$. By using the hypothesis that any such reduced cost vector satisfies $\| \bar{\cb} \|_2 \leq \chi$, we obtain the desired result. \hfill \Halmos
\endproof

\subsection{Proof of Proposition~\ref{prop:near_feasibility_bound}}
\label{proof:near_feasibility_bound}

Let $\xb^{*0}$ be an optimal solution of $P^{\feas}_{\distr}$. Define the solution $\xb'$ as
\begin{equation*}
\xb' = \frac{1}{K} \sum_{k=1}^K \frac{x^{*0}_{j_k}}{\xi_{j_k}} \cdot \eb_{j_k}.
\end{equation*}
With $\xb'$, we can bound the objective value of $P^{\feas}_J$ as follows:
\begin{align}
v(P^{\feas}_J) & \leq \| \Ab \xb' - \bb \|_1 \nonumber \\
& = \| \Ab \xb' - \Ab \xb^{*0} + \Ab \xb^{*0} - \bb \|_ 1 \nonumber \\
& \leq \| \Ab \xb' - \Ab \xb^{*0} \|_1 + \| \Ab \xb^{*0} - \bb \|_ 1 \nonumber \\
& = \| \Ab \xb' - \Ab \xb^{*0} \|_1 + v(P^{\feas}_{\distr}) \label{bound:feasibility_gap}
\end{align}
where the first step follows by the fact that $\xb'$, when restricted to the indices in $J$, is a feasible solution of $P^{\feas}_J$; the third step follows by the triangle inequality; and the fourth follows by the definition of $\xb^{*0}$ as an optimal solution of $P^{\feas}_{\distr}$. 

The only remaining step is to bound $\| \Ab \xb' - \Ab \xb^{*0}\|_1$. To do this, let us define the vector $\vb_k$ as 
\begin{equation*}
\vb_k = \frac{ x^{*0}_{j_k} }{\xi_{j_k}} \Ab_{j_k}
\end{equation*}
for each $k \in [K]$. The vectors $\vb_1, \dots, \vb_K$ are special for three reasons. First, their sample mean is exactly
\begin{align*}
\frac{1}{K} \sum_{k=1}^K \vb_k & = \frac{1}{K} \sum_{k=1}^K \frac{ x^{*0}_{j_k} }{\xi_{j_k}} \Ab_{j_k} \\
& = \frac{1}{K} \sum_{k=1}^K \frac{ x^{*0}_{j_k} }{\xi_{j_k}} \Ab \eb_{j_k}\\
& = \Ab \xb'. 
\end{align*}
Second, letting $\vb$ denote a random variable following the same distribution as each $\vb_k$, the expected value of each $\vb_k$ is
\begin{align*}
\Ebb[\vb] & = \sum_{j \in I_+} \xi_j \cdot \frac{ x^{*0}_j}{\xi_j} \Ab_j \\
& = \sum_{j \in I_+} x^{*0}_j \Ab_j \\
& = \sum_{j \in [n]} x^{*0_j} \Ab_j \\
& = \Ab \xb^{*0}
\end{align*}
where $I_+$ is the subset of indices in $[n]$ such that $\xi_j > 0$. Note that the third step is justified by observing that $\xi^{*0}_j = 0$ whenever $j \notin I_+$ (this is because of the constraint $\zerob \leq \xb \leq C \xib$ in the definition of $P^{\feas}_{\distr}$). 

Lastly, observe that each $\vb_k$ is bounded as 
\begin{equation*}
\| \vb_k \|_{\infty} = \frac{x^{*0}_{j_k}}{\xi_{j_k}} \cdot \| \Ab_{j_k} \|_{\infty} \leq C \cdot H,
\end{equation*}
where we use the hypothesis that $\| \Ab_j \|_{\infty} \leq \| \Ab \|_{\max}$ and the fact that $\xb^{*0}$ satisfies $\zerob \leq \xb^{*0} \leq C \xib$. 

With all of these properties, the quantity $\| \Ab \xb' - \Ab \xb^{*0}\|_1$ can be re-written as $\| (1/K) \sum_{k=1}^K \vb_k - \Ebb[\vb] \|_1$, which we can bound using Lemma~\ref{lemma:averaged_point_and_distance_L1} (see Section~\ref{subsec:preliminary_results}). Invoking Lemma~\ref{lemma:averaged_point_and_distance_L1}, we get that
\begin{align*}
\| \Ab \xb' - \Ab \xb^{*0} \|_1 & = \| \frac{1}{K} \sum_{k=1}^K \vb_k - \Ebb[\vb] \|_1 \\
& \leq \frac{mC \| \Ab \|_{\max}}{\sqrt{K}} \left( 1 + \sqrt{2 \log \frac{1}{\delta}} \right).
\end{align*}
with probability at least $1 - \delta$. Using this within the bound~\eqref{bound:feasibility_gap}, we obtain that
\begin{align*}
v(P^{\feas}_J) & \leq v(P^{\feas}_{\distr}) + \| \Ab \xb' - \Ab \xb^{*0} \|_1 \\
& \leq v(P^{\feas}_{\distr}) + \frac{C}{\sqrt{K}} \cdot m \cdot \| \Ab \|_{\max} \cdot \left( 1 + \sqrt{2 \log \frac{1}{\delta}} \right)
\end{align*}
holds with probability at least $1 - \delta$, which completes the proof. \hfill \Halmos

\subsection{Proof of Theorem~\ref{thm:main_largest_abs_of_all_dual_BFS_FG} and \ref{thm:main_reduced_cost_FG}}
\label{subsec:feasbility_guarantee_proof}

As with Theorems~\ref{thm:main_largest_abs_of_all_dual_BFS} and \ref{thm:main_reduced_cost}, we first establish analogs of Propositions~\ref{prop:objective_bound_by_dual_solution} and \ref{prop:objective_bound_by_dual_slack} for Algorithm~\ref{alg:main_FG}. 

\begin{proposition}
	\label{prop:objective_bound_by_dual_solution_FG}
	Let $C$ be a nonnegative constant and let $P_{J \cup J_F}$ be the column-randomized LP solved by Algorithm~\ref{alg:main_FG}. For any $\delta \in (0,1)$, with probability at least $1 - \delta$ over the sample $J$, then 
	\begin{equation*}
	\Delta v(P_{J \cup J_F}) \leq \Delta v(P_{\distr}) + \frac{C}{\sqrt{K}} \cdot (1 + \| \pb \|_{\infty} \cdot m \cdot \| \Ab \|_{\max} ) \cdot \left(1 + \sqrt{2 \log \frac{2}{\delta} } \right)
	\end{equation*}
	for any optimal solution $\pb$ of problem $D_{J \cup J_F}$.
\end{proposition}
\proof{Proof:}
The proof of Proposition~\ref{prop:objective_bound_by_dual_solution_FG} follows along similar lines as the proof of Proposition \ref{prop:objective_bound_by_dual_solution}. More specifically, we construct $\xb'$ and $\bb'$ in the same way, and Steps 1 and 2 follow through identically. In the last step, Step 3, the sequence of bounding steps is almost the same, with a few differences:
\begin{align*}
v(P_{J \cup J_F}) & = v(P_{J \cup J_F}(\bb)) \\
& \leq v(P_{J \cup J_F}(\bb')) + \pb^T (\bb - \bb') \\
& \leq \cb^T \xb' + \pb^T (\bb - \bb') \\
& = \cb^T \xb^{*0} + \cb^T (\xb' - \xb^{*0}) + \pb^T (\bb - \bb') \\
& = v(P_{\distr}) + \cb^T (\xb' - \xb^{*0}) + \pb^T (\bb - \bb') \\
& \leq v(P_{\distr}) + \| \cb \|_2 \| \xb' - \xb^{*0} \|_2 +  \| \pb \|_{\infty}  \| \bb - \bb' \|_1.
\end{align*}
In the above, there are two important, subtle differences in the bounding. First, $\pb$ is now any optimal dual solution of $P_{J \cup J_F}$, whereas in Proposition~\ref{prop:objective_bound_by_dual_solution}, we required $\pb$ to be any optimal dual solution of $P_{J}$. Additionally, the second inequality follows because $\xb'$, which we defined as
\begin{equation*}
\xb' = \frac{1}{K} \sum_{k=1}^K \frac{ x^{*0}_{j_k}}{\xi_{j_k}} \eb_{j_k},
\end{equation*}
is still a feasible solution of $P_{J \cup J_F}(\bb')$. (Note that $\xb'$ is supported on $J$, which is obviously a subset of $J \cup J_F$; additionally, $\bb'$ was defined as $\bb' \equiv \Ab \xb'$, so by construction $\xb'$ must satisfy the equality constraint, and by construction $\xb'$ is nonnegative, so it satisfies the nonnegativity constraint.)

The remaining steps, which involve applying the high probability bounds from Steps 1 and 2 to $\| \xb' - \xb^{*0} \|_2$ and $\| \bb - \bb' \|_1$, follow in the same way as in the proof of Proposition~\ref{prop:objective_bound_by_dual_solution}. \hfill \Halmos
\endproof

We now prove Theorem~\ref{thm:main_largest_abs_of_all_dual_BFS_FG}.

\proof{Proof of Theorem~\ref{thm:main_largest_abs_of_all_dual_BFS_FG}:}
As in the proof of Theorem~\ref{thm:main_largest_abs_of_all_dual_BFS}, we invoke Proposition~\ref{prop:objective_bound_by_dual_solution_FG} with $\pb$ set to an optimal basic feasible solution of $D_{J \cup J_F}$. Note that such a solution exists because Algorithm~\ref{alg:main_FG} guarantees that $\rank(\Ab_{J \cup J_F}) = m$. Since $\pb$ is a basic solution of $D_{J \cup J_F}$ it remains a basic solution of the complete dual problem $D$, and thus it obeys $\| \pb \|_{\infty} \leq \gamma$, which establishes the theorem. \hfill \Halmos
\endproof

To establish Theorem~\ref{thm:main_reduced_cost_FG}, we similarly need an analog of Proposition~\ref{prop:objective_bound_by_dual_slack} for Algorithm~\ref{alg:main_FG}.

\begin{proposition}
	\label{prop:objective_bound_by_dual_slack_FG}
	Let $C$, $P_J$ and $P_{\distr}$ be defined as in the statement of Proposition~\ref{prop:objective_bound_by_dual_solution}. For any $\delta \in (0,1)$, with probability at least $1 - \delta$ over the sample $J$, the following holds: if $P_J$ is feasible, then 
	\begin{equation*}
	\Delta v(P_J) \leq \Delta v(P_{\distr}) + \frac{C}{\sqrt{K}} \cdot  \| \cb^T - \pb^T \Ab \|_2 \cdot \left(1 + \sqrt{2 \log \frac{1}{\delta} } \right)
	\end{equation*}
	for any optimal solution $\pb$ of problem $D_J$ (the dual of problem $P_J$).
\end{proposition}

\proof{Proof:}
As with Proposition~\ref{prop:objective_bound_by_dual_solution_FG}, we construct $\xb'$ and $\bb'$ as in the proof of Proposition~\ref{prop:objective_bound_by_dual_solution}, and follow Steps 1 and 2 from that proof. We then follow the bounding procedure in the proof of Proposition~\ref{prop:objective_bound_by_dual_slack}, with some minor modifications:
\begin{align*}
v(P_{J \cup J_F}) & = v(P_{J \cup J_F}(\bb)) \\
& \leq v(P_{J \cup J_F}(\bb')) + \pb^T (\bb - \bb') \\
& \leq \cb^T \xb' + \pb^T (\bb - \bb') \\
& = \cb^T \xb^{*0} + \cb^T (\xb' - \xb^{*0}) + \pb^T (\bb - \bb') \\
& = v(P_{\distr}) + \cb^T (\xb' - \xb^{*0}) + \pb^T (\Ab \xb - \Ab \xb') \\
& = v(P_{\distr}) + (\cb^T - \pb^T \Ab) (\xb' - \xb^{*0}) \\
& \leq v(P_{\distr}) + \| \cb^T - \pb^T \Ab \|_2 \| \xb' - \xb^{*0} \|_2, 
\end{align*}
where the main difference from the proof of Proposition~\ref{prop:objective_bound_by_dual_slack} is again that $\pb$ is a dual optimal solution of $P_{J \cup J_F}$, and we use the fact that $\xb'$, which is supported on $J$, is a feasible solution of $P_{J \cup J_F}(\bb')$. From here, the rest of the proof is the same as Proposition~\ref{prop:objective_bound_by_dual_slack}. \hfill \Halmos
\endproof

We now prove Theorem~\ref{thm:main_reduced_cost_FG}.

\proof{Proof of Theorem~\ref{thm:main_reduced_cost_FG}:}
As in the proof of Theorem~\ref{thm:main_reduced_cost}, we invoke Proposition~\ref{prop:objective_bound_by_dual_slack_FG} and set $\pb$ to be an optimal basic feasible solution of the dual problem $D_{F \cup F_J}$. Since $\pb^T = \cb^T_B \Ab_{B}^{-1}$ for some set of basic variables $B \subset J \cup J_F \subset [n]$, the dual slack vector $\cb^T - \pb^T \Ab$ is the reduced cost vector $\bar{\cb}$ of the basis $B$ within the full problem $P$, and using the assumption that any such $\bar{\cb}$ obeys $\| \bar{\cb} \|_2 \leq \chi$, the result follows. \hfill \Halmos
\endproof

\section{Omitted Proofs and Other Results for Section~\ref{sec:P_distr}}
\label{sec:proofs_distributional_counterpart}

\subsection{Proof of Theorem~\ref{theorem:many_BFS_P_distr_bound}}

For the solutions $\xb^{1}, \dots, \xb^{M}$, consider the averaged solution $\tilde{\xb}$ defined as
\begin{equation*}
\tilde{x}_j = \frac{1}{M} \sum_{i=1}^M x^i_j,
\end{equation*}
for each column $j \in [n]$.
Since each column $j$ is in at most $R$ of the bases $B^1,\dots B^M$, any coordinate $j$ of $\tilde{\xb}$ is the average of $M$ values of which at most $R$ have non-zero values, and each of those at most $R$ values is upper bounded by $x_{\max}$. It thus follows that for all $j$,
\begin{equation*}
\tilde{x}_j \leq \frac{R}{M} x_{\max}.
\end{equation*}

Observe now that by setting $C = n \cdot (R / M) x_{\max}$, we obtain that 
\begin{equation*}
C \xi_j = n \cdot (R / M) x_{\max} \cdot 1/n = (R / M) x_{\max},
\end{equation*} 
which means that $\tilde{\xb}$ satisfies the constraint $x_j \leq C \xi_j$ for all $j$. Since $\tilde{\xb}$ is the convex combination of BFSs to $P$, it satisfies $\Ab \xb = \bb$ and $\xb \geq \zerob$. Thus $\tilde{\xb}$ is a feasible solution to $P_{\distr}$. We therefore have
\begin{align*}
v(P_{\distr}) & \leq \cb^T \tilde{\xb} \\
& = \frac{1}{M} \sum_{i=1}^M \cb^T \xb^i \\
& \leq \frac{1}{M} \sum_{i=1}^M ( v(P) + \epsilon ) \\
& = v(P) + \epsilon,
\end{align*}
where the first inequality follows since $\tilde{\xb}$ is feasible for $P_{\distr}$ and the second inequality follows since each of the $M$ BFSs is assumed to be within $\epsilon$ of $v(P)$. Subtracting $v(P)$ from both sides gives the desired result. \hfill \Halmos

\subsection{Additional comments on the generative model \GMDirichletNum}
\label{subsec:additional_comment_on_generative_model_Dirichlet}

We complete our remarks on the generative model \GMDirichletNum.

\subsubsection*{An alternative interpretation}

An alternative interpretation of generative model \GMDirichletNum can also be obtained in the case when $\eta = 1$ and when $\Ab$ is structured as
\begin{equation*}
\Ab = \left[ \begin{array}{c} \Ab' \\ \oneb^T \end{array} \right],
\end{equation*}
where $\oneb$ is an $n$-dimensional vector of ones, and $\Ab'$ is a $(m-1)$-by-$n$ matrix. In this case, we can see that for any $\thetab$ in the $(n-1)$ dimensional simplex, we will have 
\begin{equation*}
\bb = \Ab \thetab = \left[ \begin{array}{c} \Ab' \thetab \\ \oneb^T \thetab \end{array} \right] = \left[ \begin{array}{c} \Ab' \thetab \\ 1 \end{array} \right],
\end{equation*}
which implies that $P$ can be written as 
\begin{subequations}
	\begin{alignat}{2}
	& \underset{\xb}{\text{minimize}} & \quad  & \cb^T \xb \\
	& \text{subject to} & & \Ab' \xb = \Ab' \thetab, \label{prob:GMDirichlet_moment_matching} \\
	& & & \oneb^T \xb  = 1, \label{prob:GMDirichlet_moment_unitsum}\\
	& & & \xb \geq \zerob. \label{prob:GMDirichlet_moment_nonnegative}
	\end{alignat}
	\label{prob:GMDirichlet_moment}%
\end{subequations}
We can think of problem~\eqref{prob:GMDirichlet_moment} as an estimation problem over the space of discrete probability distributions on $[n]$. In particular, constraints~\eqref{prob:GMDirichlet_moment_unitsum} and \eqref{prob:GMDirichlet_moment_nonnegative} enforce that $\xb$ is a probability distribution, while constraint~\eqref{prob:GMDirichlet_moment_matching} can be interpreted as a constraint that enforces a set of moments of $\xb$ to match those of $\thetab$. From this perspective, generative model \GMDirichletNum can be loosely interpreted as imposing a uniform prior. The nonparametric choice estimation problem that we numerically study in Section~\ref{sec:numerics_NCME} can be regarded as an instance of the moment problem~\eqref{prob:GMDirichlet_moment} with some modifications. %

\subsubsection*{Scaling $\cb$ so that $v(P) \geq 0$.}

A key element of generative model \GMDirichletNum is that $\cb$ is selected so that $v(P) \geq 0$. We note that this can always be accomplished: since $\xb \geq \zerob$, any nonnegative choice of $\cb$ will ensure that $v(P) \geq 0$, no matter what $\Ab$ and $\bb$ are. In addition, note that the assumption of $v(P) \geq 0$ is actually without loss of generality. If $v(P) < 0$, then let $B$ be an optimal basis, for which the corresponding reduced cost vector $\bar{\cb}$ satisfies $\bar{\cb} \geq \zerob$. (Although an optimal BFS may have negative reduced costs due to degeneracy, an optimal BFS and corresponding basis $B$ with a nonnegative reduced cost vector $\bar{\cb}$ can be obtained by applying the simplex algorithm with an anticycling pivoting rule such as Bland's rule; see Chapter 3 of \citealt{bertsimas1997introduction}.) We can then re-write $P$ as
\begin{align*}
& \min\{ \cb^T \xb \mid \Ab \xb = \bb, \xb \geq \zerob \} \\
& = \min\{ \cb_B^T \xb_B + \cb_N^T \xb_N \mid \Ab \xb = \bb, \xb \geq \zerob \} \\
& = \min\{ \cb_B^T \Ab^{-1}_B (\bb - \Ab_N \xb_N) + \cb_N^T \xb_N \mid \Ab \xb = \bb, \xb \geq \zerob \} \\
& = \cb_B^T \Ab^{-1}_B \bb + \min\{  [ \cb_N^T  - \cb_B^T \Ab^{-1}_B \Ab_N] \xb_N  \mid \Ab \xb = \bb, \xb \geq \zerob \} \\
& = v(P) + \min\{  \bar{\cb}^T \xb \mid \Ab \xb = \bb, \xb \geq \zerob \},
\end{align*}
where we observe that the problem $P' \equiv \min\{  \bar{\cb}^T \xb \mid \Ab \xb = \bb, \xb \geq \zerob \}$ is such that $P'$ and $P$ have the same feasible region and optimal solutions, and $v(P') = 0$. Thus, by replacing $\cb$ with the reduced cost vector $\bar{\cb}$ we obtain an equivalent problem, up to a constant shift. By further normalizing $\bar{\cb}$ to have unit norm, we can ensure that the last step of Algorithm~\ref{algorithm:GMDirichlet} can be accomplished.

\subsection{Proof of Theorem~\ref{theorem:GMDirichlet_gap_bound_whp}}

In this section, we establish Theorem~\ref{theorem:GMDirichlet_gap_bound_whp} for generative model \GMDirichletNum. The first key result we require is Lemma~\ref{lemma:beta_bound_on_P_distr}, which allows us to bound the gap of the distributional counterpart in terms of an upper bound $\beta$ on the minimum infinity norm attainable in the polyhedron $P$. 

\begin{lemma}
Suppose that $P$ is feasible, $v(P) \geq 0$ and $\beta \geq \min \{ \| \xb \|_{\infty} \mid \Ab \xb = \bb, \xb \geq \zerob\}$. Suppose that $\xib$ is the uniform distribution over $[n]$, i.e., $\xi_j = 1/n$ for all $j \in [n]$. If $C = n \beta$, then $P_{\distr}$ is feasible and we have that 
\begin{equation*}
\Delta v(P_{\distr}) \leq \sqrt{n} \beta.
\end{equation*}
\label{lemma:beta_bound_on_P_distr}
\end{lemma}
\proof{Proof:}
If $C = n \beta$, then 
\begin{align*}
P_{\distr} & = \min \{ \cb^T \xb \mid \Ab \xb = \bb, \zerob \leq \xb \leq C \xib \} \\
& = \min \{ \cb^T \xb \mid \Ab \xb = \bb, \zerob \leq \xb \leq n \beta \cdot (1/n) \cdot \oneb \} \\
& = \min \{ \cb^T \xb \mid \Ab \xb = \bb, \zerob \leq \xb \leq \beta \cdot \oneb \},
\end{align*}
which must be feasible; this follows by the definition of $\beta$ as an upper bound on the minimum infinity norm of any feasible solution to $P$, which itself is assumed to be feasible. Now, observe that for any feasible solution $\xb$ of $P_{\distr}$, we have
\begin{align*}
\cb^T \xb & \leq \| \cb \|_2 \cdot \| \xb \|_2 \\
& = \sqrt{ \sum_{j=1}^n x^2_j } \\
& \leq \sqrt{ \sum_{j=1}^n C^2 \xi^2_j } \\
& = C \sqrt{ \sum_{j=1}^n (1 / n)^2 } \\
& = C / \sqrt{n} \\
& = \sqrt{n} \beta,
\end{align*}
where the first inequality follows by Cauchy-Schwartz, and the second inequality by the constraint $\xb \leq \beta \cdot \oneb$. 
This implies that $v(P_{\distr}) \leq \sqrt{n} \beta$. By the assumption that $v(P) \geq 0$, we thus have that
\begin{align*}
\Delta v(P_{\distr})  = v(P_{\distr}) - v(P) \leq \sqrt{n} \beta - 0  = \sqrt{n} \beta,
\end{align*}
as required. \hfill \Halmos
\endproof

Lemma~\ref{lemma:beta_bound_on_P_distr} is a general result that is independent of the generative model chosen; we shall use it later when establishing guarantees for generative models \GMGaussianNum and \GMBernoulliCoveringNum. 

The next auxiliary result we need is a result on ordered uniform spacings.  Let $n'$ be an integer, and suppose that $X_1,\dots, X_{n'}$ are independent uniformly distributed random variables on $[0,1]$. Define $X_{0,n'} = 0$, $X_{n'+1,n'} = 1$, and define $X_{1,n'}, \dots, X_{n',n'}$ as the order statistics of $X_1,\dots, X_{n'}$. Define $\Delta_{k:n'} = X_{k,n'} - X_{k-1,n'}$ for $k = 1,\dots, n'+1$ as the \emph{(uniform) spacings} of the sample $X_1,\dots, X_{n'}$. Finally, define the \emph{ordered uniform spacings} $\Delta_{1,n'}, \dots, \Delta_{n'+1,n'}$ as the order statistics of $\Delta_{1:n'}, \dots, \Delta_{n'+1:n'}$. The following lemma is a known result on ordered uniform spacings (see \citealt{bairamov2010limit}). 

\begin{lemma}
(\citealt{bairamov2010limit}, Section 3.) For any $k = 1,\dots, n'+1$,
\begin{equation*}
\Ebb[ \Delta_{k,n'}] = \frac{1}{n'+1} \sum_{i=n'+2-k}^{n'+1} \frac{1}{i}.
\end{equation*}
\label{lemma:bairamov}
\end{lemma} 

The uniform spacings $(\Delta_{1:n'}, \dots, \Delta_{n'+1:n'})$ are useful because their joint distribution is uniform on the $n'$-dimensional unit simplex, which is identical to the $\Dirichlet(1,\dots,1)$ distribution \citep[see equation 2.1 of][]{pyke1965spacings}. The ordered uniform spacings are useful because the largest such ordered uniform spacing, $\Delta_{n'+1,n'}$, is exactly the maximum value of a $\Dirichlet(1,\dots,1)$ random vector. The expected value of this largest ordered uniform spacing will be essential to being able to obtain a high probability bound on the minimum infinity norm solution of $P$, which is the focus of our next lemma.

\begin{lemma}
Suppose that $P$ is generated according to generative model \GMDirichletNum. Let $t \geq 1$. Then, with probability at least $1 - 1/t$, we have
\begin{equation*}
\min \{ \| \xb \|_{\infty} \mid \Ab \xb = \bb, \xb \geq \zerob \} \leq \frac{ t \eta (1 + \log n)}{n}.
\end{equation*}
\label{lemma:Dirichlet_maxtheta_bound_whp}
\end{lemma}
\proof{Proof:}
Observe that by the definition of generative model \GMDirichletNum, we know that $\bb = \Ab ( \eta \thetab)$ for a $\thetab$ drawn from the $\Dirichlet(1,\dots,1)$ distribution. Since this implies that $\eta \thetab$ is a feasible solution of $P$, we immediately have
\begin{align*}
& \min\{ \| \xb \|_{\infty} \mid \Ab \xb = \bb, \xb \geq \zerob \} \\
& \leq \| \eta \thetab \|_{\infty} \\
& = \eta \max_{j \in [n]} \theta_j.
\end{align*}
For the random variable $\max_{j \in [n]} \theta_j$, we can bound its expected value as
\begin{align*}
\Ebb [ \max_{j \in [n]} \theta_j ] & = \Ebb[ \Delta_{n, n-1}] \\
& = \frac{1}{n-1 + 1} \sum_{i=n-1+2-n}^{n-1+1} \frac{1}{i} \\
& = \frac{1}{n} \sum_{i=1}^{n} \frac{1}{i} \\
& \leq \frac{1 + \log n}{n}.
\end{align*}
In the above, the steps are as follows. The first step follows because given a sample of $n - 1$ i.i.d. uniform random variables, the $n$ unordered spacings $\Delta_{1:n-1}, \dots, \Delta_{n:n-1}$ are distributed in the same way as $\thetab$ (i.e., they follow a $\Dirichlet(1,\dots,1)$ distribution). Thus, the $n$th ordered spacing $\Delta_{n,n-1}$, which is the maximum of $\Delta_{1:n-1}, \dots, \Delta_{n:n-1}$, is distributed the same way as $\max_{j \in [n]} \theta_j$. The second step follows by Lemma~\ref{lemma:bairamov}. The third step follows by algebra. The last step follows by using the bound $\sum_{i=2}^{n} \frac{1}{i} \leq \int_{1}^n \frac{1}{s} ds = \log n$. 

Using this bound on the expected value, an application of Markov's inequality implies that with probability at least $1 - 1/t$, 
\begin{equation*}
\max_{j \in [n]} \theta_j \leq \frac{ t (1 + \log n)}{n}.
\end{equation*}
Thus, with probability at least $1 - 1/t$, 
\begin{align*}
& \min\{ \| \xb \|_{\infty} \mid \Ab \xb = \bb, \xb \geq \zerob \} \\
& \leq \eta \max_{j \in [n]} \theta_j \\
& \leq \frac{ t \eta (1 + \log n)}{n},
\end{align*}
as required. \hfill \Halmos
\endproof

We now prove Theorem~\ref{theorem:GMDirichlet_gap_bound_whp}. 
\proof{Proof of Theorem~\ref{theorem:GMDirichlet_gap_bound_whp}:}
We know that $P$ is feasible, since $\xb = \eta \thetab$ is a feasible solution, and that $v(P) \geq 0$, which is just by definition of generative model~\GMDirichletNum. By Lemma~\ref{lemma:Dirichlet_maxtheta_bound_whp}, we have that $\min\{ \| \xb \|_{\infty} \mid \Ab \xb = \bb, \xb \geq \zerob \}$ is bounded by $\beta = t \eta (1 + \log n) / n$ with probability at least $1 - 1/t$. Therefore, by Lemma~\ref{lemma:beta_bound_on_P_distr}, it follows that when $C = n \beta = t \eta (1 + \log n)$, we will have that with probability at least $1 - 1/t$, that $P_{\distr}$ is feasible and the following holds:
\begin{align*}
\Delta v(P_{\distr}) & \leq \sqrt{n} \beta \\
& = \sqrt{n} \cdot \frac{t \eta (1 + \log n)}{ n } \\
& = \frac{t \eta (1 + \log n)}{ \sqrt{n} },
\end{align*}
as required. \hfill \Halmos
\endproof

\subsection{Proof of Theorem~\ref{theorem:GMGaussian_gap_bound_whp}}

To prove Theorem~\ref{theorem:GMGaussian_gap_bound_whp}, we begin with two simple results on the behaviors of the random vectors $\Ab_1,\dots, \Ab_n$. As a preview of the later results, we will need to bound the expected value of the supremum of the deviation of the sample average of $(\vb^T \Ab_1)_+,\dots, (\vb^T \Ab_n)_+$, where $\vb$ is an $m$-dimensional unit norm vector, from its expected value. To do this, we will essentially use the Rademacher complexity of the class of functions of the form $f_{\vb}(\tilde{\Ab}) = \vb^T \tilde{\Ab}$ over all unit norm vectors $\vb$. The first result, Lemma~\ref{lemma:symmetry_YjAj}, will allow us to eliminate the $(\cdot)_+$ function when we eventually bound this Rademacher complexity, while the second result, Lemma~\ref{lemma:norm_mean_YjAj_bound}, will allows us to bound the simplified expression that results from Lemma~\ref{lemma:symmetry_YjAj}.

\begin{lemma}
Suppose that $\sigma_j$ is a Rademacher variable, i.e., it takes the values -1 and +1 each with probability 1/2; $\Ab_j$ is a random vector drawn from a standard multivariate normal distribution on $\Rbb^m$; and $Y_j$ is a $\Bernoulli(1/2)$ random variable. Suppose that all three random variables are independent. Then, for any vector $\vb \in \Rbb^m$, the random variables $\sigma_j (\vb^T \Ab_j)_+$ and $\vb^T Y_j \Ab_j$ have the same distribution. 
\label{lemma:symmetry_YjAj}
\end{lemma}
\proof{Proof:}
Suppose that $t < 0$. Then we have
\begin{align*}
\Pr( \sigma_j (\vb^T \Ab_j)_+ \leq t) & = \Pr( \sigma_j (\vb^T \Ab_j)_+ \leq t \mid \sigma_j = +1) \Pr(\sigma_j = +1) + \Pr( \sigma_j (\vb^T \Ab_j)_+ \leq t \mid \sigma_j = -1) \Pr(\sigma_j = -1) \\
& = \Pr( (\vb^T \Ab_j)_+ \leq t ) \cdot (1/2) + \Pr( - (\vb^T \Ab_j)_+ \leq t) \cdot (1/2) \\
& = \Pr( (\vb^T \Ab_j)_+ \geq -t) \cdot (1/2) \\
& = \Pr( \vb^T \Ab_j \geq - t) \cdot (1/2) \\
& = \Pr( \vb^T \Ab_j \leq t) \cdot (1/2),
\end{align*}
whereas
\begin{align*}
\Pr( \vb^T Y_j \Ab_j \leq t) & = \Pr( \vb^T Y_j \Ab_j \leq t \mid Y_j = 1) \Pr( Y_j = 1) + \Pr( \vb^T Y_j \Ab_j \leq t \mid Y_j = 0) \Pr( Y_j = 0) \\
& = \Pr( \vb^T \Ab_j \leq t) \cdot (1/2) + 0 \cdot (1/2) \\
& = \Pr( \vb^T \Ab_j \leq t) \cdot (1/2).
\end{align*}
On the other hand, suppose $t \geq 0$. Then we have 
\begin{align*}
\Pr( \sigma_j (\vb^T \Ab_j)_+ \leq t) & = \Pr( \sigma_j (\vb^T \Ab_j)_+ \leq t \mid \sigma_j = +1) \Pr( \sigma_j = +1) + \Pr( \sigma_j (\vb^T \Ab_j)_+ \leq t \mid \sigma_j = -1) \Pr( \sigma_j = -1)  \\
& = \Pr( (\vb^T \Ab_j)_+ \leq t ) \cdot (1/2) + \Pr( - (\vb^T \Ab_j)_+ \leq t) \cdot (1/2) \\
& = \Pr( \vb^T \Ab_j \leq t ) (1/2) + (1) (1/2),
\end{align*}
whereas
\begin{align*}
\Pr( \vb^T Y_j \Ab_j \leq t) & = \Pr( \vb^T Y_j \Ab_j \leq t \mid Y_j = 1) \Pr( Y_j = 1) + \Pr( \vb^T Y_j \Ab_j \leq t \mid Y_j = 0) \Pr( Y_j = 0) \\
& = \Pr( \vb^T \Ab_j \leq t) (1/2) + (1) (1/2),
\end{align*}
as desired. \hfill\Halmos
\endproof

\begin{lemma}
Let $Y_1,\dots,Y_n$ be sampled independently from a $\Bernoulli(1/2)$ distribution and $\Ab_1,\dots, \Ab_n$ be sampled independently from a standard multivariate normal distribution, i.e., $\Ab_1,\dots, \Ab_n \sim \Normal( \zerob, \Ib)$. Then
\begin{equation*}
\Ebb \left\| \frac{1}{n} \sum_{j=1}^n Y_j \Ab_j \right\|_2 \leq \frac{\sqrt{m}}{\sqrt{2n}}.
\end{equation*}
\label{lemma:norm_mean_YjAj_bound}
\end{lemma}
\proof{Proof:}
We have
\begin{align*}
& \Ebb \left\| \frac{1}{n} \sum_{j=1}^n Y_j \Ab_j \right\|_2 \\
& \leq \sqrt{ \Ebb \left\| \frac{1}{n} \sum_{j=1}^n Y_j \Ab_j \right\|^2_2 } \\
& = \sqrt{ \Ebb \left[  \frac{1}{n^2} \sum_{j_1=1}^n \sum_{j_2=1}^n Y_{j_1} Y_{j_2} \Ab^T_{j_1} \Ab_{j_2} \right] } \\
& = \sqrt{ \frac{1}{n^2} \sum_{j=1}^n \Ebb[ Y_j^2] \Ebb[ \| \Ab_j \|^2_2 ]   } \\
& = \sqrt{ \frac{1}{n^2} \sum_{j=1}^n (1/2) m   } \\
& = \sqrt{ \frac{ n \cdot m}{ 2 n^2} } \\
& = \frac{ \sqrt{m} }{ \sqrt{2n} },
\end{align*}
where the first step follows by Jensen's inequality; the second comes from the definition of the squared norm of a vector as the inner product of that vector with itself; the third comes from the fact that each term $Y_{j_1} Y_{j_2} \Ab^T_{j_1} \Ab_{j_2}$ has an expected value of zero when $j_1 \neq j_2$ (since $\Ab_{j_1}$ and $\Ab_{j_2}$ both have expected value $\zerob$ and are independent), and the independence of the $Y_j$ and $\Ab_j$ variables; the fourth comes from the fact $\Ebb[ Y_j^2] = \Ebb[Y_j] = 1/2$, while $\| \Ab_j \|^2_2$ is a chi-squared random variable with degrees of freedom $m$, so $\Ebb \| \Ab_j \|^2_2 = m$; and the fifth and sixth steps follow by algebra. \hfill \Halmos
\endproof

The next result we will need is a bound on the aforementioned expected supremum of the deviation of the sample average of $(\vb^T \Ab_1)_+,\dots, (\vb^T \Ab_n)_+$ from its expected value.
\begin{lemma}
Let $\tilde{\Ab}, \Ab_1,\dots, \Ab_n$ be sampled independently from a standard multivariate normal distribution, i.e., $\tilde{\Ab}, \Ab_1,\dots, \Ab_n \sim \Normal( \zerob, \Ib)$. Then 
\begin{equation*}
\Ebb \left[ \sup_{\vb : \| \vb \|_2 = 1} \left| \Ebb (\vb^T \tilde{\Ab})_+ - \frac{1}{n} \sum_{j=1}^n ( \vb^T \Ab_j )_+ \right| \right] \leq \frac{\sqrt{2m}}{\sqrt{n}}
\end{equation*}
where $(\cdot)_+ = \max\{ \cdot, 0\}$. 
\label{lemma:exp_sup_v_bound}
\end{lemma}

\proof{Proof:}
To prove this, we will use a classical symmetrization argument from statistical learning theory (see for example \citealt{mohri2018foundations}, Theorem~3.1). Such techniques are typically used to bound an expected value of the form $\Ebb[ \sup_{f \in \Fcal} ( \frac{1}{n} \sum_{j=1}^n f( X_j) - \Ebb[ f(X) ] )]$, where $X_1,\dots, X_n, X$ are i.i.d. random variables and $\Fcal$ is a class of functions, by the Rademacher complexity of $\Fcal$, which is defined as $R( \Fcal) = \Ebb[ \sup_{f \in \Fcal} \frac{1}{n} \sum_{j=1}^n \sigma_j f(X_j) ]$, where $\sigma_1,\dots, \sigma_n$ are i.i.d. Rademacher random variables, that is, random variables that are either $+1$ or $-1$ with probability 1/2. (For our purposes, it will not be necessary to formally define the Rademacher complexity, because as we will see, our assumption that $\Ab_1,\dots, \Ab_n$ are standard Gaussian random vectors will allow us to bound it directly.)

We have:
\begin{align*}
& \Ebb_{ \{ \Ab_j\} } \left[ \sup_{\vb : \| \vb \|_2 = 1} \left| \Ebb_{\tilde{\Ab}} ( \vb^T \tilde{\Ab} )_+ - \frac{1}{n} \sum_{j=1}^n ( \vb^T \Ab_j )_+ \right| \right] \\
& = \Ebb_{ \{ \Ab_j\} } \left[ \sup_{\vb : \| \vb \|_2 = 1} \left| \Ebb_{ \{ \tilde{\Ab}_j \} }[ \frac{1}{n} \sum_{j=1}^n (\vb^T \tilde{\Ab}_j)_+ ] - \frac{1}{n} \sum_{j=1}^n (\vb^T \Ab_j )_+ \right| \right] \\
& \leq \Ebb_{ \{ \Ab_j\} } \left[ \sup_{\vb : \| \vb \|_2 = 1} \Ebb_{ \{ \tilde{\Ab}_j \} } \left|  \frac{1}{n} \sum_{j=1}^n (\vb^T \tilde{\Ab}_j)_+ - \frac{1}{n} \sum_{j=1}^n (\vb^T \Ab_j)_+ \right| \right] \\
& \leq \Ebb_{ \{ \Ab_j\}, \{ \tilde{\Ab}_j \}}  \sup_{\vb : \| \vb \|_2 = 1} \left|  \frac{1}{n} \sum_{j=1}^n ( (\vb^T \tilde{\Ab}_j)_+  -  ( \vb^T \Ab_j)_+ ) \right|  \\
& = \Ebb_{ \{ \Ab_j\}, \{ \tilde{\Ab}_j \}, \{ \sigma_j \} }  \sup_{\vb : \| \vb \|_2 = 1} \left|  \frac{1}{n} \sum_{j=1}^n \sigma_j ( (\vb^T \tilde{\Ab}_j)_+ - (\vb^T \Ab_j)_+ ) \right|  \\
& \leq \Ebb_{ \{ \Ab_j\}, \{ \tilde{\Ab}_j \}, \{ \sigma_j \} } \sup_{\vb : \| \vb \|_2 = 1} \left\{  \left|  \frac{1}{n} \sum_{j=1}^n \sigma_j (\vb^T \tilde{\Ab}_j)_+ \right|  +  \left|  \frac{1}{n} \sum_{j=1}^n \sigma_j (\vb^T \Ab_j)_+  \right| \right\}  \\
& \leq \Ebb_{ \{ \tilde{\Ab}_j \}, \{ \sigma_j \} } \sup_{\vb : \| \vb \|_2 = 1} \left|  \frac{1}{n} \sum_{j=1}^n \sigma_j  (\vb^T \tilde{\Ab}_j)_+ \right| + \Ebb_{ \{ \Ab_j\}, \{ \sigma_j \} } \sup_{\vb : \| \vb \|_2 = 1} \left|  \frac{1}{n} \sum_{j=1}^n \sigma_j  (\vb^T \Ab_j)_+ \right| \\
& = 2 \Ebb_{ \{ \Ab_j\}, \{ \sigma_j \} } \sup_{\vb : \| \vb \|_2 = 1} \left|  \frac{1}{n} \sum_{j=1}^n \sigma_j  (\vb^T \Ab_j)_+ \right| \\
&=  2 \Ebb_{ \{ \Ab_j\}, \{Y_j\}  } \sup_{\vb : \| \vb \|_2 = 1} \left|  \frac{1}{n} \sum_{j=1}^n \vb^T Y_j \Ab_j \right| \\
& \leq 2 \Ebb_{ \{ \Ab_j\}, \{ Y_j \}  } \sup_{\vb : \| \vb \|_2 = 1}  \left\{ \| \vb \|_2  \cdot  \left\| \frac{1}{n} \sum_{j=1}^n Y_j \Ab_j \right\|_2 \right\} \\ 
& = 2 \Ebb \left\| \frac{1}{n} \sum_{j=1}^n Y_j \Ab_j \right\|_2 \\
& \leq 2 \frac{ \sqrt{m}}{\sqrt{2n}} \\
& = \frac{\sqrt{2m}}{\sqrt{n}}
\end{align*}
In the above derivation, the steps are as follows. The first step follows by introducing another i.i.d. sample of columns, $\tilde{\Ab}_1,\dots, \tilde{\Ab}_n$, which follow the same standard multivariate normal distribution as $\Ab_1,\dots, \Ab_n$. The second step follows by Jensen's inequality. The third step follows by elementary properties of $\sup$ and expectation. The fourth step follows by observing that the random variables $\Ab_1,\dots, \Ab_n, \tilde{\Ab}_1,\dots, \tilde{\Ab}_n$ are exchangeable, and so multiplying the difference $( |\vb^T \tilde{\Ab}_j| - |\vb^T \Ab_j| )$ by $\sigma_j$, which is equally likely to be $+1$ (leaving the term unchanged) or $-1$ (flipping the difference), will leave the overall expectation unchanged. The fifth follows by the triangle inequality and elementary properties of $\sup$, and the sixth by linearity of expectation. The seventh follows by observing that the two expectations in the prior step are identical. 

From here, the remaining steps rely on the properties of the distribution of $\Ab_1,\dots,\Ab_n$. In particular, the eighth step follows by applying Lemma~\ref{lemma:symmetry_YjAj} to assert that the random variable $\sigma_j ( \vb^T \Ab_j )_+$ is identically distributed to $\vb^T Y_j \Ab_j$, allowing us to replace the former random variable with the latter random variable and leave the expectation unchanged. The ninth step follows by the Cauchy-Schwartz inequality, and the tenth step by the fact that each $\vb$ is unit norm. The tenth and eleventh steps follows by applying Lemma~\ref{lemma:norm_mean_YjAj_bound} and algebra. \hfill \Halmos

\endproof

The last auxiliary result we will need is to characterize in closed form the expected value of $\Ebb( \pb^T \tilde{\Ab})_+$, where $\tilde{\Ab}$ is a standard normal random vector (i.e., a column of the matrix $\Ab$). 

\begin{lemma}
Suppose that $\tilde{\Ab} \sim \Normal(\zerob, \Ib)$. Then for any vector $\pb \in \Rbb^m$, 
\begin{equation*}
\Ebb (\pb^T \tilde{\Ab})_+ = \frac{\sqrt{2}}{2 \sqrt{\pi}} \| \pb \|_2. 
\end{equation*}
\label{lemma:closed_form_pbTAb_plus}
\end{lemma}
\proof{Proof:}
We have
\begin{align*}
\Ebb[ (\pb^T \tilde{\Ab})_+ ] & = \Ebb[ (\pb^T \tilde{\Ab})_+ \mid \pb^T \tilde{\Ab} < 0] \Pr(\pb^T \tilde{\Ab} < 0) + \Ebb[ (\pb^T \tilde{\Ab})_+ \mid \pb^T \tilde{\Ab} \geq 0] \Pr(\pb^T \tilde{\Ab} \geq 0) \\
& = (0) \cdot (1/2) + \Ebb[ (\pb^T \tilde{\Ab})_+ \mid \pb^T \tilde{\Ab} \geq 0]  \cdot (1/2) \\
& = \frac{\sqrt{2}}{\sqrt{\pi}} \sqrt{ \pb^T \Ib \pb} \cdot (1/2) \\
& = \frac{\sqrt{2}}{2 \sqrt{\pi}} \| \pb \|_2,
\end{align*}
where the first step follows by conditioning; the second step follows by the fact that $\pb^T \tilde{\Ab}$ follows a normal distribution with mean 0, and that $(\pb^T \tilde{\Ab})_+ = 0$ when $\pb^T \tilde{\Ab} < 0$; the third step follows by recognizing that the random variable $\left( (\pb^T \tilde{\Ab})_+ \mid \pb^T \tilde{\Ab} \geq 0 \right)$ follows the same distribution as $| \pb^T \tilde{\Ab} |$, and $|\pb^T \tilde{\Ab}|$ follows a \emph{half-normal} distribution, whose mean is $\sigma \sqrt{2} / \sqrt{\pi}$, where $\sigma$ is the standard deviation of $\pb^T \tilde{\Ab}$; and the final step follows by algebra. \hfill \Halmos
\endproof

With these auxiliary results in hand, we can now establish the following major result, which provides a high probability bound on the minimum infinity norm of any feasible solution of $P$. 

\begin{theorem}
Suppose that $P$ is generated according to generative model \GMGaussianNum. Let $t \geq 1$ and suppose that $n > 4 \pi t^2 m$. With probability at least $1 - 1/t$, we have that $P$ is feasible and that
\begin{equation*}
\min \{ \| \xb \|_{\infty} \mid \Ab \xb = \bb, \xb \geq \zerob\} < \frac{ \| \bb \|_2 }{n} \cdot \frac{1}{ \frac{\sqrt{2}}{2\sqrt{\pi}} - \frac{ t \sqrt{2m}}{\sqrt{n}}}.
\end{equation*}
\label{theorem:GMGaussian_infinity_norm_bound_whp}
\end{theorem}

\proof{Proof:}
Let $\tilde{\Ab}$ be a random vector that follows the same distribution as $\Ab_1,\dots, \Ab_n$. 
We have that 
\begin{align*}
& \min \left\lbrace \| \xb \|_{\infty} \mid \Ab \xb = \bb, \xb \geq \zerob \right\rbrace \\
& = \max \left\lbrace \pb^T \bb \mid \sum_{j=1}^n (\pb^T \Ab_j)_+ \leq 1 \right\rbrace \\
& = \max \left\lbrace \pb^T \bb \mid  \frac{1}{n} \sum_{j=1}^n (\pb^T \Ab_j)_+  \leq \frac{1}{n} \right\rbrace \\
& = \max \left\lbrace \pb^T \bb \mid  \Ebb ( \pb^T \tilde{\Ab} )_+ - \Ebb ( \pb^T \tilde{\Ab} )_+ + \frac{1}{n} \sum_{j=1}^n ( \pb^T \Ab_j )_+ \leq \frac{1}{n} \right\rbrace \\
& = \max \left\lbrace \pb^T \bb \mid  \Ebb ( \pb^T \tilde{\Ab} )_+ \leq \frac{1}{n} + \Ebb ( \pb^T \tilde{\Ab} )_+ - \frac{1}{n} \sum_{j=1}^n ( \pb^T \Ab_j )_+ \right\rbrace \\
& \leq \max \left\lbrace \pb^T \bb \mid  \Ebb ( \pb^T \tilde{\Ab} )_+ \leq \frac{1}{n} + \| \pb \|_2 \cdot \sup_{\vb : \| \vb \|_2 = 1} \left| \Ebb ( \vb^T \tilde{\Ab} )_+ - \frac{1}{n} \sum_{j=1}^n ( \vb^T \Ab_j )_+ \right|   \right\rbrace \quad (*)
\end{align*}
where the first step follows by strong duality; the second, third and fourth step by algebra; and the sixth step by recognizing that 
\begin{align*}
& \Ebb ( \pb^T \tilde{\Ab} )_+ - \frac{1}{n} \sum_{j=1}^n ( \pb^T \Ab_j )_+ \\
& \leq \| \pb \|_2 \cdot \sup_{\vb : \| \vb \|_2 = 1} \left| \Ebb ( \vb^T \tilde{\Ab} )_+ - \frac{1}{n} \sum_{j=1}^n ( \vb^T \Ab_j )_+ \right|  
\end{align*}
holds trivially when $\pb = \zerob$, whereas when $\pb \neq \zerob$, we have
\begin{align*}
& \Ebb ( \pb^T \tilde{\Ab} )_+ - \frac{1}{n} \sum_{j=1}^n ( \pb^T \Ab_j )_+ \\
& = \| \pb \|_2 \cdot \left[ \Ebb \left( \left( \frac{\pb}{\| \pb \|_2 } \right)^T \tilde{\Ab} \right)_+ - \frac{1}{n} \sum_{j=1}^n \left( \left( \frac{\pb}{\| \pb \|_2 } \right)^T \Ab_j  \right)_+ \right] \\
& \leq \| \pb \|_2 \cdot \sup_{\vb : \| \vb \|_2 = 1} \left| \Ebb ( \vb^T \tilde{\Ab} )_+ - \frac{1}{n} \sum_{j=1}^n ( \vb^T \Ab_j )_+ \right|,
\end{align*}
where the equality follows because $(\cdot)_+$ is positively homogenous, and the inequality follows because $\pb / \| \pb \|_2$ is a unit norm vector. 

To proceed from here, we will now use Lemma~\ref{lemma:exp_sup_v_bound}. Recall by Lemma~\ref{lemma:exp_sup_v_bound} that 
\begin{equation*}
\Ebb \sup_{\vb : \| \vb \|_2 = 1} \left| \Ebb(\vb^T \tilde{\Ab})_+ - \frac{1}{n} \sum_{j=1}^n (\vb^T \Ab_j)_+ \right|  \leq \frac{\sqrt{2m}}{\sqrt{n}}.
\end{equation*}
Consider the event $E$ defined as 
\begin{equation*}
E = \left\{ \sup_{\vb : \| \vb \|_2 = 1} | \Ebb(\vb^T \tilde{\Ab})_+ - \frac{1}{n} \sum_{j=1}^n (\vb^T \Ab_j)_+ | \leq \frac{t \sqrt{2m}}{\sqrt{n}} \right\}.
\end{equation*}
By Markov's inequality, we have that $\Pr( E ) \geq 1 - 1/t$. Thus, with probability $1 - 1/t$, $(*)$ is bounded from above as 
\begin{align*}
(*) & \leq \max \{ \pb^T \bb \mid  \Ebb ( \pb^T \tilde{\Ab} )_+ \leq \frac{1}{n} + \| \pb \|_2 \cdot \frac{ t \sqrt{2m}}{\sqrt{n}}  \}.
\end{align*}
We now have 
\begin{align*}
& \max \{ \pb^T \bb \mid  \Ebb ( \pb^T \tilde{\Ab} )_+ \leq \frac{1}{n} + \| \pb \|_2 \cdot \frac{ t \sqrt{2m}}{\sqrt{n}}  \} \\
& = \max \{ \pb^T \bb \mid  \frac{\sqrt{2}}{2\sqrt{\pi}} \cdot \| \pb \|_2 \leq \frac{1}{n} + \| \pb \|_2 \cdot \frac{ t \sqrt{2m}}{\sqrt{n}}  \} \\
& = \max \{ \pb^T \bb \mid  \left( \frac{\sqrt{2}}{2\sqrt{\pi}} - \frac{ t \sqrt{2m}}{\sqrt{n}} \right) \cdot \| \pb \|_2 \leq \frac{1}{n}   \} \\
& = \max \{ \pb^T \bb \mid   \| \pb \|_2 \leq \frac{1}{n} \cdot \frac{1}{\frac{\sqrt{2}}{2\sqrt{\pi}} - \frac{ t \sqrt{2m}}{\sqrt{n}}}   \} \\
& = \frac{ \| \bb \|_2 }{ n } \cdot \frac{1}{\frac{\sqrt{2}}{2\sqrt{\pi}} - \frac{ t \sqrt{2m}}{\sqrt{n}}} \\
\end{align*}
where the first step follows by the closed form expression for $\Ebb ( \pb^T \tilde{\Ab} )_+$ from Lemma~\ref{lemma:closed_form_pbTAb_plus}; the second step follows by algebra; the third step follows by algebra and also by our assumption on $n$; and the final step by the fact that $\max\{ \db^T \xb \mid \| \xb \|_2 \leq r\} = r \| \db \|_2$. Note that in the third step, we are using the hypothesis that $n > 4 \pi t^2 m$ to ensure that the coefficient $\left( \frac{\sqrt{2}}{2\sqrt{\pi}} - \frac{ t \sqrt{2m}}{\sqrt{n}} \right)$ is positive, and that the direction of the inequality in the constraint is unchanged.

We thus have, that with probability at least $1 - 1/t$, that 
\begin{equation*}
\min \{ \| \xb \|_{\infty} \mid \Ab \xb = \bb, \xb \geq \zerob \} \leq \frac{ \| \bb \|_2 }{ n } \cdot \frac{1}{\frac{\sqrt{2}}{2\sqrt{\pi}} - \frac{ t \sqrt{2m}}{\sqrt{n}}},
\end{equation*}
which establishes the required bound on the infinity norm. 

To see why $P$ must be feasible, that is, why $\{ \xb \mid \Ab \xb = \bb, \xb \geq \zerob\}$ is non-empty, observe that $P$ is feasible if and only if $\min \{ \| \xb \|_{\infty} \mid \Ab \xb = \bb, \xb \geq \zerob \}$ is feasible. Observe that this latter problem is feasible if and only if its dual problem $\max \{ \pb^T \bb \mid \sum_{j=1}^n (\pb^T \Ab_j)_+ \leq 1 \}$, which is always feasible, is bounded. By our reasoning above, this problem is bounded with probability at least $1 - 1/t$, and therefore $P \equiv \min \{ \cb^T \xb \mid \Ab \xb = \bb, \xb \geq \zerob\}$ is feasible with probability at least $1 - 1/t$, as required. \hfill \Halmos
\endproof

We can now prove Theorem~\ref{theorem:GMGaussian_gap_bound_whp}. 
\proof{Proof of Theorem~\ref{theorem:GMGaussian_gap_bound_whp}:}
Let $\beta$ be defined as 
\begin{equation*}
\beta = \frac{ \| \bb \|_2 }{ n } \cdot \frac{1}{\frac{\sqrt{2}}{2\sqrt{\pi}} - \frac{ t \sqrt{2m}}{\sqrt{n}}}.
\end{equation*}
Observe that the given $C$ in the statement of Theorem~\ref{theorem:GMGaussian_gap_bound_whp} is exactly $C = n \beta$.

By Theorem~\ref{theorem:GMGaussian_infinity_norm_bound_whp}, with probability at least $1 - 1/t$, it follows that $P$ is feasible, and that
\begin{align*}
\min \{ \| \xb \|_{\infty} \mid \Ab \xb = \bb, \xb \geq \zerob \} \leq \beta.
\end{align*}
Recall by the definition of generative model~\GMGaussianNum that whenever $P$ is feasible, we set $\cb$ so that $v(P) \geq 0$. Thus, with probability at least $1 - 1/t$, by Lemma~\ref{lemma:beta_bound_on_P_distr}, it follows that $P_{\distr}$ is feasible and that
\begin{align*}
\Delta v(P_{\distr}) & \leq \frac{C}{\sqrt{n}} \\
& = \sqrt{n} \beta \\
& = \frac{ \| \bb \|_2 }{ \sqrt{n} } \cdot \frac{1}{\frac{\sqrt{2}}{2\sqrt{\pi}} - \frac{ t \sqrt{2m} }{\sqrt{n}}},
\end{align*}
exactly as required. \hfill \Halmos
\endproof

\subsection{Generative model \GMUniformSphereNum: uniform distribution on the unit sphere}
\label{subsec:proofs_GMUniformSphere}

In this section, we discuss an alternate generative model, generative model~\GMUniformSphereNum, which is closely related to generative model~\GMGaussianNum. In this new generative model, the columns $\Ab_1,\dots, \Ab_n$ are independently randomly generated from the uniform distribution on the unit sphere $S^m = \{ \vb \in \Rbb^m \mid \| \vb \|_2 = 1\}$ in $\Rbb^m$. Upon generating the columns, we then fix $\bb$ and $\cb$. 

\begin{algorithm}
\caption{Generative Model \GMUniformSphereNum}
\begin{algorithmic}[1]
	\STATE Generate $n$ i.i.d. random vectors $\Ab_1,\dots,\Ab_n \sim \Uniform(S^m)$, where $\Uniform(S^m)$ denotes the uniform distribution on the unit sphere $S^m = \{ \vb \in \Rbb^m \mid \| \vb \|_2 = 1\}$. 
	\STATE Set $\Ab = [\Ab_1 \ \Ab_2 \ \cdots \Ab_n]$. 
	\STATE Fix any right-hand side vector $\bb \in \Rbb^m$.
	\STATE If $\{ \xb \mid \Ab \xb = \bb, \xb \geq \zerob\}$ is nonempty, fix any $\cb \in \{ \vb \in \Rbb^n \mid \| \vb \|_2 = 1 \}$ such that $v(P) \geq 0$; otherwise, fix any $\cb \in \{ \vb \in \Rbb^n \mid \| \vb \|_2 = 1 \}$. 
	\RETURN $(\Ab, \bb, \cb)$.
\end{algorithmic}
\label{algorithm:GMUniformSphere}
\end{algorithm}

This model is closely related to generative model~\GMGaussianNum, because for any standard normal random vector $\Ab'$, the random vector $\Ab' / \| \Ab' \|_2$ is uniformly distributed on the unit sphere $S^m$. Additionally, as mentioned in Section~\ref{subsec:P_distr_analysis}, this model is universal, in the sense that any LP of the form $\min \{ \cb^T \xb \mid \Ab \xb = \bb, \xb \geq \zerob \}$ with non-zero columns can be transformed into an equivalent LP where all the columns have unit norm. 

The main theoretical result of this section is Theorem~\ref{theorem:GMUniformSphere_gap_bound_whp}, which asserts that with high probability, the distributional counterpart gap under this generative model is $O(1/\sqrt{n})$. In the statement of the theorem below, $\Gamma$ is the gamma function, i.e., $\Gamma(z) = \int^{\infty}_{0} t^{z-1} e^{-t} \, dt$.

\begin{theorem}
Suppose that $P$ is generated according to generative model \GMUniformSphereNum. Assume that $\xib$ is the uniform distribution over $[n]$, that is, $\xi_j = 1 /n$ for all $j \in [n]$. Let $t \geq 1$, and suppose that $n > 4 \pi t^2 \mu_m^2$, where $\mu_m = \sqrt{2} \cdot \Gamma((m+1)/2)/\Gamma(m/2)$. Suppose that $C$ is set as 
\begin{equation*}
C = \| \bb \|_2  \cdot \frac{1}{ \frac{\sqrt{2}}{2\sqrt{\pi} \mu_m} - \frac{ t \sqrt{2}}{\sqrt{n}}}.
\end{equation*}
Then, with probability at least $1 - 1/t$, we have that $P$ and $P_{\distr}$ are feasible, and
\begin{equation*}
\Delta v(P_{\distr}) \leq \frac{ \| \bb \|_2 }{ \sqrt{n} } \cdot \frac{1}{ \frac{\sqrt{2}}{2\sqrt{\pi} \mu_m} - \frac{ t \sqrt{2}}{\sqrt{n}}}.
\end{equation*}
\label{theorem:GMUniformSphere_gap_bound_whp}
\end{theorem}
This result is very similar to Theorem~\ref{theorem:GMGaussian_gap_bound_whp}. The main difference is in the factor which multiplies $\| \bb \|_2$ in the particular choice of $C$ and $\| \bb \|_2 / \sqrt{n}$ in the upper bound on $\Delta v(P_{\distr})$. In the standard Gaussian case, this factor is $[ \sqrt{2} / (2\sqrt{\pi}) - t \sqrt{2m} / \sqrt{n}]^{-1}$, whereas in the case of the uniform distribution on the unit sphere, the factor is $[ \sqrt{2}/(2\sqrt{\pi} \mu_m) - t \sqrt{2} / \sqrt{n} ]^{-1}$. Note that $\mu_m$ appears because this is the mean of a chi distributed random variable with $m$ degrees of freedom; this distribution, in turn, appears because this is the distribution of the norm of a standard normal random vector. By Jensen's inequality, $\mu_m$ is lower than $\sqrt{m}$, which is the square root of the mean of a chi-squared distributed random variable with $m$ degrees of freedom, but numerically $\mu_m$ is actually very close to $\sqrt{m}$. Thus, comparing the two factors, the factor for generative model \GMUniformSphereNum is roughly $\sqrt{m}$ larger. This makes sense, because in generative model \GMGaussianNum, the columns will have norm that is on average larger by a factor of $\mu_m \approx \sqrt{m}$ than the columns in generative model \GMUniformSphereNum, so the decision variable vector $\xb$ should be correspondingly scaled by $\mu_m$ to ensure $\Ab \xb = \bb$.

We now turn our attention to proving Theorem~\ref{theorem:GMUniformSphere_gap_bound_whp}. As in the case of generative model~\GMGaussianNum, we will require a number of auxiliary results. Our first such auxiliary result is an analog of Lemma~\ref{lemma:symmetry_YjAj}, which is a technical result needed to bound the expected supremum of the deviation of the sample average of $(\pb^T \Ab_j)_+$ from its expected value. The proof of this result follows along very similar lines to the proof of Lemma~\ref{lemma:symmetry_YjAj}. The key is that like in the case where $\Ab_j$ is a standard normal random vector, when $\Ab_j$ is uniformly distributed on the unit sphere, the distribution of $\vb^T \Ab_j$ is symmetric about zero. For brevity, we omit the proof. 

\begin{lemma}
Suppose that $\Ab_j \sim \Uniform(S^m)$; suppose that $\sigma_j$ is a Rademacher random variable (i.e., $\sigma_j$ is either +1 or -1, both with probability 1/2); and suppose that $Y_j \sim \Bernoulli(1/2)$. Suppose that $\Ab_j$, $Y_j$ and $\sigma_j$ are independent. Then for any vector $\vb \in \Rbb^m$, the random variables $\sigma_j (\vb^T \Ab_j)_+$ and $\vb^T Y_j \Ab_j$ follow the same distribution.
\label{lemma:symmetry_YjAj_UniformSphere}
\end{lemma}

We next have an analog of Lemma~\ref{lemma:norm_mean_YjAj_bound}, which bounds the expected Euclidean norm of the sample average of a particular collection of i.i.d. random vectors. The proof of this lemma is omitted as it follows along essentially the same lines as the proof of Lemma~\ref{lemma:norm_mean_YjAj_bound}. 
\begin{lemma}
Suppose that $\Ab_1,\dots, \Ab_n \sim \Uniform(S^m)$ are independent random variables and that $Y_1,\dots, Y_n \sim \Bernoulli(1/2)$ are independent random variables. Then we have
\begin{equation*}
\Ebb \left\| \frac{1}{n} \sum_{j=1}^n Y_j \Ab_j \right\|_2 \leq \frac{1}{\sqrt{2n}}.
\end{equation*}
\label{lemma:norm_mean_YjAj_bound_UniformSphere}
\end{lemma}

With Lemma~\ref{lemma:symmetry_YjAj_UniformSphere} and \ref{lemma:norm_mean_YjAj_bound_UniformSphere} in hand, we can prove the following lemma, which is an analog of Lemma~\ref{lemma:exp_sup_v_bound}. This lemma allows us to bound the expected supremum of the deviation of the sample average of $(\vb^T \Ab_1)_+, \dots, (\vb^T \Ab_n)_+$ from its expected value, over all unit vectors $\vb$. %

\begin{lemma}
Suppose that $\tilde{\Ab}, \Ab_1,\dots, \Ab_n \sim \Uniform(S^m)$ are independent random variables distributed uniformly on the unit sphere $S^m$. Then 
\begin{equation*}
\Ebb \sup_{\vb : \| \vb \|_2 = 1} \left| \frac{1}{n} \sum_{j=1}^n (\pb^T \Ab_j)_+ - \Ebb (\pb^T \tilde{\Ab})_+ \right| \leq \frac{\sqrt{2}}{\sqrt{n}}
\end{equation*}
\label{lemma:exp_sup_v_bound_UniformSphere}
\end{lemma}
\proof{Proof:}
The proof follows essentially the same initial steps as that of Lemma~\ref{lemma:exp_sup_v_bound}. Following those steps, we obtain the bound 
\begin{align*}
& \Ebb \sup_{\vb : \| \vb \|_2 = 1} \left| \frac{1}{n} \sum_{j=1}^n (\pb^T \Ab_j)_+ - \Ebb (\pb^T \tilde{\Ab})_+ \right| \\
& \leq 2 \Ebb \sup_{\vb : \| \vb \|_2 = 1} \left| \frac{1}{n} \sum_{j=1}^n \sigma_j (\vb^T \Ab_j)_+ \right| \\
& = 2 \Ebb \sup_{\vb : \| \vb \|_2 = 1} \left| \frac{1}{n} \sum_{j=1}^n \vb^T Y_j \Ab_j \right| \\
& \leq  2 \Ebb \left[ \sup_{\vb : \| \vb \|_2 = 1} \| \vb \|_2  \left\| \frac{1}{n} \sum_{j=1}^n Y_j \Ab_j \right\|_2 \right] \\ 
& \leq  2 \Ebb  \left\| \frac{1}{n} \sum_{j=1}^n Y_j \Ab_j \right\|_2 \\ 
& \leq 2 \cdot \frac{1}{\sqrt{2n}} \\
& = \frac{\sqrt{2}}{\sqrt{n}}
\end{align*}
where the first equality follows by Lemma~\ref{lemma:symmetry_YjAj_UniformSphere}; the second inequality by Cauchy-Schwartz; the third inequality by the fact that each $\vb$ in the sup is of unit norm; and the fourth inequality by Lemma~\ref{lemma:norm_mean_YjAj_bound_UniformSphere}. \hfill \Halmos
\endproof

The last auxiliary result we will need is the following lemma, which provides the closed form expression for $\Ebb (\pb^T \tilde{\Ab})_+$ when $\tilde{\Ab}$ is uniformly distributed on the unit sphere.
\begin{lemma}
Suppose that $\tilde{\Ab} \sim \Uniform(S^m)$. Then for any vector $\pb \in \Rbb^m$,
\begin{equation*}
\Ebb (\pb^T \tilde{\Ab} )_+ = \frac{\sqrt{2}}{2 \sqrt{\pi} \mu_m} \cdot \| \pb \|_2.
\end{equation*}
where $\mu_m = \sqrt{2} \Gamma( (m+1)/2) / \Gamma( m/2)$, where $\Gamma(\cdot)$ is the gamma function. 
\label{lemma:closed_form_pbTAb_plus_UniformSphere}
\end{lemma}
\proof{Proof:}
Let $\Zb \sim \Normal(\zerob,\Ib)$ be a standard normal random vector in $\Rbb^m$. Let $\theta$ be a random variable that follows the chi distribution with $m$ degrees of freedom, and suppose that $\theta$ is independent of $\tilde{\Ab}$. Then $\Zb$ and $\theta \tilde{\Ab}$ have the same distribution. We therefore have 
\begin{align*}
\Ebb (\pb^T \Zb)_+ & = \Ebb (\pb^T \theta \tilde{\Ab})_+ \\
& = \Ebb[ \theta \cdot (\pb^T \tilde{\Ab})_+ ] \\
& = \Ebb[ \theta ] \cdot \Ebb (\pb^T \tilde{\Ab})_+ \\
& = \mu_m \Ebb (\pb^T \tilde{\Ab})_+ 
\end{align*}
where the first step follows by the distributional equivalence of $\Zb$ and $\theta \tilde{\Ab}$; the second by the fact that the $(\cdot)_+$ function is positively homogeneous; the third by the independence of $\theta$ and $\tilde{\Ab}$; and the fourth by the fact that $\mu_m$ is precisely the mean of a chi-distributed random variable. Using the fact (Lemma~\ref{lemma:closed_form_pbTAb_plus}) that 
\begin{equation*}
\Ebb (\pb^T \tilde{\Zb}_+) = \frac{\sqrt{2}}{2 \sqrt{\pi}} \| \pb \|_2,
\end{equation*}
we obtain that 
\begin{align*}
\Ebb (\pb^T \tilde{\Ab})_+ & = \frac{1}{\mu_m} \cdot \Ebb (\pb^T \Zb)_+ \\
& = \frac{\sqrt{2}}{2 \sqrt{\pi} \mu_m} \| \pb \|_2,
\end{align*}
as required. \hfill \Halmos
\endproof

We are now in a position to prove Theorem~\ref{theorem:GMUniformSphere_infinity_norm_bound_whp}, which is an analog of Theorem~\ref{theorem:GMGaussian_infinity_norm_bound_whp}. 
\begin{theorem}
Suppose that $P$ is generated according to generative model \GMUniformSphereNum. 
Let $t \geq 1$. Suppose that $n > 4 t^2 \pi \mu^2_m$, where $\mu_m = \sqrt{2} \Gamma( (m+1)/2) / \Gamma(m/2)$ is the mean of a chi-distributed random variable with $m$ degrees of freedom. Then, with probability at least $1 - 1/t$, we have that $P$ is feasible and
\begin{equation*}
\min\{ \| \xb \|_{\infty} \mid \Ab \xb = \bb, \xb \geq \zerob \} \leq \frac{ \| \bb \|_2 }{n} \cdot \frac{ 1 }{ \frac{\sqrt{2}}{2 \sqrt{\pi} \mu_m} - \frac{t \sqrt{2}}{\sqrt{n}} }.
\end{equation*}
\label{theorem:GMUniformSphere_infinity_norm_bound_whp}
\end{theorem}
\proof{Proof:}
Following similar steps as in the proof of Theorem~\ref{theorem:GMGaussian_infinity_norm_bound_whp}, we have
\begin{align*}
& \min\{ \| \xb \|_{\infty} \mid \Ab \xb = \bb, \xb \geq \zerob \} \\
& = \max \left\{ \pb^T \bb \mid \sum_{j=1}^n (\pb^T \Ab_j)_+ \leq 1 \right\} \\
& = \max \left\{ \pb^T \bb \mid \Ebb (\pb^T \tilde{\Ab})_+ \leq \frac{1}{n} + \Ebb (\pb^T \tilde{\Ab})_+ - \frac{1}{n} \sum_{j=1}^n (\pb^T \Ab_j)_+ \right\} \\
& \leq \max \left\{ \pb^T \bb \mid \Ebb (\pb^T \tilde{\Ab})_+ \leq \frac{1}{n} + \| \pb \|_2 \cdot \sup_{\vb : \| \vb \|_2 = 1} \left| \Ebb (\vb^T \tilde{\Ab})_+ - \frac{1}{n} \sum_{j=1}^n (\vb^T \Ab_j)_+ \right|   \right\}. \quad (*) 
\end{align*}
From here, by using Lemma~\ref{lemma:exp_sup_v_bound_UniformSphere} in combination with Markov's inequality, we have with probability at least $1 - 1/t$ that 
\begin{align*}
(*) & \leq \max \left\{ \pb^T \bb \mid \Ebb (\pb^T \tilde{\Ab})_+ \leq \frac{1}{n} + \| \pb \|_2 \cdot \frac{t \sqrt{2}}{\sqrt{n}} \right\} \\
& = \max \left\{ \pb^T \bb \mid \frac{ \sqrt{2}}{ 2 \sqrt{\pi} \mu_m} \| \pb \|_2 \leq \frac{1}{n} + \| \pb \|_2 \cdot \frac{t \sqrt{2}}{\sqrt{n}} \right\} \\
& = \max \left\{ \pb^T \bb \mid \| \pb \|_2 \leq \frac{1}{n} \cdot \frac{1}{ \frac{ \sqrt{2}}{ 2 \sqrt{\pi} \mu_m} -  \frac{t \sqrt{2}}{\sqrt{n}} } \right\} \\
& = \frac{ \| \bb \|_2 }{n} \cdot \frac{1}{ \frac{ \sqrt{2}}{ 2 \sqrt{\pi} \mu_m} -  \frac{t \sqrt{2}}{\sqrt{n}} }, 
\end{align*}
where the first equality follows by applying Lemma~\ref{lemma:closed_form_pbTAb_plus_UniformSphere}. With regard to the feasibility of $P$, this again follows by the fact that the dual of $\min\{ \| \xb \|_{\infty} \mid \Ab \xb = \bb, \xb \geq \zerob \}$, which is always feasible, is bounded with probability at least $1 - 1/t$. \hfill \Halmos
\endproof

\subsection{Proof of Theorem~\ref{theorem:GMBernoulliCovering_gap_bound_whp}}

To establish Theorem~\ref{theorem:GMBernoulliCovering_gap_bound_whp}, we first require a simple adaptation of Lemma~\ref{lemma:beta_bound_on_P_distr}. The proof is straightforward, and omitted for brevity.
\begin{lemma}
Suppose that $P^{\covering}$ is feasible, $v(P) \geq 0$ and $\beta \geq \min \{ \| \xb \|_{\infty} \mid \Ab \xb = \bb, \xb \geq \zerob\}$. Suppose that $\xib$ is the uniform distribution over $[n]$, i.e., $\xi_j = 1/n$ for all $j \in [n]$. If $C = n \beta$, then $P^{\covering}_{\distr}$ is feasible and we have that 
\begin{equation*}
\Delta v(P^{\covering}_{\distr}) \leq \sqrt{n} \beta.
\end{equation*}
\label{lemma:beta_bound_on_P_covering_distr}
\end{lemma}

We next require the following lemma, which is a concentration result for the minimum of a collection of independent binomial random variables.
\begin{lemma}
Let $Y_1,\dots, Y_m$ be independent random variables, with each $Y_i \sim \Binomial(n, q_i)$. Let $\delta \in (0,1)$. Then with probability at least  $1 - \delta$, we have
\begin{equation*}
\min_{i \in [m]} Y_i \geq n \cdot \left( \min_{i \in [m]} q_i - \sqrt{ \frac{1}{2n} \log \frac{m}{\delta} } \right).
\end{equation*}
\label{lemma:minYi_bound_whp}
\end{lemma}
\proof{Proof:}
Let $\epsilon > 0$. Then for any $i \in [m]$, %
\begin{align*}
& \Pr( Y_i < \min_{i' \in [m]} \Ebb[Y_{i'}] - n \epsilon) \\
& \leq \Pr( Y_i  < \Ebb[Y_i] - n\epsilon)  \\
& \leq \exp( - \frac{2 n^2 \epsilon^2 }{ \sum_{j=1}^n (1 - 0)^2 } ) \\
& = \exp( - \frac{2 n^2 \epsilon^2 }{ n } ) \\
& = \exp( -2 n \epsilon^2 ).
\end{align*}
where the second inequality follows by Hoeffding's inequality. 

Now, observe that 
\begin{align*}
& \Pr( \min_{i \in [m]} Y_i < \min_{i' \in [m]} \Ebb[Y_{i'}] - n \epsilon) \\
& = \Pr\left( \bigcup_{i \in [m]} \{ Y_i < \min_{i' \in [m]} \Ebb[Y_{i'}] - n \epsilon \} \right) \\
& \leq \sum_{i = 1}^m \Pr( Y_i < \min_{i' \in [m]} \Ebb[Y_{i'}] - n \epsilon ) \\
& \leq m \cdot \exp( - 2 n \epsilon^2 ),
\end{align*}
where the first inequality follows by the union bound. This implies that 
\begin{align*}
\Pr\left( \min_{i \in [m]} Y_i \geq \min_{i' \in [m]} \Ebb[Y_{i'}] - n \epsilon \right) \geq 1 - m \cdot \exp( -2 n \epsilon^2 ).
\end{align*}
Note that $\epsilon$ was arbitrary; to make the right hand side of the previous bound equal to $1 - \delta$, we can solve for $\epsilon$ as
\begin{align*}
\delta & = m \cdot \exp(- 2 n \epsilon^2 ) \\
\Rightarrow \log \delta & = \log m - 2 n \epsilon^2 \\
\Rightarrow \epsilon^2 & = \frac{1}{2n} \log \frac{m}{\delta} \\
\Rightarrow \epsilon & = \sqrt{ \frac{1}{2n} \log \frac{m}{\delta} }.
\end{align*}
Thus, with probability at least $1 - \delta$, we have that 
\begin{align*}
\min_{i \in [m]} Y_i \geq \min_{i' \in [m]} \Ebb[Y_{i'}] - \sqrt{ \frac{n}{2} \log \frac{m}{\delta} } 
\end{align*}
and noting that $\Ebb[Y_{i'}] = n q_{i'}$, this is equivalent to
\begin{align*}
\min_{i \in [m]} Y_i \geq n \min_{i' \in [m]} q_{i'} - \sqrt{ \frac{n}{2} \log \frac{m}{\delta} },
\end{align*}
as required. \hfill \Halmos 
\endproof

With this concentration result in hand, we can now prove Theorem~\ref{theorem:GMBernoulli_infinity_norm_bound_whp}, which states that the minimum infinity norm of any feasible solution of $P$ is $O(1/n)$ with high probability. 

\begin{theorem}
Suppose that $P$ is generated according to generative model \GMBernoulliCoveringNum. Let $\delta \in (0,1)$. Suppose that $n > \log(m/\delta) / [ 2 (\min_{i \in [m]} q_i)^2 ]$. Then with probability at least $1 - \delta$, $P^{\covering}$ is feasible and
\begin{equation*}
\min \{ \| \xb \|_\infty \mid \Ab \xb \geq \bb, \xb \geq \zerob\} \leq \frac{ \max_{i \in [m]} b_i }{n} \cdot \frac{1}{ \min_{i' \in [m]} q_{i'} - \sqrt{ \frac{1}{2n} \log \frac{m}{\delta}} }
\end{equation*}
\label{theorem:GMBernoulli_infinity_norm_bound_whp}
\end{theorem}
\proof{Proof:}
We have
\begin{align} 
& \min \{ \| \xb \|_\infty \mid \Ab \xb \geq \bb, \xb \geq \zerob\} \nonumber \\
& = \max \{ \pb^T \bb \mid \sum_{j=1}^n ( \pb^T \Ab_j )_+ \leq 1, \pb \geq \zerob \} \nonumber \\
& = \max \{ \pb^T \bb \mid \sum_{j=1}^n \pb^T \Ab_j \leq 1, \pb \geq \zerob \} \nonumber \\
& = \max \{ \pb^T \bb \mid \sum_{i=1}^m p_i Y_i \leq 1, \pb \geq \zerob \}, \label{prob:cover_dual_Yi}
\end{align}
where $Y_i = \sum_{j=1}^n A_{i,j}$ for each $i \in [m]$. In the first step, we have simply taken the dual of the original problem; in the second step, we use the fact that $\pb \geq \zerob$ and $\Ab \geq \zerob$ to assert that $(\pb^T \Ab_j)_+ = \pb^T \Ab_j$; and in the third, we use the definition of the $Y_i$'s. 

By the definition of the generative model, we have that each $Y_i \sim \Binomial(n, q_i)$. Therefore, by Lemma~\ref{lemma:minYi_bound_whp}, we have that with probability at least $1 - \delta$,
\begin{equation}
\min_{i \in [m]} Y_i \geq n \min_{i' \in [m]} q_{i'} - \sqrt{ \frac{n}{2} \log \frac{m}{\delta} }. \label{eq:min_Yi_LB}
\end{equation}
In addition, by the assumption on $n$, it follows that the right hand side of \eqref{eq:min_Yi_LB} is positive, which implies that $Y_i > 0$ for all $i \in [m]$. Thus, when \eqref{eq:min_Yi_LB} holds, we can determine the optimal solution of problem~\eqref{prob:cover_dual_Yi} as follows: the optimal solution is given by $p_{i^*} = 1/ Y_{i^*}$ for $i^* = \arg \max_{i \in [m]} b_i / Y_i$ and $p_i = 0$ for all $i \neq i^*$. (We remind the reader here that the definition of generative model~\GMBernoulliCoveringNum requires $\bb$ to be nonnegative. We also note in the case that the arg max is not a singleton, we can set $i^*$ to be any maximizing index $i$.) 

When \eqref{eq:min_Yi_LB} holds, the objective value of \eqref{prob:cover_dual_Yi} can therefore be further refined as
\begin{align*}
& \max \{ \pb^T \bb \mid \sum_{i=1}^m p_i Y_i \leq 1, \pb \geq \zerob \}, \\
& = p_{i^*} \cdot b_{i^*} \\
& = \frac{ b_{i^*}}{Y_{i^*}} \\
& \leq \frac{ \max_{i \in [m]} b_i }{ \min_{i \in [m]} Y_i } \\
& \leq \frac{ \max_{i \in [m]} b_i }{ n \min_{i' \in [m]} q_{i'} - \sqrt{ \frac{n}{2} \log \frac{m}{\delta} } } \\ 
& = \frac{ \max_{i \in [m]} b_i }{n} \cdot \frac{1}{ \min_{i' \in [m]} q_{i'} - \sqrt{ \frac{1}{2n} \log \frac{m}{\delta}} },
\end{align*}
which holds with probability at least $1 - \delta$. This establishes the required bound on the objective value of the infinity norm problem in the theorem statement. 

To see that $P^{\covering}$ is feasible, observe that $P^{\covering} \equiv \min\{ \cb^T \xb \mid \Ab \xb \geq \bb, \xb \geq \zerob\}$ is feasible if and only if $\min \{ \| \xb \|_{\infty} \mid \Ab \xb \geq \bb, \xb \geq \zerob\}$ is feasible. The latter problem is feasible if and only if its dual problem $\max \{ \pb^T \bb \mid \sum_{j=1}^n ( \pb^T \Ab_j )_+ \leq 1, \pb \geq \zerob \}$, which is always feasible, is bounded. Our reasoning above establishes that this problem is bounded with probability at least $1 - \delta$, which implies that $P^{\covering}$ is feasible with probability at least $1 - \delta$. This completes the proof. \hfill \Halmos
\endproof

We can now complete the proof of Theorem~\ref{theorem:GMBernoulliCovering_gap_bound_whp}.
\proof{Proof of Theorem~\ref{theorem:GMBernoulliCovering_gap_bound_whp}:}
Let $\beta$ be defined as
\begin{equation*}
\beta = \frac{ \max_{i \in [m]} b_i }{n} \cdot \frac{1}{ \min_{i' \in [m]} q_{i'} - \sqrt{ \frac{1}{2n} \log \frac{m}{\delta}} }
\end{equation*}
By Theorem~\ref{theorem:GMBernoulli_infinity_norm_bound_whp}, we have that $P^{\covering}$ is feasible and that $\min \{ \| \xb \|_{\infty} \mid \Ab \xb \geq \bb, \xb \geq \zerob \} \leq \beta$ with probability at least $1 - \delta$. By Lemma~\ref{lemma:beta_bound_on_P_covering_distr}, observe that by setting $C$ as 
\begin{align*}
C & = n \beta \\
& = \max_{i \in [m]} b_i  \cdot \frac{1}{ \min_{i' \in [m]} q_{i'} - \sqrt{ \frac{1}{2n} \log \frac{m}{\delta}} },
\end{align*}
we have that $P^{\covering}_{\distr}$ is feasible and 
\begin{align*}
\Delta v(P^{\covering}_{\distr}) & \leq \sqrt{n} \beta \\
& = \frac{ \max_{i \in [m]} b_i }{\sqrt{n}} \cdot \frac{1}{ \min_{i' \in [m]} q_{i'} - \sqrt{ \frac{1}{2n} \log \frac{m}{\delta}} },
\end{align*}
as desired. \hfill \Halmos
\endproof

\section{Special Structures and Extensions}
\label{sec:special_structures}

In this section, we demonstrate how the results of Sections~\ref{sec:main_results} and \ref{sec:proofs} can be applied to LPs with specific problem structures, including LPs with totally unimodular constraints (Section~\ref{subsec:special_structures_TU}), Markov decision processes (Section~\ref{subsec:special_structures_MDP}), covering problems (Section~\ref{subsec:special_structures_covering}) and packing problems (Section~\ref{subsec:special_structures_packing}). In Section~\ref{subsec:special_structures_portfolio}, we consider the portfolio optimization problem, which is in general not an LP, but is amenable to the same type of analysis.

\subsection{LPs with Totally Unimodular Constraints}
\label{subsec:special_structures_TU}

Consider a linear program with a totally unimodular constraint matrix, i.e., every square submatrix of $\Ab$ has determinant $0$, $1$, or $-1$. Such LPs appear in various applications, such as minimum cost network flow problems and assignment problems \citep{bertsekas1998network}. 
In such problems, it is not uncommon to encounter the situation where the number of variables is much larger than the number of constraints. For example, in a minimum cost network flow problem, each constraint corresponds to a flow-balance constraint at a given node, while each variable corresponds to the flow over an edge; in a graph of $n$ nodes, one will therefore have $n$ constraints and as many as ${n}\choose{2}$ decision variables. We can thus consider solving the problem using the column randomization method. We obtain the following guarantee on the objective value of the column randomization method when applied to linear programs with totally unimodular constraints.
\begin{proposition}
	When $\Ab$ is totally unimodular, then 
	\begin{align*}
	\gamma = m \| \cb \|_{\infty}
	\end{align*}
	is a valid upper bound on $\| \pb \|_{\infty}$ for every basic solution $\pb$ of $D$. 
\end{proposition}
\proof{Proof:}
Any basic solution $\pb$ to the dual problem $D$ can be written as $\pb^T =  \cb_B^T \Ab_B^{-1} $, where $B$ is a basis. In addition, since $\Ab$ is totally unimodular, any element of $\Ab_B^{-1}$ is either $1$, $-1$, or $0$. Therefore, the $i$th component of $\pb$ satisfies $p_i = \sum_{j=1}^m [\Ab_B^{-1}]_{ji} (\cb_B)_j \leq m \cdot \| \cb \|_\infty$ for all $i \in [m]$. Thus, $\gamma = m \| \cb \|_\infty$ is a valid upper bound on $\| \pb \|_{\infty}$. \hfill \Halmos
\endproof

Using this result together with the observation that $\| \Ab \|_{\max} = 1$ for any totally unimodular matrix $\Ab$, we can invoke Theorem~\ref{thm:main_largest_abs_of_all_dual_BFS} to obtain the following performance guarantee for column randomization when applied to LPs with totally unimodular matrices.
\begin{corollary}
	\label{prop:totally_unimodular}
	Assume the constraint matrix of $\Ab$ of the complete problem $P$ is totally unimodular. Define $C$, $P_\distr$, $P_J$ and $\Ab_J$ as in Theorem~\ref{thm:main_largest_abs_of_all_dual_BFS}. For any $\delta \in (0,1)$, with probability at least $1 - \delta$ over the set $J$, the following holds: if $P_J$ is feasible and $\rank(\Ab_J) = m$, then
	\begin{align}
	\label{eq:thm_convergence_totally_unimodular}
	\Delta v(P_J) \leq \Delta v(P_\distr) + \frac{C \left( 1 + m^2  \| \cb \|_\infty  \right)}{\sqrt{K}} \left(  1 + \sqrt{2 \log \frac{2}{\delta}}  \right).
	\end{align}
\end{corollary}

\subsection{Markov Decision Processes}
\label{subsec:special_structures_MDP}

Consider a discounted infinite horizon MDP, with $n_s$ states and $n_a$ actions. The cost function $c(s,a)$ represents the immediate cost of taking action $a$ in state $s$. The transition probability $P_{s}(s',a)$ represents the probability of being in state $s'$ after taking action $a$ in state $s$. Let $\theta \in (0,1)$ be the discount factor. One can solve the MDP by formulating a linear program \citep{manne1960linear}:
\begin{alignat*}{3}
\label{problem:LP_MDP_0}
\underset{\xb_1,\ldots,\xb_{n_s} \in \mathbb{R}^{n_a}}{\text{minimize}}  \quad & \cb_1^T\xb_1 + \ldots + \cb_s^T \xb_s + \ldots +  \cb_{n_s}^T \xb_{n_s}  \\   \text{such that}  \quad & (\Eb_1 - \theta \Pb_1) \xb_1 + \ldots + (\Eb_s - \theta \Pb_s) \xb_j + \ldots + (\Eb_{n_s} - \theta \Pb_{n_s}) \xb_{n_s} = \oneb,
\\   & \xb_1,\ldots,\xb_s,\ldots,\xb_{n_s} \geq \zerob,
\end{alignat*}
where %
$\Eb_j$ is a $n_s \times n_a$ matrix such that the $j$th row is all ones and every other entry is zero. The vector $\cb_s$ is of size $n_a$ such that its $a$th component is equal to $c(s,a)$. The matrix $\Pb_s$ is of size $n_s \times n_a$ such that its $(s',a)$-th component represents $P_s(s',a)$. Notice that matrix $\Pb_s$ is a column stochastic matrix, i.e., $\oneb^T \Pb_s = \oneb^T$ and $\Pb_s \geq \zerob$ for all $s \in [n_s]$. The decision variable vector $\xb_s$ is of size $n_a$, where the $a$th entry represents the expected discounted long-run frequency of the system being in state $s$ and action $a$ being taken. If one sorts the decision variables by actions \citep{ye2005new}, then the linear program can be re-written as:
\begin{subequations}
	\label{problem:LP_MDP}
	\begin{alignat}{3}
	\underset{\tilde{\xb}^1,\ldots,\tilde{\xb}^{n_a} \in \mathbb{R}^{n_s}}{\text{minimize}}  \quad & \tilde{\cb}_1^T\tilde{\xb}_1 + \ldots + \tilde{\cb}_a^T \tilde{\xb}_a + \ldots +  \tilde{\cb}_{n_a}^T \tilde{\xb}_{n_a}  \\   \text{such that}  \quad & (\Ib - \theta \tilde{\Pb}_1) \tilde{\xb}_1 + \ldots + (\Ib - \theta \tilde{\Pb}_a) \tilde{\xb}_a + \ldots + (\Ib - \theta \tilde{\Pb}_{n_a}) \tilde{\xb}_{n_a} = \oneb,
	\\   & \tilde{\xb}_1,\ldots,\tilde{\xb}_a,\ldots,\tilde{\xb}_{n_a} \geq \zerob,
	\end{alignat}
\end{subequations}
where $\tilde{\cb}_a = [c(1,a); \ldots; c(s,a); \ldots; c(n_s,a)]$ for $a \in [n_a]$ and $\tilde{\Pb}_a$ is a $n_s \times n_s$ matrix such that its $(s',s)$-th element is equal to $P_s(s',a)$. 
Note that problem~\eqref{problem:LP_MDP} is a standard form LP and has more columns than rows. 
We can therefore apply the column randomization method to solve problem~\eqref{problem:LP_MDP}. To adapt our performance guarantee from Section~\ref{subsec:sampling_columns_code_and_theorems}, we establish a bound $\gamma$ on $\| \pb \|_{\infty}$ for every dual basic solution $\pb$ that is specific to problem~\eqref{problem:LP_MDP}. 
\begin{proposition}
	For the infinite horizon discounted MDP problem~\eqref{problem:LP_MDP}, then 
	\begin{align*}
	\gamma = \frac{\| \cb \|_{\infty} }{1 - \theta}
	\end{align*}
	is a valid upper bound on $\| \pb \|_{\infty}$ for any basic solution $\pb$ of the dual of problem~\eqref{problem:LP_MDP}. 
\end{proposition}
\proof{Proof:}
Any basic solution $\pb$ of the dual has the form $\pb^T = \cb_B^T \Ab^{-1}_B$, where $B$ is a basis of the linear program~\eqref{problem:LP_MDP}. Note that $\Ab_B$ has the form $\Ab_B = \Ib - \theta \Pb$, where $\Pb$ is an $n_s \times n_s$ matrix such that each of its columns is selected from the columns of $[\tilde{\Pb}_1,\ldots,\tilde{\Pb}_{n_a}]$ \citep[see][]{ye2005new}. In addition, a standard property of $\Ab^{-1}_B$ is that it can be written as the following infinite series:
\begin{align*}
\Ab^{-1}_B = \Ib + \theta \Pb + \theta^2 \Pb^2 + \dots = \Ib + \sum_{n=1}^{\infty} \theta^n \cdot \Pb^n.
\end{align*}
Thus, we can bound $\| \pb \|_{\infty}$ as $\| \pb^T \|_\infty \leq \| \cb^T_B \|_\infty + \sum_{n=1}^\infty \theta^n \cdot \| \cb^T_B \Pb^n  \|_\infty$. Note that for any $n \in \Nbb$ and vector $\vb \in \Rbb^{n_s}$, we have
\begin{align*}
\| \vb^T \Pb^n \|_{\infty}  & = \max_{s \in [n_s] } \left| \sum_{s' \in [n_s]} v_{s'} \Pb^n_{(s', s)} \right|  \\
& \leq \max_{s \in [n_s]} \sum_{s' \in [n_s]} |v_{s'}| \cdot \Pb^n_{(s',s)} \\
& \leq  \| \vb \|_{\infty} \cdot \max_{s \in [n_s]}\sum_{s' \in [n_s]} \Pb^n_{(s',s)} \\
& = \| \vb \|_{\infty},
\end{align*}
where $\Pb^n_{(s',s)}$ is the $(s',s)$th entry of $\Pb^n$. Therefore, we obtain that 
\begin{align*}
\| \pb^T \|_\infty & = \| \cb^T_B \Ab^{-1}_B \|_{\infty} \\
& \leq \| \cb_B \|_\infty + \sum_{n=1}^\infty \theta^n \cdot \| \cb^T_B \Pb^ n \|_{\infty} \\
& \leq \| \cb_B \|_\infty / (1 - \theta) \\
& \leq \| \cb \|_{\infty} / (1 - \theta).
\end{align*}
Since $\pb$ was an arbitrary basic solution of the complete dual of problem~\eqref{problem:LP_MDP}, we can therefore set $\gamma = \| \cb \|_\infty / (1 - \theta)$. \hfill \Halmos
\endproof

With this result in hand, and observing that $\| \Ab \|_{\max} \leq 1$, we can apply Theorem~\ref{thm:main_largest_abs_of_all_dual_BFS} to obtain the following performance guarantee for column randomization in the case of discounted infinite horizon MDPs. 

\begin{corollary}
	\label{thm:MDP}
	Consider solving a discounted infinite horizon MDP with $n_s$ states and $n_a$ actions by the column randomization method. Define $C$, $P_\distr$, $P_J$ and $\Ab_J$ as in Theorem~\ref{thm:main_largest_abs_of_all_dual_BFS}. For any $\delta \in (0,1)$, with probability at least $1 - \delta$, the following holds: if $P_J$ is feasible and $\rank(\Ab_J) = n_s$, then
	\begin{align}
	\label{eq:thm_convergence_totally_MDP}
	\Delta v(P_J) \leq  \Delta v(P_\distr) + \frac{C}{\sqrt{K}} \cdot \left( 1 + \frac{n_s \| \cb \|_\infty}{1 - \theta} \right) \cdot \left(  1 + \sqrt{2 \log \frac{2}{\delta}}  \right).
	\end{align}
\end{corollary}

\subsection{Covering Problems}
\label{subsec:special_structures_covering}

A covering linear program can be formulated as 
\begin{subequations}
	\begin{alignat}{2}
	P^{\covering} \ : \ & \underset{\xb}{\text{minimize}} & & \cb^T \xb \\
	& \text{subject to} & \quad & \Ab \xb \geq \bb, \\
	& & & \xb \geq \zerob,
	\end{alignat}
\end{subequations}
where $\Ab$, $\bb$ and $\cb$ are all nonnegative, and we additionally assume that for every $i \in [m]$, there exists a $j \in [n]$ such that $A_{i,j} > 0$. This type of problem arises in numerous applications such as facility location \citep{owen1998strategic}. The column-randomized counterpart of this problem and its dual can be written as
\begin{align*}
P^{\covering}_J \ & : \ \min\{ \cb_J^T \tilde{\xb} \mid \Ab_J \tilde{\xb} \geq \bb, \tilde{\xb} \geq \zerob \}, \\
D^{\covering}_J \ & : \ \max\{ \pb^T \bb \mid \pb^T \Ab_J \leq \cb^T_J, \pb \geq \zerob \}.
\end{align*}
Although $P^{\covering}$ is not a standard form LP, it is straightforward to extend Proposition~\ref{prop:objective_bound_by_dual_solution} to this problem, leading to the following result. We omit the proof for brevity.
\begin{proposition}
	\label{prop:covering_objective_bound_by_dual_solution}
	Let $C$ be a nonnegative constant and define $P^{\covering}_{\distr}$ as 
	\begin{equation*}
	P^{\covering}_{\distr} \equiv \min\{ \cb^T \xb \mid \Ab \xb \geq \bb, \zerob \leq \xb \leq C \xib \}.
	\end{equation*}
	For any $\delta \in (0,1)$, with probability at least $1 - \delta$ over the sample $J$, the following holds: if $P^{\covering}_J$ is feasible, then
	\begin{equation*}
	\Delta v(P^{\covering}_J) \leq \Delta v(P^{\covering}_{\distr} ) + \frac{C}{\sqrt{K}} \cdot (1 + \| \pb \|_{\infty} \cdot m \cdot \|\Ab\|_{\max}) \cdot \left(1 + \sqrt{ 2 \log \frac{2}{\delta} } \right)
	\end{equation*}
	for any optimal solution $\pb$ of $D^{\covering}_J$. 
\end{proposition}
To now use this result, we need to be able to bound $\| \pb \|_{\infty}$ for any solution $\pb$ of any dual $D^{\covering}_J$ of the column-randomized problem. Let us define the quantity $U^{\covering}$ as 
\begin{equation*}
U^{\covering} = \max_{i,j} \left\{ \frac{c_j}{A_{i,j}} \ \vline \ A_{i,j} > 0 \right\}.
\end{equation*}
We then have the following result. 
\begin{proposition}
	Let $J \subseteq [n]$, and suppose that $P^{\covering}_J$ is feasible. Then for any feasible solution $\pb$ of $D^{\covering}_J$, $\| \pb \|_{\infty} \leq U^{\covering}$. 
\end{proposition}
\proof{Proof:}
Fix an $i \in [m]$, and consider the LP
\begin{equation}
D^{\text{B}-\covering}_J \ : \ \max \{ p_i \mid \pb^T \Ab_J \leq \cb_J^T, \ \pb \geq \zerob \}. 
\end{equation}
The optimal objective value of this problem, $v(D^{\text{B}-\covering}_J)$, is an upper bound on $p_i$ for any  feasible solution $\pb$ of $D^{\covering}_J$ (and thus, it is an upper bound on $p_i$ for any optimal solution $\pb$ of $D^{\covering}_J$). 
Consider the dual of this problem:
\begin{equation}
P^{\text{B}-\covering}_J \ : \ \min \{ \cb_J^T \tilde{\xb} \mid \Ab_J \tilde{\xb} \geq \eb_i, \ \tilde{\xb} \geq \zerob\},
\end{equation}
where $\eb_i$ is the $i$th standard basis vector for $\Rbb^m$. By weak duality, the objective value of any feasible solution of $P^{\text{B}-\covering}_J$ is an upper bound on $v(D^{\text{B}-\covering}_J)$. 

We now construct a particular feasible solution. Let $j'$ be any column in $J$ such that $A_{i,j'} > 0$; such a column is guaranteed to exist by our assumption on the matrix $\Ab$. Define a solution $\tilde{\xb}$ as 
\begin{equation*}
\tilde{x}_{j} = \left\{ \begin{array}{lll} 1 / A_{i,j'} & \text{if}\ j = j', \\
0 & \text{otherwise}. \end{array} \right.
\end{equation*}
It is easy to see that $\tilde{\xb}$ is a feasible solution of $P^{\text{B}-\covering}_J$, and that its objective value is $\cb^T_J \tilde{\xb} = c_{j'} / A_{i,j'}$. Since this objective value is bounded by $U^{\covering}$, it follows that $U^{\covering} \geq \max\{ p_i \mid \pb^T \Ab_J \leq \cb^T_J, \ \pb \geq \zerob\}$. 

Since our choice of $i$ was arbitrary, it follows that $\| \pb \|_{\infty} \leq U^{\covering}$ for any feasible solution of $D^{\covering}_J$. \hfill \Halmos
\endproof

Using this result together with Proposition~\ref{prop:covering_objective_bound_by_dual_solution} yields the following guarantee. 
\begin{corollary}
	\label{thm:covering_objective_bound_by_c_over_Aij}
	Let $C$ and $P^{\covering}_{\distr}$ be defined as in Proposition~\ref{prop:covering_objective_bound_by_dual_solution}. For any $\delta \in (0,1)$, with probability at least $1 - \delta$ over the sample $J$, the following holds: if $P^{\covering}_J$ is feasible, then
	\begin{equation*}
	\Delta v(P^{\covering}_J) \leq \Delta v(P^{\covering}_{\distr} ) + \frac{C}{\sqrt{K}} \cdot (1 + U^{\covering} \cdot m \cdot \|\Ab\|_{\max}) \cdot \left(1 + \sqrt{ 2 \log \frac{2}{\delta} } \right).
	\end{equation*}
\end{corollary}

\subsection{Packing Problems}
\label{subsec:special_structures_packing}

A packing linear program is defined as 
\begin{subequations}
	\begin{alignat}{2}
	P^{\packing} \ : \ & \underset{\xb}{\text{maximize}} & & \cb^T \xb \\
	& \text{subject to} & \quad & \Ab \xb \leq \bb, \\
	& & & \xb \geq \zerob,
	\end{alignat}
\end{subequations}
where we assume that $\cb \geq \zerob$, $\bb > \zerob$, and that $\Ab$ is such that for every column $j \in [n]$, there exists an $i \in [m]$ such that $A_{i,j} > 0$. Packing problems have numerous applications, such as network revenue management \citep{talluri2006theory}.

The column-randomized counterpart of this problem and its dual can be written as
\begin{align*}
P^{\packing}_J \ & : \ \max\{ \cb_J^T \tilde{\xb} \mid \Ab_J \tilde{\xb} \leq \bb, \tilde{\xb} \geq \zerob \}, \\
D^{\packing}_J \ & : \ \min\{ \pb^T \bb \mid \pb^T \Ab_J \geq \cb^T_J, \pb \geq \zerob \}.
\end{align*}
As with covering problems, the packing problem $P^{\packing}$ is not a standard form LP, but we can derive a counterpart of Proposition~\ref{prop:objective_bound_by_dual_solution} for $P^{\packing}$. Note that in this guarantee, for a problem $P'$ with the same feasible region as $P^{\packing}$, the optimality gap $\Delta v(P')$ is defined as $\Delta v(P') = v(P^{\packing}) - v(P')$, since the complete problem $P^{\packing}$ is a maximization problem. As with Proposition~\ref{prop:covering_objective_bound_by_dual_solution}, the proof is straightforward, and thus omitted. 
\begin{proposition}
	\label{prop:packing_objective_bound_by_dual_solution}
	Let $C$ be a nonnegative constant and define $P^{\covering}_{\distr}$ as 
	\begin{equation*}
	P^{\packing}_{\distr} \equiv \max\{ \cb^T \xb \mid \Ab \xb \leq \bb, \zerob \leq \xb \leq C \xib \}.
	\end{equation*}
	For any $\delta \in (0,1)$, with probability at least $1 - \delta$ over the sample $J$, the following holds: if $P^{\packing}_J$ is feasible, then
	\begin{equation*}
	\Delta v(P^{\packing}_J) \leq \Delta v(P^{\packing}_{\distr} ) + \frac{C}{\sqrt{K}} \cdot (1 + \| \pb \|_{\infty} \cdot m \cdot \|\Ab\|_{\max}) \cdot \left(1 + \sqrt{ 2 \log \frac{2}{\delta} } \right)
	\end{equation*}
	for any optimal solution $\pb$ of $D^{\packing}_J$. 
\end{proposition}

To obtain a more specific guarantee, define for each $i$ the following quantities:
\begin{align*}
r_i & = \max \left\{ \frac{c_j}{A_{i,j}} \ \vline \ A_{i,j} > 0 \right\}, \\
j^*_i & = \arg \max_j \left\{ \frac{c_j}{A_{i,j}} \ \vline \ A_{i,j} > 0 \right\}.
\end{align*}
These two quantities can be understood by interpreting each $i$ as a resource constraint, and $b_i$ as the available amount of resource $i$. The column $j^*_i$ is the column that has the best rate of objective value garnered per unit of resource $i$ consumed, and the quantity $r_i$ is that corresponding rate. Define now $W$ as 
\begin{equation*}
W = \sum_{i' = 1}^m r_{i'} b_{i'},
\end{equation*}
and $U^{\packing}$ as the maximum over $i$ of $W / b_i$, i.e.,
\begin{equation*}
U^{\packing} = \max_{i \in [m]} \frac{ W }{b_i}  = \frac{ \sum_{i'=1}^m r_{i'} b_{i'}}{ \min_{i \in [m]} b_i }.
\end{equation*}
With these definitions, we can establish that $U^{\packing}$ is an upper bound on the infinity norm of any dual optimal solution $\pb$ of $P^{\packing}_J$. 
\begin{proposition}
	Let $J \subseteq [n]$. Then any optimal solution $\pb$ of $D^{\packing}_J$ satisfies $\| \pb \|_{\infty} \leq U^{\packing}$. 
\end{proposition}
\proof{Proof:}
We first establish a useful property of $W$: the quantity $W$ is actually an upper bound on $v(P)$. To see this, define the solution $\tilde{\xb}^{(i)}$ for each $i$ as 
\begin{equation*}
\tilde{\xb}^{(i)} = \frac{b_i}{A_{i, j^*_i} } \cdot \eb_{j^*_i},
\end{equation*}
and define $\tilde{\xb} = \sum_{i=1}^m \tilde{\xb}^{(i)}$. Let $\xb$ be any feasible solution of the complete problem $P^{\packing}$. Note that for each $\tilde{\xb}^{(i)}$, we have:
\begin{align*}
\cb^T \tilde{\xb}^{(i)} & = \frac{ c_{j^*_i} b_i}{ A_{i, j^*_i} } \\
& \geq \frac{c_{j^*_i} }{ A_{i, j^*_i}} \left[ \sum_{j =1}^n A_{i,j} x_j \right] \\
& = \frac{c_{j^*_i} }{ A_{i, j^*_i}} \left[ \sum_{j: A_{i,j} > 0} A_{i,j} x_j \right] \\
& \geq \sum_{j: A_{i,j} > 0} A_{i,j} \cdot \frac{c_{j}}{A_{i,j}} \cdot x_j \\
& = \sum_{j: A_{i,j} > 0} c_j x_j.
\end{align*}
where the first inequality follows because $\xb$ satisfies $\Ab \xb \leq \bb$, and the second follows by the definition of $j^*_i$. Using this bound, we have
\begin{align*}
\cb^T \tilde{\xb} & = \sum_{i=1}^m \cb^T \tilde{\xb}^{(i)} \\
& \geq \sum_{i=1}^m \left[ \sum_{j : A_{i,j} > 0} c_j x_j \right] \\
& \geq \sum_{j=1}^n c_j x_j \\
& = \cb^T \xb,
\end{align*}
where the second inequality follows by our assumption that for each $j$, there exists an $i$ such that $A_{i,j} > 0$.

Now, let us fix an $i \in [m]$. We wish to bound $|p_i|$ for an optimal solution $\pb$ of $D^{\packing}_J$. We can compute a bound on $| p_i|$ by solving the following LP:
\begin{equation*}
D^{\text{B}-\packing}_J \ : \ \max \{ p_i \mid  \pb^T \bb \leq v(P^{\packing}_J),\  \pb^T \Ab_J \geq \cb_J^T, \ \pb \geq \zerob \}. 
\end{equation*}
Note that by weak duality, the feasible region of $D^{\text{B}-\covering}_J$ is exactly the set of all optimal solutions to the sampled dual problem, $D^{\packing}_J$. Observe that for any $J$, $v(P^{\packing}_J) \leq v(P^{\packing}) \leq W$. Thus, a valid upper bound on $v(D^{\text{B}-\packing}_J)$ can be obtained by solving the following relaxation of $D^{\text{B}-\packing}_J$:
\begin{equation*}
D^{\text{B}-\packing-\text{rlx}}_J \ : \ \max \{ p_i \mid  \pb^T \bb \leq W, \ \pb \geq \zerob \}. 
\end{equation*}
This problem is a valid relaxation, because we have simply removed the constraint $\pb^T \Ab_J \geq \cb_J^T$, and we have replaced the value $v(P^{\packing}_J)$ with the larger value of $W$. The optimal objective value of this relaxation is simply $W / b_i$. Therefore, we obtain that for any dual optimal solution $\pb$ of $D^{\packing}_J$, $| p_i | \leq W / b_i$. It follows that $\| \pb \|_{\infty} \leq \max_{i \in [m]} (W / b_i) \equiv U^{\packing}$, for any optimal solution $\pb$ of $D^{\packing}_J$. \hfill \Halmos
\endproof

By combining this result with Proposition~\ref{prop:packing_objective_bound_by_dual_solution}, we obtain the following specific guarantee for packing LPs. 
\begin{corollary}
	\label{thm:packing_objective_bound_by_U}
	Let $C$ and $P^{\packing}_{\distr}$ be defined as in Proposition~\ref{prop:packing_objective_bound_by_dual_solution}. For any $\delta \in (0,1)$, with probability at least $1 - \delta$ over the sample $J$, the following holds: if $P^{\packing}_J$ is feasible, then
	\begin{equation*}
	\Delta v(P^{\packing}_J) \leq \Delta v(P^{\packing}_{\distr} ) + \frac{C}{\sqrt{K}} \cdot (1 + U^{\packing} \cdot m \cdot \|\Ab\|_{\max}) \cdot \left(1 + \sqrt{ 2 \log \frac{2}{\delta} } \right).
	\end{equation*}
\end{corollary}
With regard to $U^{\packing}$ which appears in this guarantee, we note that this constant depends on the constant $W$. Our choice of $W$ is special only in that it bounds $v(P^{\packing}_J)$. For particular packing problems, if one has access to a problem-specific bound $W'$ on $v(P^{\packing}_J)$, one could define $U^{\packing}$ with $W'$ instead to obtain a more refined bound.

\subsection{Portfolio Optimization}
\label{subsec:special_structures_portfolio}
In this last section, we deviate slightly from our previous examples by showing how our approach can be applied to problems that are not linear programs. The specific problem that we consider is the portfolio optimization problem, which is defined as
\begin{subequations}
	\label{problem:portfolio}
	\begin{alignat}{2}
	P^{\text{portfolio}}: \quad \underset{\xb \in \Rbb^n, \rb \in \Rbb^m}{\text{minimize}}  \quad & f(r_1,\ldots,r_m) \\
	\text{such that}  \quad & \sum_{j=1}^n \alpha_{ij} x_j = r_i, \quad\forall i \in [m] \\
	& \sum_{j=1}^n x_j = 1,\\
	& \xb \geq \zerob,
	\end{alignat}
\end{subequations}
where both $\xb$ and $\rb$ are decision variables. Problem~\eqref{problem:portfolio} can be interpreted as follows: a decision maker seeks an optimal portfolio, which is a distribution over instruments, according to some objectives. The decision variable $x_j$ represents the fraction of allocation committed to instrument $j$, the constraint parameter $\alpha_{ij}$ represents the return of instrument $j$ in scenario $i$, and $r_i$ is the total return in $i$th scenario. The objective function $f$ is a function measuring the risk of the returns $r_1,\dots,r_m$. Unlike the optimization problems we discussed so far, we assume that $f$ is any Lipschitz continuous function with Lipschitz constant $L$, and is not necessarily a linear function of $\rb$.

Although problem $P^{\text{portfolio}}$ is not in general a linear program, we can still apply the column randomization method to solve the problem. We describe the procedure in Algorithm~\ref{alg:portfolio}. Notice that, unlike Algorithm~\ref{alg:main} which samples columns associated with all variables, here we only sample columns associated with  $\xb$.

\begin{algorithm}
	\SingleSpacedXI
	\caption{The Column Randomization Method - Portfolio Optimization}
	\label{alg:portfolio}
	\begin{algorithmic}[1]
		\STATE Sample $K$ i.i.d. indices in $[n]$ as $J \equiv \{ J_1,\ldots,J_K \}$ according to a randomization scheme $\rho$.
		\STATE Solve the sampled optimization problem:
		\begin{equation}
		\label{problem:sampled_portfolio}
		P^{\text{portfolio}}_J  : \quad \min \left\{ f(\rb) \ \vline \  \sum_{j \in J} \alpha_{ij} \tilde{x}_j = r_i, \ \forall \ i \in [m], \ \sum_{j \in J} \tilde{x}_j = 1, \ \tilde{\xb} \geq \zerob \right\} 
		\end{equation}
		\RETURN optimal solution $\left( \tilde{\xb}^*,\rb^* \right)$ and optimal objective value $f(\rb^*)$
	\end{algorithmic}
\end{algorithm}

For $P^{\port}_J$ that is produced and solved by Algorithm~\ref{alg:portfolio}, we have the following performance guarantee. 
\begin{proposition}
	\label{prop:portfolio}
	Assume vectors $\alphab_j = (\alpha_{ij})_{i \in [m]}$ in problem $P^{\text{portfolio}}$ satisfying $\| \alphab_j \|_2 \leq H$ for all $j \in [n]$. Let $C \geq 1$ be an arbitrary constant and define the optimization problem
	\begin{align}
	\label{problem:portfolio_restricted}
	P^{\port}_{\distr}: \quad \underset{\xb,\rb}{\min} \left\lbrace f(\rb) \ \vline \ \sum_{j \in [n]} \alphab_j x_j = \rb, \,\, \oneb^T \xb = 1 , \,\, \zerob \leq \xb \leq C \xib \right\rbrace.
	\end{align}
	Denote $F$, $F_\distr$, and $F_J$ as optimal objective values of problems $P^\port$, $P^\port_{\distr}$, and $P^\port_J$, respectively. Define $\Delta F_J \equiv F_J - F$ and $\Delta F_\distr = F_\distr - F$. For any $\delta \in (0,1)$, with probability at least $1 - \delta$, the following statement holds:
	\begin{align}
	\label{eq:thm_convergence_portfolio}
	\Delta F_J \leq \Delta F_\distr + \frac{CLH}{\sqrt{K}} \left(  1 + 3 \sqrt{\frac{1}{2}\log \frac{4}{\delta}}  \right).
	\end{align}
\end{proposition}
While the proof (see below) is similar to that of Proposition~\ref{prop:objective_bound_by_dual_solution} in the construction of a random solution that is close to the solution of the distributional counterpart problem~$P^{\port}_{\distr}$, the main difference is that it relies on Lipschitz continuity, rather than LP duality. 

It is worthwhile to point out several aspects about this result and the portfolio optimization problem. First, the portfolio optimization problem~\eqref{problem:portfolio} is not required to be a convex optimization problem; the objective function $f$ can be non-convex, so long as it is Lipschitz continuous. Second, this result is related to a more specific result from our prior work \citep{chen2019decision}. In that paper, we consider the problem of estimating the decision forest choice model, which is a probability distribution over a collection of decision trees, and show that by solving an optimization problem over a random sample of trees, one can obtain a gap on the $\ell_1$ training error of the model that decays with rate $1 / \sqrt{K}$ (Theorem 5 of \citealt{chen2019decision}). Proposition~\ref{prop:portfolio} is a generalization of that result to more general optimization problems outside of choice model estimation, and allows for objective functions more general than those based on $\ell_1$ distance.

\proof{Proof of Proposition~\ref{prop:portfolio}:}

Let $(\xb^{*0}, \rb^{*0})$ be an optimal solution of $P^\port_\distr$. Consider the solution $(\xb', \rb')$ defined relative to the sample $J$: 
\begin{align}
\xb' & = \frac{1}{K} \sum_{k=1}^K \frac{x^{*0}_{j_k}}{ \xi_{j_k} } \eb_{j_k},\\
\rb' & = \sum_{j \in [n]} \alphab_{j} x'_j = \frac{1}{K} \sum_{k=1}^K ( x^{*0}_{j_k} / \xi_{j_k} ) \alphab_{j_k}.
\end{align}
The significance of $(\xb', \rb')$ is that we will be able to show that $\rb'$ will be close to $\rb^{*0}$, and that $f(\rb')$ will be close to $f(\rb^{*0}) = F_{\distr}$. However, $(\xb', \rb')$ is not necessarily a feasible solution to problem $P^\port$, because $\xb'$ will in general not satisfy the unit sum constraint. To turn it into a feasible solution for problem $P^\port$, we consider the solution $(\xb'', \rb'')$ obtained by normalizing $\xb'$ by its sum:
\begin{align}
\xb'' = \frac{\xb'}{ \oneb^T \xb'}, \\
\rb'' = \frac{ \rb'}{ \oneb^T \xb'}.
\end{align}
Note that $(\xb'', \rb'')$ is a feasible solution of $P^\port_J$. 

To understand why we consider $(\xb', \rb')$ and $(\xb'', \rb'')$, we show how these two solutions can be used to bound the difference between $F_J$ and $F_\distr$. Let $(\xb, \rb)$ be an optimal solution of $P^\port_J$. We now bound $F_J - F_\distr$ as follows:
\begin{align}
F_J - F_\distr & = f(\rb) - f(\rb^{*0}) \nonumber \\
& \leq f(\rb'') - f(\rb^{*0}) \nonumber \\
& = f(\rb'') - f(\rb') + f(\rb') - f(\rb^{*0}) \nonumber \\
& \leq | f(\rb'') - f(\rb')| + |f(\rb') - f(\rb^{*0})| \nonumber \\
& \leq L\| \rb'' - \rb' \|_2 + L \| \rb' - \rb^{*0} \|_2 \label{bound:portfolio_Lipschitz}
\end{align}
where the first step follows by the definitions of $(\xb, \rb)$ and $(\xb^{*0}, \rb^{*0})$; the second step follows because $(\xb'', \rb'')$ is a feasible solution of $P^{\port}_J$; the third and fourth step follow by algebra and basic properties of absolute values; and the last step follows by the fact that $f$ is Lipschitz continuous with constant $L$. 

We now proceed to show that $\| \rb' - \rb^{*0} \|_2$ and $\| \rb'' - \rb'\|_2$ can be bounded with high probability. \\

\noindent \textbf{Bounding $\| \rb' - \rb^{*0} \|_2$}: To bound this term, let us define for each $k \in [K]$ the random vector $\rb_{j_k}$ as 
\begin{equation*}
\rb_{j_k} = \frac{x^{*0}_{j_k} }{\xi_{j_k}} \alphab_{j_k}.
\end{equation*}
We make three important observations about $\rb_{j_1}, \dots, \rb_{j_K}$. First, for each $k$, the norm of $\rb_{j_k}$ is bounded as
\begin{align*}
\| \rb_{j_k} \|_2 = \left\| \frac{x^{*0}_{j_k} }{\xi_{j_k}} \alphab_{j_k} \right\|_2 \leq \frac{x^{*0}_{j_k} }{\xi_{j_k}} \cdot \left\| \alphab_{j_k} \right\|_2 \leq \frac{C \xi_{j_k}}{\xi_{j_k}}\cdot H = CH.
\end{align*}
Second, observe that $\rb'$ is just the sample mean of $\rb_{j_1}, \dots, \rb_{j_K}$, i.e., $\rb' = (1/K) \sum_{k=1}^K \rb_{j_k}$. Lastly, we observe that the expected value of each $\rb_{j_k}$ is
\begin{align*}
\Ebb[ \rb_{j_k} ] & = \sum_{j \in [n]: \xi_j > 0} \xi_j \cdot \frac{ x^{*0}_{j} }{\xi_j} \alphab_j  \\
& = \sum_{j \in [n]: \xi_j > 0} x^{*0}_{j} \alphab_j  \\
& = \sum_{j \in [n]} x^{*0}_j \alphab_j \\
& = \rb^{*0},
\end{align*}
where the third step uses the fact that $x^{*0}_j = 0$ when $\xi_j = 0$ (by virtue of the constraint $\zerob \leq \xb \leq C \xib$). Therefore, the term $\| \rb' - \rb^{*0} \|_2$ is just the distance between the sample mean of an i.i.d. collection of random vectors from its expected value, where the $\ell_2$ norm of each random vector is bounded. We can therefore invoke Lemma~\ref{lemma:averaged_point_and_distance} to assert that
\begin{equation}
\| \rb' - \rb^{*0} \|_2 \leq \frac{CH}{\sqrt{K}} \left( 1 + \sqrt{2 \log \frac{2}{\delta} } \right) \label{bound:concentration_rb1_to_rb0}
\end{equation}
with probability at least $1 - \delta/2$. \\

\noindent \textbf{Bounding $\| \rb'' - \rb' \|_2$}: For this term, observe first that since $\rb'' = \rb' / (\oneb^T \xb')$, we can re-arrange this to obtain that $\rb' = (\oneb^T \xb') \rb''$. Let us use $s$ to denote the normalization constant, i.e., $s = \oneb^T \xb'$. We can now bound $\| \rb'' - \rb' \|_2$ in the following way:
\begin{align*}
\| \rb'' - \rb' \|_2 & = \left\| \rb'' -  s \rb''  \right\|_2 \\
& = | s - 1| \cdot \left\| \rb'' \right\|_2.
\end{align*}
We now bound $|s - 1|$. Note that $s$ can be written as 
\begin{equation*}
s = \oneb^T \xb' = \frac{1}{K} \sum_{k=1}^K \frac{x^{*0}_{j_k}}{ \xi_{j_k} } \oneb^T \eb_{j_k} = \frac{1}{K} \sum_{k=1}^K \frac{x^{*0}_{j_k}}{ \xi_{j_k} }. 
\end{equation*}
Letting $w_k = (x^{*0}_{j_k} / \xi_{j_k})$, we obtain $s = (1/K) \sum_{k=1}^K w_k$; in other words, $s$ is the average of $K$ i.i.d. random variables, $w_1, \dots, w_K$. Note that each $w_k$ has expected value $\Ebb[w_k] = \sum_{j \in [n]: \xi_j > 0 } (x^{*0}_{j} / \xi_{j}) \cdot \xi_j = \sum_{j \in [n]} x^{*0}_j = 1$; therefore, the term $|s - 1|$ represents how much the sample mean $s$ deviates from its expected value of 1. We also observe that each $w_k$ is contained in the interval $[0, C]$. Therefore, using Hoeffding's inequality, we obtain that 
\begin{equation}
\Pr[ |s - 1| > \epsilon] =\Pr[ |s - \Ebb[s]| > \epsilon] \leq 2 \cdot \exp\left( -\frac{2 K \epsilon^2}{C^2} \right),
\end{equation}
for any $\epsilon > 0$; by setting $\epsilon = C \sqrt{\log(4 / \delta) / (2K) }$, we obtain that
\begin{equation}
|s - 1| \leq C \sqrt{ \frac{1}{2K} \log \frac{4}{\delta} },
\end{equation}
with probability at least $1 - \delta/2$. 

With this bound in hand, let us now bound $\| \rb'' \|_2$. Observe that
\begin{equation*}
\| \rb' \|_2 \leq \frac{1}{K} \cdot \sum_{k=1}^K \left( \frac{x^{*0}_{j_k}}{ \xi_{j_k} } \right) \cdot \| \alphab_{j_k} \|_2 \leq \frac{1}{K} \cdot \sum_{k=1}^K \left( \frac{x^{*0}_{j_k}}{ \xi_{j_k} } \right) \cdot H = s \cdot H,
\end{equation*}
so it follows that $\| \rb'' \|_2 = (1/s) \| \rb' \|_2 \leq H$. We therefore have that $\| \rb'' - \rb' \|_2$ satisfies
\begin{equation*}
\| \rb'' - \rb' \|_2 \leq \frac{CH}{\sqrt{K}} \sqrt{ \frac{1}{2} \log \frac{4}{\delta} }, \label{bound:concentration_rb2_to_rb1}
\end{equation*}
with probability at least $1 - \delta/2$. \\

\noindent \textbf{Completing the proof}: We now put these two bounds together to complete the bound in \eqref{bound:portfolio_Lipschitz}. Combining inequalities~\eqref{bound:concentration_rb2_to_rb1} and \eqref{bound:concentration_rb1_to_rb0} together using the union bound, we have that with probability at least $1 - \delta$, 
\begin{align*}
F_J - F_{\distr} & \leq L\| \rb'' - \rb' \|_2 + L \| \rb' - \rb^{*0} \|_2 \\
& \leq L \cdot \frac{CH}{\sqrt{K}} \sqrt{ \frac{1}{2} \log \frac{4}{\delta} } + L \cdot \frac{CH}{\sqrt{K}} \left( 1 + \sqrt{2 \log \frac{2}{\delta} } \right) \\
& \leq \frac{CHL}{ \sqrt{K} }  \left( 1 + 3 \sqrt{ \log \frac{4}{\delta} } \right).
\end{align*}
By moving $F_{\distr}$ to the right hand side, and subtracting $F$ from both sides, we obtain the desired inequality. \hfill \Halmos

\endproof

\section{Statistically-Dependent Columns}
\label{sec:dependent_columns}
So far we have assumed that each column in the column-randomized linear program is sampled independently. In this section, we show how this assumption can be relaxed. We state our main performance guarantee in Section~\ref{subsec:dependent_columns_theorem}. In Section~\ref{subsec:dependent_columns_group_sampling}, we consider a specific non-i.i.d. column sampling scheme --  \emph{groupwise sampling} -- which has natural applications in problems such as Markov decision processes, and apply our guarantee from Section~\ref{subsec:dependent_columns_theorem} to this sampling scheme.Finally, in Section~\ref{subsec:dependent_columns_worep}, we develop a different type of guarantee for the case when columns are uniformly sampled without replacement.

\subsection{Performance Guarantees via Dependency Graph and Forest Complexity}
\label{subsec:dependent_columns_theorem}

We begin by assuming that the randomization scheme $\rho$ is such that $j_1, \dots, j_K$ still follow the distribution $\xib$, i.e., $\Pr[ j_k = t] = \xi_t$ for $k \in [K]$ and $t \in [n]$, but they are no longer independent. Thus, the indices $j_1,\dots, j_K$ are no longer an i.i.d. sample from $\xib$, and we require a different set of tools to analyze Algorithm~\ref{alg:main} and $\Delta v(P_J)$ in this setting. 

To analyze the column randomization method, we will make use of a specific concentration inequality from \cite{liu2019mcdiarmid}, which requires specifying the dependence structure of a collection of random variables through a specific type of graph. We thus begin by briefly defining the relevant graph-theoretic concepts. 

Given an undirected graph $G$, we use $V(G)$ to denote the vertices of $G$, and $E(G)$ to denote the edges of $G$. Given two vertices $u, v \in V(G)$, the edge between $u$ and $v$ is denoted by $\langle u, v \rangle$. We say that $u$ and $v$ are adjacent if $\langle u, v \rangle \in E(G)$. We say that $u$ and $v$ are non-adjacent if they are not adjacent. For two sets of nodes $U, V \subseteq V(G)$, we say that $U$ and $V$ are non-adjacent if $u$ and $v$ are non-adjacent for every $u \in U$ and $v \in V$. Lastly, a graph $G$ is a \emph{forest} if it does not contain any cycles, and is a \emph{tree} if it does not contain any cycles and consists of a single connected component. 

With this definitions, we now define the dependency graph, which is a representation of the dependency structure within a collection of random variables. 

\begin{definition}\it (Dependency graph)
	An undirected graph $G$ is called a \emph{dependency graph} of a set of random variables $X_1,X_2,\ldots,X_K$ if it satisfies the following two properties:
	\begin{enumerate}
		\item $V(G) = [K]$.
		\item For every $I,J \subseteq [K]$, $I \cap J = \emptyset$ such that $I$ and $J$ are non-adjacent, $\{ X_i\}_{i \in I}$ and $\{  X_j  \}_{j \in J}$ are independent.
	\end{enumerate}
\end{definition}

We now introduce the concept of a forest approximation from \cite{liu2019mcdiarmid}.

\begin{definition}\it (Forest approximation, \cite{liu2019mcdiarmid})
	Given a graph $G$, a forest $F$, and a mapping $\phi: V(G) \rightarrow V(F)$, we say that $(\phi, F)$ is a \emph{forest approximation} of $G$ if, for any $u, v \in V(G)$ such that $\langle u, v \rangle \in E(G)$, either $\phi(u) = \phi(v)$ or $\langle \phi(u), \phi (v) \rangle \in E(F)$.

\end{definition}
In words, a forest approximation is a mapping of a general graph $G$ to a smaller forest $F$ that is obtained by merging nodes in $G$. For a given node $v \in V(F)$, the set $\phi^{-1}(v)$ corresponds to the set of nodes in $V(G)$ that were merged to obtain the node $v$. Using the notion of a forest approximation, we can now define the forest complexity of a graph $G$. 

\begin{definition} \it (Forest complexity, \cite{liu2019mcdiarmid})
Let $\Phi(G)$ denote the set of all forest approximations of $G$. Given a forest approximation $(\phi, F)$, define $\lambda_{(\phi, F)}$ as
	\begin{align*}
	\lambda_{(\phi,F)} = \sum_{ \left\langle u,v \right\rangle \in E(F)} \left(  |  \phi^{-1}(u)  | + |  \phi^{-1}(v)  |  \right)^2 + \sum_{i=1}^k \min_{ u \in V(T_i)} |  \phi^{-1}(u)  |^2
	\end{align*}
	where $T_1, \dots, T_k$ is the collection of trees that comprise $F$. 	We call $\Lambda(G) = \min_{ (\phi,F) \in \Phi(G)  } \lambda_{(\phi,F)}$ the \emph{forest complexity} of $G$.
	\end{definition}

The forest complexity $\Lambda(G)$ quantifies how much the graph $G$ looks like a forest. Notice that $\Lambda(G) \geq |V(G)|$ for any graph $G$. In practice, we only need an upper bound on $\Lambda(G)$, rather than its exact value; we refer readers to \cite{liu2019mcdiarmid} for several examples on how $\Lambda(G)$ can be bounded. 

Given a dependency graph $G$ for the random indices in the set $J$, we now bound the optimality gap of the column-randomized linear program.

\begin{theorem}
	\label{thm:main_dependent_columns}
	Let $C$ be a nonnegative constant, define $P_\distr$ as in Theorem~\ref{thm:main_largest_abs_of_all_dual_BFS} and assume the random indices in $J$ follow the dependency graph $G$ with forest complexity $\Lambda(G)$. For any $\delta \in (0,1)$, with probability at least $1 - \delta$ over the sample $J$, the following holds: if $P_J$ is feasible and $\rank(\Ab_J) = m$, then
	\begin{align}
	\Delta v(P_J)   \leq \Delta  v(P_\distr) +  C \cdot \left(1 + m \gamma \| \Ab \|_{\max}  \right) \cdot \left( \sqrt{\frac{K + 2 | E(G)  |}{K^2}} +  \sqrt{  \frac{2 \Lambda(G) \log(2 / \delta) }{K^2}}  \right),
	\end{align}
	where $\gamma$ and $\| \Ab \|_{\max}$ are defined as in Theorem~\ref{thm:main_largest_abs_of_all_dual_BFS}. 
	
	Under the same conditions, with probability at least $1 - \delta$ over the sample $J$, the following holds: if $P_J$ is feasible and $\rank(\Ab_J) = m$, then
	\begin{align}
	\Delta v(P_J)   \leq \Delta  v(P_\distr) +  C \cdot \chi \cdot \left( \sqrt{\frac{K + 2 | E(G)  |}{K^2}} +  \sqrt{  \frac{2 \Lambda(G) \log(1 / \delta) }{K^2}}  \right),
	\end{align}
	where $\chi$ is defined as in Theorem~\ref{thm:main_reduced_cost}.
\end{theorem}

The proof (see below) follows by utilizing the McDiarmid inequality for dependent random variables from \cite{liu2019mcdiarmid}. We note that Theorem~\ref{thm:main_dependent_columns} is a generalization of Theorems~\ref{thm:main_largest_abs_of_all_dual_BFS} and \ref{thm:main_reduced_cost}. If $j_1,j_2,\ldots,j_K$ are independent, then the dependency graph $G$ has no edges, and thus $|E(G)| = 0$ and $\Lambda(G) = K$. Therefore, when each column is generated independently, the upper bounds in Theorem~\ref{thm:main_dependent_columns} are equivalent to the bounds in Theorem~\ref{thm:main_largest_abs_of_all_dual_BFS} and \ref{thm:main_reduced_cost}. 

We close this section by now proving Theorem~\ref{thm:main_dependent_columns}. Before we can prove Theorem~\ref{thm:main_dependent_columns}, we need to establish two auxiliary results. The first result is the analog of Lemma~\ref{lemma:averaged_point_and_distance} for a collection of possibly dependent random variables, formulated in terms of forest complexity. 

\begin{lemma}
	\label{lemma:mean_concentration_dependent_samples}
	Let $\wb_1,\wb_2,\ldots,\wb_K$ be $K$ random vectors with same distribution. Let $G$ be the dependency graph of $\wb_1,\wb_2,\ldots,\wb_K$. In addition, assume $\| \wb_k\|_2 \leq C$ for $k = 1,\ldots,K$. Let $\bar{\wb} = (1/K) \cdot \sum_{k=1}^K \wb_k$. Then for any $\delta \in (0,1)$, we have, with probability at least $1 - \delta$,
	\begin{align*}
	\|  \bar{\wb} - \Ebb \bar{\wb}  \|_2 \leq C \cdot \left(\sqrt{\frac{K + 2 \cdot | E(G)| }{K^2}} + \sqrt{   \frac{2 \cdot  \Lambda(G)}{K^2}  \cdot \log \frac{1}{\delta}   }  \right).
	\end{align*}
\end{lemma}
\proof{Proof of Lemma~\ref{lemma:mean_concentration_dependent_samples}:} 
Define a space $\Wcal \equiv \left\lbrace \zb \mid \| \zb \|_2 \leq C \right\rbrace$. Consider a scalar function $f: \Wcal^K \rightarrow \Rbb$ defined as
\begin{align*} 
f(\zb_1,\zb_2,\ldots,\zb_K) = \left\| \frac{1}{K} \left( \zb_1 + \zb_2 + \ldots + \zb_K \right)  - \Ebb \bar{ \wb } \right\|_2
\end{align*}
For any $ k \in [K]$ and any $\zb_1,\ldots,\zb_k,\ldots,\zb_K,\zb_k' \in \Wcal$, we have
\begin{align*}
| f(\zb_1,\ldots,\zb_k,\ldots,\zb_K) - f(\zb_1,\ldots,\zb_k',\ldots,\zb_K)  | \leq \frac{\|  \zb_k - \zb_k'  \|}{K}  \leq \frac{ 2 C }{K}.
\end{align*}
Therefore, $f$ has the bounded differences property (note that in \citealt{liu2019mcdiarmid}, this is referred to as the  $\mathbf{c}$-Lipschitz property; see Definition 2.1 of that paper). By Theorem 3.6 of \cite{liu2019mcdiarmid}, for any $\epsilon>0$, we have
\begin{align*}
\Pr \left[  f(\wb_1,\ldots,\wb_K)  - \Ebb f(\wb_1,\ldots,\wb_K) \geq \epsilon  \right] \leq \exp \left(  - \frac{K^2 \epsilon^2}{2C^2 \cdot \Lambda(G) }  \right)
\end{align*}

On the other hand, define $\ub_i = \wb_i - \Ebb \wb_i$. Then
\begin{equation*}
\Ebb \left[\ub_i^T \ub_j \right]= \left \{
\begin{aligned}
& \Ebb \left[ \wb_i^T \wb_j \right] - \| \Ebb \wb_i \|_2^2 \leq \Ebb \left[\| \wb_i\|_2 \| \wb_j\|_2 \right] \leq C^2   , && \text{if}\ i = j \text{ or } \left\langle i, j \right\rangle \in E(G), \\
&0, && \text{otherwise}.
\end{aligned} \right.
\end{equation*}
Therefore,
\begin{align*}
\Ebb \left[ f(\wb_1,\ldots,\wb_K)^2   \right] & = \left\| \frac{1}{K} \left( \wb_1 + \ldots + \wb_K \right)  - \Ebb \bar{ \wb } \right\|_2^2 \\
& = \frac{1}{K^2} \left(  \sum_{i,j \in [K]} \Ebb \left[ \ub_i^T \ub_j  \right]  \right) \\ & = \frac{1}{K^2} \left(  \sum_{i \in [K]} \Ebb \left[\ub_i^T \ub_i \right]  + \sum_{ \left\langle i,j \right\rangle \in E(G)} \Ebb \left[\ub_i^T \ub_j \right]  \right) \\
& \leq C^2 \cdot \frac{K + 2 |E(G)|}{K^2}.
\end{align*}
As a result,
\begin{align*}
\Ebb f(\wb_1,\ldots,\wb_K)  \leq \sqrt{ \Ebb f(\wb_1,\ldots,\wb_K)^2  } \leq C \cdot \sqrt{ \frac{K + 2 |E(G)|}{K^2} },
\end{align*}
where the first inequality comes from the concavity of square root function. With all the results above, we have
\begin{align*}
\Pb \left[  f(\wb_1,\ldots,\wb_K) -  C \cdot \sqrt{ \frac{K + 2 |E(G)|}{K^2} } \geq \epsilon \right] & \leq \Pb \left[   f(\wb_1,\ldots,\wb_K) -   \Ebb  f(\wb_1,\ldots,\wb_K)  \geq \epsilon \right] \\ & \leq \exp \left(  - \frac{K^2 \epsilon^2}{2C^2 \cdot \Lambda(G) }  \right)
\end{align*}
Let $\epsilon = \sqrt{ 2C^2 \Lambda(G) \log(1 / \delta) /K^2 }$. Then with probability at least $1 - \delta$, we have
\begin{align*}
f(\wb_1,\dots, \wb_K) \leq C \cdot \sqrt{ \frac{K + 2 |E(G)|}{K^2} } +  C \sqrt{  \frac{2 \cdot \Lambda(G)}{K^2} \log \left( \frac{1}{\delta}  \right)  }.
\end{align*} 
We thus prove the statement. \hfill \Halmos

From Lemma~\ref{lemma:mean_concentration_dependent_samples}, we can also straightforwardly prove the following result, which is the analog of Lemma~\ref{lemma:averaged_point_and_distance_L1} for possibly dependent random variables. 

\begin{corollary}
	\label{lemma:mean_concentration_dependent_samples_l1}
	Let $\wb_1,\wb_2,\ldots,\wb_K$ be $K$ random vectors of size $m$ and with same distribution. Let $G$ be the dependency graph of $\wb_1,\wb_2,\ldots,\wb_K$. In addition, assume $\| \wb_k\|_\infty \leq C$ for $k = 1,\ldots,K$. Let $\bar{\wb} = (1/K) \cdot \sum_{k=1}^K \wb_k$. Then for any $\delta \in (0,1)$, we have, with probability at least $1 - \delta$,
	\begin{align*}
	\|  \bar{\wb} - E\bar{\wb}  \|_1 \leq \sqrt{m} \cdot C \cdot \left(\sqrt{\frac{K + 2 \cdot | E(G)| }{K^2}} + \sqrt{   \frac{2 \cdot  \Lambda(G)}{K^2}  \cdot \log \frac{1}{\delta}   }  \right).
	\end{align*}
\end{corollary}

With these two results, we can now proceed with proving Theorem~\ref{thm:main_dependent_columns}. 

\proof{Proof of Theorem~\ref{thm:main_dependent_columns}:} We define $\xb^{*0}$ and construct random vectors $\wb_{j_1},\ldots,\wb_{j_K}$, $\bb_{j_1},\ldots,\bb_{j_K}$ as in the proof of Proposition~\ref{prop:objective_bound_by_dual_solution}; we note that this construction is valid even if there exists dependency between the indices $j_1$, $\ldots$, and $j_K$. We further define $\xb'$ as the sample mean of $\wb_{j_1},\ldots,\wb_{j_K}$ and $\bb'$ as the sample mean of $\bb_{j_1},\ldots,\bb_{j_K}$. By Proposition~\ref{prop:objective_bound_by_dual_solution} and Expression~\eqref{eq:objective_bound_on_both_x_and_b_abs}, we have
\begin{align}
\label{eq:bound_optimality_gap_by_x_and_b_dependent_samples}
\Delta v(P_J) \leq  \Delta v(P_\distr) + \|  \xb' - \xb^{*0}\|_2 + \| \pb_J^* \|_\infty \cdot \| \bb' - \bb \|_1.
\end{align}
By invoking Lemma~\ref{lemma:mean_concentration_dependent_samples}, with probability at least $1- \delta$,
\begin{align}
\label{eq:x_concentration_dependent_samples}
\|  \xb' - \xb^{*0}  \|_2 \leq C \cdot \left(\sqrt{\frac{K + 2 \cdot | E(G)| }{K^2}} + \sqrt{   \frac{2 \cdot  \Lambda(G)}{K^2}  \cdot \log \frac{1}{\delta}   }  \right).
\end{align}
Similarly, by Corollary~\ref{lemma:mean_concentration_dependent_samples_l1}, with probability at least $1- \delta$,
\begin{align}
\label{eq:b_concentration_dependent_samples}
\|  \bb' - \bb  \|_1 \leq \sqrt{m} \cdot C \cdot \| \Ab \|_{\max} \cdot \left(\sqrt{\frac{K + 2 \cdot | E(G)| }{K^2}} + \sqrt{   \frac{2 \cdot  \Lambda(G)}{K^2}  \cdot \log \frac{1}{\delta}   }  \right).
\end{align}
Combining inequalities~\eqref{eq:bound_optimality_gap_by_x_and_b_dependent_samples}, \eqref{eq:x_concentration_dependent_samples}, and \eqref{eq:b_concentration_dependent_samples} and applying the union bound, we conclude that, with probability at least $1 - \delta$, the following holds: if $P_J$ is feasible and $\rank(\Ab_J) = m$,
\begin{align}
\Delta v(P_J)   \leq \Delta  v(P_\distr) +  C \cdot \left(1 + m \gamma \| \Ab \|_{\max}  \right) \cdot \left( \sqrt{\frac{K + 2 | E(G)  |}{K^2}} +  \sqrt{  \frac{2 \Lambda(G) \log(2 / \delta) }{K^2}}  \right).
\end{align}
Similarly, by Proposition~\ref{prop:objective_bound_by_dual_solution} and inequality~\eqref{eq:key_step_of_proof_of_theorem_2}, we have
\begin{align}
\label{eq:bound_optimality_gap_with_only_x_dependent_samples}
\Delta v(P_J) \leq \Delta v(P_\distr) +  \chi \cdot \| \xb' - \xb^{*0}\|_2.
\end{align}
Combining with inequality~\eqref{eq:x_concentration_dependent_samples}, we conclude that, with probability $1 - \delta$, the following holds: if $P_J$ is feasible and $\rank(\Ab_J) = m$,
\begin{align}
\Delta v(P_J)   \leq \Delta  v(P_\distr) +  C \cdot \chi \cdot \left( \sqrt{\frac{K + 2 | E(G)  |}{K^2}} +  \sqrt{  \frac{2 \Lambda(G) \log(1 / \delta) }{K^2}}  \right),
\end{align}
which completes the proof. \hfill \Halmos
\endproof

\subsection{Groupwise Column Sampling}
\label{subsec:dependent_columns_group_sampling}

In many linear programs, we can naturally rearrange and group related columns together. For example, in the LP formulation of an MDP, one can collect columns associated with state $s$ into a set $\Gcal(s)$; the collection of all columns is simply the disjoint union $\bigcup_{s=1}^{n_s} \Gcal(s) $, where $n_s$ is number of states in the MDP and each $\Gcal(s) = \{ (s,a) \mid a \in [n_a] \}$. For such a problem, sampling $J = \{j_1, \dots, j_K\}$ independently from the complete collection of columns, i.e., from $[n_s] \times [n_a]$, may not be attractive. The reason for this is that we may sample the columns in such a way that we do not sample any columns corresponding to a particular state $\tilde{s}$; in such a scenario, the sampled problem $P_J$ will automatically be infeasible. 

In the presence of a natural group structure of the columns, rather than sampling columns in total across all $n$ columns, one could consider sampling $n_r$ columns from each group. In the MDP example, this would correspond to sampling $n_r$ columns (which correspond to state-action pairs) for each state $s$. The resulting column-randomized linear program $P_J$ corresponds to an MDP where there is a random set of $n_r$ actions out of the complete set of $n_a$ actions available in each state $s$. Most importantly, $P_J$ is guaranteed to be feasible. 

It turns out that our results for dependent columns can be used to study column-randomized LPs where columns are sampled by groups. We refer to such a mechanism as a \emph{groupwise randomization scheme} and define it formally below.

\begin{definition}
\it (Groupwise Randomization Scheme) Assume the set of indices $[n]$ can be organized into $n_\Gcal$ groups, i.e., $[n]$ is the disjoint union of sets $\Gcal_g$ for $g=1,2,\ldots,n_\Gcal$. Consider a randomization scheme $\rho$ such that (i) it samples indices in $n_r$ rounds of sampling; (ii) in each round, it samples $n_\Gcal$ indices as follows: for $i = 1,\ldots,n_\Gcal$, it first uniformly at random chooses an index $g_i$ from $[n_\Gcal] \setminus  \{ g_j \mid j \in [i-1] \}$ then samples an index from group $\Gcal_{g_i}$ according to a distribution $\xib^{g_i}$. We refer to such a randomization scheme $\rho$ as a \emph{groupwise randomization scheme}. 
	\label{definition:group_sampling}
	\end{definition}

Note that the randomization scheme $\rho$ samples $K = n_r n_\Gcal$ indices in total, and samples $n_r$ columns in each group. By design, each random index $j$ follows the distribution $\xib$, whose probabilities are given by 
\begin{align*}
\xi_t \equiv \Pr \left[  j = t \right] = \frac{1}{n_\Gcal} \sum_{ g \in [n_\Gcal]} \Ibb\{ t \in \Gcal_g \} \cdot \xi^{g}_t = \frac{1}{n_\Gcal} \cdot \xi^{\Gcal(t)}_{t}
\end{align*}
where $\Gcal(t)$ is the group to which column $t \in [n]$ belongs to. 

By using our general result for dependent columns (Theorem~\ref{thm:main_dependent_columns}), we obtain a specific guarantee for column-randomized LPs obtained by groupwise randomization schemes.

\begin{figure}
	\centering
	\includegraphics[width = 0.9\textwidth]{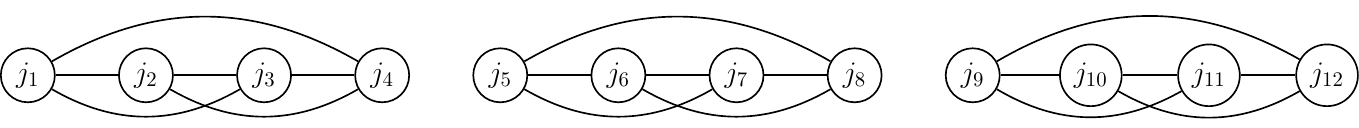}
	\caption{Dependency graph of random indices sampled by the groupwise randomization scheme with $n_\Gcal = 4$ and $n_r = 3$. 	\label{fig:group_sampling}}
\end{figure}

\begin{theorem}
	\label{thm:main_dependent_columns_group_sampling}
	Let $J$ be a sample of $K = n_r n_\Gcal$ indices sampled according to a groupwise randomization scheme $\rho$. Let $C$ be a nonnegative constant and define $P_\distr$ as in Theorem~\ref{thm:main_largest_abs_of_all_dual_BFS}. For any $\delta \in (0,1)$, with probability at least $1 - \delta$, the following holds: if $P_J$ is feasible and $\rank(\Ab_J) = m$, then
	\begin{align*}
	\Delta v(P_J)   \leq \Delta  v(P_\distr) + \frac{C \left( 1 + m \gamma \| \Ab \|_{\max} \right)}{\sqrt{n_r}} \left(  1 + \sqrt{2 \log \frac{2}{\delta}}  \right),
	\end{align*}
	where $\gamma$ and $\| \Ab \|_{\max}$ are defined as in Theorem~\ref{thm:main_largest_abs_of_all_dual_BFS}. Under the same assumption, with probability at least $1 - \delta$, the following holds: if $P_J$ is feasible and $\rank(\Ab_J) = m$, then
	\begin{align*}
	\Delta v(P_J)   \leq \Delta  v(P_\distr) + \frac{C \cdot \chi }{\sqrt{n_r}} \left(  1 + \sqrt{2 \log \frac{1}{\delta}}  \right),
	\end{align*}
	where $\chi$ is defined as in Theorem~\ref{thm:main_reduced_cost}.
\end{theorem}
{\it Proof:} The dependency graph $G$ of $K = n_r n_\Gcal$ random indices that are sampled by $\rho$ consists of $n_r$ cliques of size $n_\Gcal$; Figure~\ref{fig:group_sampling} provides an example of the dependency graph for $n_r = 3$ and $n_\Gcal = 4$. Therefore, $| E(G) | = n_r n_\Gcal (n_\Gcal -1) / 2$ and $\Lambda(G) \leq \lambda(\phi,F) = n_r n_\Gcal^2$ for a forest approximation $(\phi,F)$ that maps each clique in $G$ as a node in $F$. By upper bounding $\Lambda(G)$ by $n_r n_\Gcal^2$ in Theorem~\ref{thm:main_dependent_columns}, and using the fact that $K = n_r n_\Gcal$, we complete the proof. \hfill \Halmos

Theorem~\ref{thm:main_dependent_columns_group_sampling} can be interpreted as a guarantee on the optimality gap as a function of the number of columns sampled \emph{per group}: for a groupwise randomization scheme, the gap decreases at a rate of $1 / \sqrt{n_r}$, where $n_r$ is the number of columns sampled per group. Compared to Theorem~\ref{thm:main_largest_abs_of_all_dual_BFS} and \ref{thm:main_reduced_cost}, the rate of convergence in Theorem~\ref{thm:main_dependent_columns_group_sampling} in terms of the \emph{total} number of columns sampled, which is $K = n_r n_\Gcal$, is slower; Theorem~\ref{thm:main_largest_abs_of_all_dual_BFS} and \ref{thm:main_reduced_cost} both have a rate of $1 / \sqrt{K}$, while Theorem~\ref{thm:main_dependent_columns_group_sampling} has a rate of $1/\sqrt{n_r} \equiv \sqrt{n_\Gcal / K}$.

\subsection{Sampling without replacement}
\label{subsec:dependent_columns_worep}

The final extension of our methodology to the non-i.i.d. case that we shall consider is when the columns $\Ab_{j_1}, \dots, \Ab_{j_K}$ are sampled without replacement. For simplicity, we shall restrict our analysis to the case where this sampling is carried out uniformly over the set of columns $[n]$. Stated differently, a random sample of size $K$ drawn uniformly without replacement from $[n]$ is the set of columns $\{j_1,\dots,j_K\}$, where $\{j_1,j_2,\dots,j_n\}$ is a random permutation of the set of columns $[n]$, with all $n!$ permutations having equal probability.

For this sampling method, we begin with an analog of Lemma~\ref{lemma:averaged_point_and_distance}. This lemma uses results from the paper of \cite{el2009transductive}, which develops a version of McDiarmid's inequality that applies to the sampling without replacement case, and which may be of independent interest to readers.

\begin{lemma}
	\label{lemma:permutation_mcdiarmid}
	Assume $\wb_1,\ldots,\wb_n$ are vectors satisfying $\| \wb_j \|_2 \leq C $ for $j \in [n]$. Let $\{i_1,\ldots,i_n\}$ be a random permutation of $[n]$ and $\zb_{j} = \wb_{i_j}$ for $j \in [n]$. Define $\bar{\zb}_K = \sum_{j=1}^K \zb_j / K$ and $\bar{\zb} = \bar{\zb}_n = \Ebb \left[ \zb_1 \right]$. Then for any $\delta \in (0,1)$, we have, with probability at least $1 - \delta$,
	\begin{align}
	\label{eq:permutation_mcdiarmid}
	\| \bar{\zb}_K - \bar{\zb} \|_2 \leq \frac{C}{\sqrt{K}} \cdot \left(  \sqrt{ \frac{n-K}{n-1}} +  \sqrt{ \frac{2}{{H_{n,K}}} \log \left(  \frac{1}{\delta} \right)  }    \right),
	\end{align}
	where
	\begin{align*}
	H_{n,k} \equiv \frac{n - 1/2}{n-K} \cdot \left( 1 - \frac{1}{2 \max (K,n-K)} \right).
	\end{align*}
	\end{lemma}

\proof{Proof:} Call $\Zb = (\zb_1,\ldots,\zb_n) $. Define the function $f(\Zb) = \| \bar{\zb}_K - \bar{\zb}   \|_2 $, which is a $(K,n-K)$ \emph{permutation symmetric} function: that is, if we permute the first $K$ or the last $n-K$ vectors of $\zb_1,\ldots,\zb_n$, the value of $f(\Zb)$ remains the same. 

Given $\Zb$, let us use $\Zb^{ij}$ to denote the ordered collection that results from swapping the $i$th and $j$th vectors in $\Zb$. For $i \in \{  1,\ldots,K \}$ and $j \in \{ K+1,\ldots,n  \}$, we then have
\begin{align*}
| f(\Zb) - f(\Zb^{ij})  | \leq \frac{2C}{K},
\end{align*}
by the triangle inequality. Therefore, by Lemma 2 of \cite{el2009transductive}, we have
\begin{align*}
\Pr \left[  f(Z) - \Ebb f(Z) \geq \epsilon   \right] \leq \exp \left(  -  \frac{K\epsilon^2 }{2C^2} \cdot \frac{(n - 1/2)}{(n-K)} \cdot \left( 1 - \frac{1}{2 \max (K,n-K)} \right)    \right).
\end{align*}
Define $H_{n,k} $ as $\frac{(n - 1/2)}{(n-K)} \cdot \left( 1 - \frac{1}{2 \max (K,n-K)} \right)$. Therefore, the above inequality implies that with probability at least $1 - \delta$, we have
\begin{align}
\label{eq:sample_wo_replace_lemma_concentration}
f(Z) \leq \Ebb f(Z)  + \frac{C}{\sqrt{K}} \cdot \sqrt{ \frac{2}{{H_{n,K}}} \log \left(  \frac{1}{\delta} \right)  }.
\end{align} 

Now we will bound $\Ebb f(Z)$. We first define $\Ebb [ \| \zb_1 \| ^2] = a$ and $\Ebb [\zb_1^T \zb_2] = b$. Then
\begin{align}
\label{eq:variance_in_sample_wo_replacement}
\Ebb \left[  \| \bar{\zb}_K - \bar{\zb}  \|_2^2   \right] = \Ebb \left[ \bar{\zb}_K^T \bar{\zb}_K  \right] -\bar{\zb}^T \bar{\zb} = \frac{1}{K} \cdot a + \frac{K-1}{K} \cdot b - \bar{\zb}^T \bar{\zb}.
\end{align}
Notice that when $K = n$, the left-hand side of Equation \eqref{eq:variance_in_sample_wo_replacement} is zero. This leads to
\begin{align*}
b = \frac{n}{n-1} \cdot \bar{\zb}^T \bar{\zb} - \frac{a}{n-1}.
\end{align*}
Plugging this expression for $b$ back into equation~\eqref{eq:variance_in_sample_wo_replacement} and noticing that $a \leq C^2$, we have
\begin{align*}
\Ebb \left[  \| \bar{\zb}_K - \bar{\zb}  \|_2^2   \right] = \frac{a}{K} \cdot \left(  1 - \frac{K-1}{n-1} \right) -  \bar{\zb}^T \bar{\zb}  \cdot \left( 1 - \frac{n(k-1)}{k(n-1)}  \right) \leq \frac{C^2}{K} \cdot \left(  1 - \frac{K-1}{n-1} \right)
\end{align*}
Combining with Jensen's inequality, we have
\begin{align}
\label{eq:sample_wo_replace_Jensen}
\Ebb f(Z) \leq \sqrt{ \Ebb f^2(\Zb)    } = \sqrt{ \Ebb  \| \bar{\zb}_K - \bar{\zb}  \|_2^2  } \leq \frac{C}{\sqrt{K}} \cdot \sqrt{1 - \frac{K-1}{n-1}}.
\end{align}
Finally, we plug inequality~\eqref{eq:sample_wo_replace_Jensen} into inequality~\eqref{eq:sample_wo_replace_lemma_concentration}, which completes the proof. \hfill $\square$ \\

Using this lemma, we can now establish an analog of Proposition~\ref{prop:objective_bound_by_dual_solution}. In this proposition, we work with the distributional counterpart problem $P^{\worep}_{\distr} = \min \{  \cb^T \xb \mid \Ab \xb = \bb, \zerob \leq \xb \leq C/n \oneb  \}$, which is the distributional counterpart corresponding to the uniform distribution on $[n]$ (i.e., with $\xi_j = 1/n$ for all $j \in [n]$).

\begin{proposition}
	Let $C$ be a nonnegative constant and define $P_\distr^{\worep}$ as the linear program $\min \{  \cb^T \xb \mid \Ab \xb = \bb, \zerob \leq \xb \leq C/n  \}$. Let $Q = \{ q_1,\ldots,q_K  \} \subset [n]$ be a set of $K$ indices that are sampled uniformly at random from $[n]$ without replacement. For any $\delta \in (0,1)$, with probability at least $1-\delta$, the following statement holds: if $P_Q$ is feasible, then
	\begin{align*}
	\Delta v (P_Q) \leq \Delta v \left( P_\distr^{\worep}  \right) + \frac{C}{\sqrt{K}} \cdot \left(  1 + \| \pb \|_\infty \cdot m \cdot \| \Ab \|_{\max}    \right) \cdot \left( \sqrt{  \frac{n-K}{n-1}  } + \sqrt{  \frac{2}{H_{n,K}} \log \frac{ 2}{\delta}  }  \right),
	\end{align*}
	for any optimal dual solution $\pb$ of $P_Q$. \label{prop:objective_bound_by_dual_solution_worep}
	\end{proposition}

	\proof{Proof:} The proof follows a similar argument for the i.i.d. case (Proposition~\ref{prop:objective_bound_by_dual_solution}). Let $\xb^{0*}$ be an optimal solution to $P_\distr^{\worep}$. Consider the solution
	\begin{align*}
	\xb' \equiv \frac{1}{K} \sum_{k = 1}^K  n x^{0*}_{q_k} \eb_{q_k}  \equiv \frac{1}{K} \cdot \sum_{k=1}^K \zb_{k},
	\end{align*}
	where $\zb_{k} = n x^{0*}_{q_k} \eb_{q_k} $ for $k \in [K]$. We also define $\bb' = \Ab \xb'$.
	
	The vectors $\{  \zb_k \}_{k=1}^K$ have the following properties. First, for all $k \in [K]$, $\Ebb \left[  \zb_k  \right] = \xb^*$. Second, $\| \zb_k  \|_2 \leq C $ for all $k \in [K]$ since $0 \leq x^{0*}_{q_k} \leq C/n$. With these properties and recognizing that $\xb' = \bar{\zb}_K$, we can invoke Lemma~\ref{lemma:permutation_mcdiarmid} and assert that, with probability at least $1 - \delta/2$,
	\begin{align*}
	\|  \xb' - \xb^{0*}  \|_2 \leq \frac{C}{\sqrt{K}} \cdot  \left( \sqrt{  \frac{n-K}{n-1}  } + \sqrt{  \frac{2}{H_{n,K}} \log \frac{2}{\delta}     }  \right).
	\end{align*}
	With the similar argument in Step 2 of the proof of Proposition~\ref{prop:objective_bound_by_dual_solution}, it can be easily shown that with probability at least $1 - \delta/2$,
	\begin{align*}
	\| \bb' - \bb \|_1 \leq \frac{m \cdot C \cdot \| \Ab\|_{\max}}{\sqrt{K}} \cdot \left( \sqrt{  \frac{n-K}{n-1}  } + \sqrt{  \frac{2}{H_{n,K}} \log \frac{2}{\delta}     }  \right).
	\end{align*}
	With the concentration inequalities in hand, we can bound the objective value of $P_Q$ following the procedure in Step 3 of the proof of Proposition~\ref{prop:objective_bound_by_dual_solution}. \hfill \Halmos
\endproof

With this result, the following analog of Theorem~\ref{thm:main_largest_abs_of_all_dual_BFS} can be established for the uniform sampling without replacement case. The proof is identical to Theorem~\ref{thm:main_largest_abs_of_all_dual_BFS} and is omitted for brevity.

\begin{theorem}
	Let $C$ be a nonnegative constant and define $P_\distr^{\worep}$ as the linear program $\min \{  \cb^T \xb \mid \Ab \xb = \bb, \zerob \leq \xb \leq C/n  \}$. Let $Q = \{ q_1,\ldots,q_K  \} \subset [n]$ be a set of $K$ indices that are sampled uniformly at random from $[n]$ without replacement. For any $\delta \in (0,1)$, with probability at least $1-\delta$, the following statement holds: if $P_Q$ is feasible and $\rank( \Ab_Q) = m$, then
	\begin{align*}
	\Delta v (P_Q) \leq \Delta v \left( P_\distr^{\worep}  \right) + \frac{C}{\sqrt{K}} \cdot \left(  1 + \| \pb \|_\infty \cdot m \cdot \| \Ab \|_{\max}    \right) \cdot \left( \sqrt{  \frac{n-K}{n-1}  } + \sqrt{  \frac{2}{H_{n,K}} \log \frac{ 2}{\delta}  }  \right),
	\end{align*}
	where $\gamma$ and $\| \Ab \|_{\max}$ are defined as in Theorem~\ref{thm:main_largest_abs_of_all_dual_BFS}.
	\label{thm:main_largest_abs_of_all_dual_BFS_worep}
	\end{theorem}

Alternatively, we can also establish an analog of Theorem~\ref{thm:main_reduced_cost}. To do so, we require an analog of Proposition~\ref{prop:objective_bound_by_dual_slack}, which we formalize below. The proof of this result follows by straightforwardly combining elements of the proof of Proposition~\ref{prop:objective_bound_by_dual_slack} and Proposition~\ref{prop:objective_bound_by_dual_solution_worep} above, and is thus omitted. 

\begin{proposition}
	Let $C$ be a nonnegative constant and define $P_\distr^{\worep}$ as the linear program $\min \{  \cb^T \xb \mid \Ab \xb = \bb, \zerob \leq \xb \leq C/n  \}$. Let $Q = \{ q_1,\ldots,q_K  \} \subset [n]$ be a set of $K$ indices that are sampled uniformly at random from $[n]$ without replacement. For any $\delta \in (0,1)$, with probability at least $1-\delta$, the following statement holds: if $P_Q$ is feasible, then
	\begin{align*}
	\Delta v (P_Q) \leq \Delta v \left( P_\distr^{\worep}  \right) + \frac{C}{\sqrt{K}} \cdot \| \cb^T -\pb^T \Ab\|_2 \cdot \left( \sqrt{  \frac{n-K}{n-1}  } + \sqrt{  \frac{2}{H_{n,K}} \log \frac{  1}{\delta}  }  \right),
	\end{align*}
	for any optimal dual solution $\pb$ of $P_Q$. \label{prop:objective_bound_by_dual_slack_worep}
	\end{proposition}

Using this proposition, we can then easily obtain the following counterpart of Theorem~\ref{thm:main_reduced_cost} for the uniform sampling without replacement case. 

\begin{theorem}
	Let $C$, $P_\distr^{\worep}$, and $Q$ be as defined in Theorem~\ref{thm:main_largest_abs_of_all_dual_BFS_worep}. For any $\delta \in (0,1)$, with probability at least $1-\delta$, the following statement holds: if $P_Q$ is feasible and $\rank( \Ab_Q) = m$, then
	\begin{align*}
	\Delta v (P_Q) \leq \Delta v \left( P_\distr^{\worep}  \right) + \frac{C}{\sqrt{K}} \cdot \chi \cdot \left( \sqrt{  \frac{n-K}{n-1}  } + \sqrt{  \frac{2}{H_{n,K}} \log \frac{ 1}{\delta}  }  \right),
	\end{align*}
	where $\chi$ is an upper bound on $\| \bar{\cb} \|_2$ for every basic solution of the complete problem $P$. 
	\label{thm:main_reduced_cost_worep}
\end{theorem}

We conclude this section by offering a remark on how the bounds we have developed here compare to our earlier bounds for the i.i.d. case. In particular, we focus on Lemma~\ref{lemma:permutation_mcdiarmid}, which is the main building block of these results. In the i.i.d. case, the counterpart of Lemma~\ref{lemma:permutation_mcdiarmid} is Lemma~\ref{lemma:averaged_point_and_distance} (Lemma 4 of \cite{rahimi2009weighted}):
\begin{align*}
\| \bar{\wb}_K - \bar{\wb} \|_2 \leq \frac{C}{\sqrt{K}} \cdot \left(  1 +  \sqrt{ 2 \log \left(  \frac{1}{\delta} \right)  }    \right).
\end{align*}
We numerically compare the bound in Lemma~\ref{lemma:averaged_point_and_distance} (``i.i.d. bound'') to that of Lemma~\ref{lemma:permutation_mcdiarmid} (``permutation bound'') in Figure~\ref{Fig:Permutation_McDiarmid} below. We set $\delta = 0.1$, $n = 100$ and vary $K$. From this figure, we can see that (i) the permutation bound \eqref{eq:permutation_mcdiarmid} is always tighter than the standard McDiarmid inequality bound, which is under the i.i.d. assumption; and (ii) as $K$ gets closer to $n$, the improvement becomes larger. 

\begin{figure}[h!]
	\centering
	\includegraphics{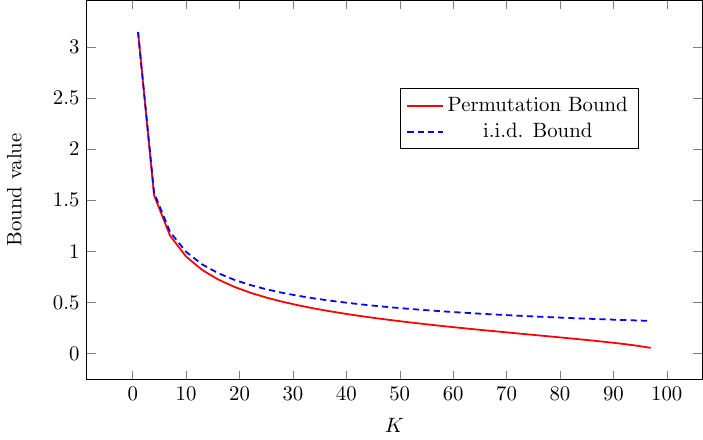}
	\caption{Comparison of the bounds of Lemma~\ref{lemma:averaged_point_and_distance} (based on the standard i.i.d. McDiarmid inequality) and Lemma~\ref{lemma:permutation_mcdiarmid} (which assumes uniform sampling without replacement), as $K$ varies.}
	\label{Fig:Permutation_McDiarmid}
\end{figure}

\section{Cutting Stock Problem Experiments (continued)}
\label{sec:cutting_stock_continued}

This section continues the numerical experiments with the cutting stock problem in Section~\ref{sec:numerics_CS}.

\subsection{Experiment \#2: comparison of incremental randomization and uniform randomization}
\label{subsec:numerics_CS_uniform_vs_incremental}

In this section, we explore the effect of changing the randomization scheme in the column randomization method. In particular, we compare the incremental randomization scheme $\rho_{\Incremental}$ of the previous section, and the uniform randomization scheme, which we will refer to by $\rho_{\UniformRS}$. 

The randomization scheme $\rho_{\UniformRS}$ samples from the set $\Acal = \{  \ab \in \Nbb^m_+ \mid \sum_{i \in [m]} a_i w_i \leq W  \}$ uniformly at random. This can be accomplished by rejection sampling. Specifically, we sample uniformly from the set $\bar{\Acal} = \{  \ab \in \Nbb^m_+ \mid 0 \leq a_i \leq \lfloor W / w_i \rfloor \}$, which can be done by sampling each component $a_i$ uniformly from the set $\{0,1,\dots, \lfloor W / w_i \rfloor \}$, and then check if $\sum_{i \in [m]} a_i w_i \leq W$. If this inequality is satisfied, we return $\ab$; otherwise, we discard $\ab$ and repeat the procedure again with a new candidate column from $\bar{\Acal}$. 

In this experiment, we set $W = 10^5$ again and draw each demand $b_i \sim U(\{1,\dots, 100\})$ for $i \in [m]$. Due to the poor scaling of rejection sampling, we restrict our focus to the case $m = 5$. 

For the widths $w_1,\dots, w_m$ of the demand types, we consider two different setups:
\begin{enumerate}
	\item \textbf{Setup 1}: we set each width $w_i \sim \Uniform( \{W/1000, W/1000 + 1, W/1000+2, \dots, W/2\})$. In this case, the widths can differ greatly, with two widths $w_i$ and $w_{i'}$ possibly differing by up to a factor of $(W/2) / (W/1000) = 500$. 
	\item \textbf{Setup 2}: we set each width $w_i \sim \Uniform( \{W/10, W/10 + 1, W/10+2, \dots, W/4\})$. In this case, the widths are generated to be closer to each other, with two widths $w_i$ and $w_{i'}$ only differing by up to a factor of $(W/4) / (W/10) = 2.5$. 
\end{enumerate}
We generate 100 random cutting stock instances in the manner described above for each setup. Then, for a fixed $K \in \{50, 100, 200, 400\}$, we run the column randomization method ten times with each of $\rho_{\UniformRS}$ and $\rho_{\Incremental}$.  

Table~\ref{tb:cutting_stock_uniform_sampling_comparison_setup12} below shows the optimality gap of column randomization with $\rho_{\Incremental}$ and $\rho_{\UniformRS}$ under Setups 1 and 2. The reported optimality gap is the average over the 100 cutting stock instances and the ten repetitions of the column randomization method. (Note that for $\rho_{\UniformRS}$ in Setup 1, there is no value shown for $K = 50$, as in one replication, the sampled problem was infeasible.) From this table, we can see that when there is high variability in the widths (Setup 1), $\rho_{\Incremental}$ outperforms $\rho_{\UniformRS}$ significantly. When there is a lower variability in the widths (Setup 2), $\rho_{\Incremental}$ generally outperforms $\rho_{\UniformRS}$, although the improvement is smaller. 

\begin{table}[!ht]
	\centering
	\begin{tabular}{rrrrr} \toprule
		& \multicolumn{2}{c}{Setup 1} &  \multicolumn{2}{c}{Setup 2} \\
		$K$   & $\Delta_{\rho_{\UniformRS}}$ (\%)  & $\Delta_{\rho_{\Incremental}}$ (\%)  & $\Delta_{\rho_{\UniformRS}}$ (\%)  & $\Delta_{\rho_{\Incremental}}$ (\%)  \\ \midrule
		50  & -- & 1.51 & 5.80 & 2.94 \\
		100 & 8.83 & 0.64 & 2.57 & 1.62  \\
		200 & 3.96 & 0.37 & 1.13 & 1.14  \\
		400 & 1.67 &  0.20 & 0.28 & 0.30 \\ \bottomrule
	\end{tabular}
	\caption{Comparison of $\rho_{\Incremental}$ and $\rho_{\UniformRS}$ on the cutting stock problem (Setups 1 and 2). 
		\label{tb:cutting_stock_uniform_sampling_comparison_setup12} }
	
\end{table}

The edge of $\rho_{\Incremental}$ over $\rho_{\UniformRS}$ is not surprising, when one considers the structure of the two randomization schemes. In particular, $\rho_{\Incremental}$ always produces patterns that are maximal, in the sense that no $a_i$ can be further incremented without violating the requirement $\sum_{i=1}^m w_i a_i \leq W$. (Note that such patterns will also be generated by the column generation subproblem~\eqref{problem:sub_cutting_stock}, as it seeks to maximize $\sum_{i=1}^m p_i a_i$ for some nonnegative dual vector $\pb$.) On the other hand, since $\rho_{\UniformRS}$ samples uniformly from $\Acal$, it frequently generates patterns that are inefficient, in that there is space on the large roll for more units of the demand types to be cut from it. Although such patterns can be used to meet the demands, one needs to cut more large rolls according to these patterns (i.e., the corresponding $x_j$'s need to be larger), resulting in a larger objective value. As a result, $\rho_{\Incremental}$ should yield lower optimality gaps than $\rho_{\UniformRS}$ for a fixed $K$. Nevertheless, this experiment is useful in showing that the choice of randomization scheme is important, and can substantially affect the performance of the column randomization method. We will further underscore this point in our next experiment, where we will also show how the demand vector $\bb$ can be used to guide the randomization scheme.

\subsection{Experiment \#3: comparison of incremental randomization and biased incremental randomization}
\label{subsec:numerics_CS_biased_vs_incremental}

In this next experiment, we compare incremental randomization with a more sophisticated scheme that we call biased incremental randomization and denote by $\rho_{\BiasedIncremental}$. This new scheme is presented as Algorithm~\ref{alg:sampling_cutting_stock_v2}. This scheme is the same as the incremental randomization scheme, with the key modification that at each iteration, the index $i$ is sampled with probability proportional to $\sqrt{b_i}$. 

The rationale behind this modification is as follows. Suppose that there is significant variability in the demands for different widths, e.g., for a width $i$, the demand $b_i$ could be very large, but for a different $i'$, the demand $b_{i'}$ could be very small. In such a situation, it may be advantageous to sample patterns where $a_i$ will tend to be large for highly demanded widths, while $a_i$ will tend to be small for less demanded widths, as patterns that are structured in this way are likely to be the most efficient patterns for meeting the demand. (Conversely, it is likely inefficient to use a pattern that yields a few units of the highly demand width and many units of the less demanded width.) The scheme $\rho_{\BiasedIncremental}$ is designed exactly for this case, and ensures that the patterns that are produced are such that $a_i$ will be larger when $b_i$ is large. 

\begin{algorithm}
	
	\begin{algorithmic}[1]
		\STATE Column $\ab$ is a zero vector of length $m$ and $\zeta \leftarrow W$.
		\WHILE {$\zeta > 0$}
		\STATE $I \leftarrow \{ i \mid w_i \leq \zeta \} $.
		\IF {$| I | \geq 1$}
		\STATE Sample an index $i$ from $I$ with probability $\sqrt{b_i} / \sum_{i \in I} \sqrt{b_i}$.
		\STATE Update $a_i \leftarrow a_i + 1$ and $\zeta \leftarrow \zeta - w_i$.
		\ELSE
		\STATE Break the while loop
		\ENDIF
		\ENDWHILE
		\RETURN Column $\ab$.
	\end{algorithmic}
	\caption{Biased incremental randomization scheme $\rho_{\BiasedIncremental}$ for the cutting stock problem. 	\label{alg:sampling_cutting_stock_v2}}
	
\end{algorithm}

In this experiment, we set $W = 10^5$. For each $i \in [m]$, we sample $w_i \sim \Uniform\left( \{ W/10,W/10+1,\dots,W/4 \}  \right)$. We then consider two different setups: {\bfseries Setup 1}, where each $b_i \sim \Uniform(\{ 25,\dots,100 \})$; and {\bfseries Setup 2}, where $b_i \sim \Uniform(\{ 50,\dots,100 \})$. For each setup, we test $m = 50$, $K \in \{100, 150, 200, 250, 300\}$, and $m = 100$, $K \in \{200, 250, 300, 350, 400\}$. For each $m$ and each setup, we generate 100 problem instances, and for each $K$, we then run column randomization with $\rho_{\Incremental}$ and $\rho_{\BiasedIncremental}$ ten times.

Table~\ref{tb:new_randomization_scheme_cutting_stock_setup12} shows the average optimality gap for the two randomization schemes $\rho_{\Incremental}$ and $\rho_{\BiasedIncremental}$ under the two different setups. As expected, we can see that in Setup 1, in which the demands exhibit greater variability, the biased scheme $\rho_{\BiasedIncremental}$ leads to a lower optimality gap than the ordinary incremental scheme $\rho_{\Incremental}$. In Setup 2, where there is less variability, $\rho_{\BiasedIncremental}$ continues to perform better, although the improvement is smaller.

\begin{table}[h]
	\centering
	\begin{tabular}{rrrrrr}
		\toprule
		& & \multicolumn{2}{c}{Setup 1} & \multicolumn{2}{c}{Setup 2} \\
		$m$   & $K$    & $\Delta_{\rho_{\Incremental}}$(\%)  & $\Delta_{\rho_{\BiasedIncremental}}$(\%) & $\Delta_{\rho_{\Incremental}}$(\%)  & $\Delta_{\rho_{\BiasedIncremental}}$(\%) \\ \midrule
		50 & 100 & 9.41 & 5.37 & 6.43 & 5.90 \\
		& 150 & 4.46 & 2.70 & 2.63 & 2.33 \\
		& 200  & 2.70 & 1.59 &1.71 &1.45  \\ 
		& 250 & 1.66 & 1.16 & 1.30 & 1.21 \\
		& 300 & 1.43 & 0.99 & 1.00 & 0.87 \\ \midrule 
		100 & 200 & 10.46 & 5.89 & 6.32 & 5.34  \\
		& 250 & 7.01 & 3.39 & 3.78 & 3.26 \\
		& 300 & 4.50 & 2.30 & 2.47 & 2.22 \\
		& 350 & 3.18 & 1.95 &  2.13 & 1.82 \\
		& 400 & 2.38& 1.56  & 1.61 & 1.50  \\ \bottomrule
	\end{tabular}
	\caption{Performance of the randomization schemes $\rho_{\Incremental}$ and $\rho_{\BiasedIncremental}$ on the cutting stock problem (Setups 1 and 2). \label{tb:new_randomization_scheme_cutting_stock_setup12}}
\end{table}

As we saw in our previous experiment in Section~\ref{subsec:numerics_CS_uniform_vs_incremental}, this experiment illustrates how the choice of randomization scheme can affect the performance of the column randomization method. It also illustrates how the structure of the problem and the nature of the problem data can affect the performance of column randomization and in the same vein, how the problem data can be used in the design of the randomization scheme (in this case specifically, how $\bb$ is used in $\rho_{\BiasedIncremental}$). 

\subsection{Experiment \#4: combining column randomization and column generation}
\label{subsec:numerics_CS_CRthenCG}

In this final experiment, we investigate the potential benefit of combining column randomization and column generation. In particular, we consider a hybrid method, where one first performs column randomization to obtain an initial solution, and then executes column generation starting from that initial solution. The hope in such a method is that column randomization can be used to quickly obtain a good solution with a low optimality gap, and that column generation can then be used to close that gap to zero. 

We set up this experiment as follows. We set $W = 10^5$. We vary $m \in \{250,500, 750, 1000, 1500\}$. For each $i \in [m]$, we draw $w_i \sim \Uniform(\{ W/10,\ldots,W/4 \})$ and $b_i \sim \Uniform(\{ 1,\ldots,100 \})$. For simplicity, we set number of sampled columns $K$ in the column randomization as $K=10m$. For the column randomization method, we use the incremental randomization scheme $\rho_{\Incremental}$. 

Table~\ref{tb:CR_then_CG_on_cutting_stock} displays the results. The columns labeled ``$\Delta_{\rho_{\BiasedIncremental}}$ (\%)'' and ``$T_{\CR}$ (s)'' show the optimality gap of the column randomization solution and the associated computation time. The next column, ``$T_{\text{CR-then-CG}}$ (s)'', shows the total time required to reach an optimality gap of zero when one executes column generation from the column randomization solution. The last column, ``$T_{\text{CG-only}}$ (s)'', shows the time required to reach an optimality gap of zero when one applies pure column generation. All values reported are averages over 100 randomly generated instances, and in the case of the columns that involve the column randomization method, are additionally averaged over ten repetitions of the column randomization method. 

\begin{table}[h!]
	\centering
	\small
	\begin{tabular}{ccrrrr}
		\toprule
		\multicolumn{1}{c}{{$m$}}    & \multicolumn{1}{c}{{$K$}} & \multicolumn{1}{c}{{$\Delta_{\rho_{\Incremental}}$(\%)}}     & \multicolumn{1}{c}{{$T_{\text{CR}}$ (s)}}          & \multicolumn{1}{c}{{$T_{\text{CR-then-CG}}$ (s)}}     & \multicolumn{1}{c}{{$T_{\text{CG-only}}$ (s)}}    \\
		\midrule  
		250 & 2500 & 1.733 & 0.403 & 35.391 & 52.138 
		\\
		500 & 5000 & 1.673 & 1.774& 97.370 & 130.223 \\
		750 & 7500 & 1.675 & 3.496 & 153.199 & 241.230 \\
		1000 & 10000 & 1.630 & 6.121  & 268.871  & 532.083 \\ 
		1500 & 15000 & 1.689 & 14.515 & 489.697 & 1107.038 
		\\ \bottomrule
	\end{tabular}
	\caption{Performance of the CG method on the cutting stock problem with and without the CR warm start}
	\label{tb:CR_then_CG_on_cutting_stock}
\end{table}

From this table, we can see that there is a benefit to combining column generation with column randomization. In particular, when $m = 500$, the improvement in the overall time required to reach a 0\% gap is small (at most about 30 seconds). For $m = 1000$ and $m = 1500$ the improvement is larger, with the combined method requiring roughly half of the time of the pure column generation method. This experiment illustrates that column randomization can serve as a simple and effective way to obtain an initial solution as an input to column generation, allowing the overall time to be significantly shortened.

\subsection{Experiment \#5: Exploration of optimal and near-optimal solution for a small instance}
\label{subsec:numerics_CS_small_instance}

In this section, we provide some more insight into why the column randomization method performs well on the cutting stock problem. We consider a small instance with $m = 8$ demand types, with large roll width $W = 200$ and the following widths and demands for the small rolls:
\begin{align*}
\wb & = (3,5,7,10,17,22,30,50), \\
\bb & = (1200,1000,1000,400,500,400,600,200).
\end{align*}
For this instance, the optimal objective value of problem $P^{\text{CS}}$ is 324.5. Our implementation of column generation returns the following solution consisting of 8 patterns that achieves this objective value (note that all $x_j$ values are given to four decimal places):
\begin{alignat*}{2}
\ab_1 & = (0, 40, 0, 0, 0, 0, 0, 0), &\quad\quad& x_1 = 19.7545 \\
\ab_2 & = (5, 1, 0, 0, 0, 0, 6, 0) && x_2 = 100 \\
\ab_3 & = (2, 0, 1, 0, 11, 0, 0, 0) && x_3 = 45.4545 \\
\ab_4 & = (0, 0, 0, 20, 0, 0, 0, 0) && x_4 = 20 \\
\ab_5 & = (0, 0, 0, 0, 0, 0, 0, 4) && x_5 = 50 \\
\ab_6 & = (30, 0, 0, 0, 0, 5, 0, 0) && x_6 = 18.1371 \\
\ab_7 & = (0, 2, 2, 0, 0, 8, 0, 0) && x_7 = 38.6643 \\
\ab_8 & = (2, 1, 27, 0, 0, 0, 0, 0) && x_8 = 32.4895
\end{alignat*}
However, this is not the only possible solution. To understand this better, we run the column randomization method 20,000 times, each time with $K = 100$ columns sampled according to $\rho_{\Incremental}$. For each run of the column randomization method, we solve for the optimal basic feasible solution of the sampled LP and we save the set $S \subseteq \{j_1,\dots, j_{100}\}$ of patterns for which $x_j > 0$. (Note that while $S$ can contain up to $m = 8$ indices, it could have fewer than 8 indices, because $P^{\text{CS}}$ is not a standard form LP.)

Over the 20,000 runs, we obtain 20,000 \emph{unique} sets of columns $S^1, \dots, S^{20000}$. Of these unique sets, \emph{5946 of them are optimal}. To illustrate, we list below five alternate optimal solutions:
\begin{alignat*}{2}
\text{Solution 1}: \quad & \ab_1 = (3, 6, 2, 2, 1, 0, 2, 1) &\quad\quad & x_1 = 68.9688\\
& \ab_2 = (5, 5, 2, 3, 2, 1, 2, 0) && x_2 = 37.4687\\
& \ab_3 = (6, 1, 2, 1, 3, 1, 1, 1) && x_3 = 59.4375\\
& \ab_4 = (3, 0, 2, 3, 1, 0, 1, 2) && x_4 = 1.2969\\
& \ab_5 = (3, 1, 5, 0, 1, 2, 3, 0) && x_5 = 88.3281\\
& \ab_6 = (2, 4, 4, 1, 2, 1, 1, 1) && x_6 = 26.75\\
& \ab_7 = (3, 4, 1, 2, 2, 0, 2, 1) && x_7 = 17.3281\\
& \ab_8 = (3, 3, 4, 1, 0, 4, 0, 1) && x_8 = 24.9219\\[0.5em]
\text{Solution 2}: \quad& \ab_1 = (3, 1, 1, 3, 1, 1, 2, 1) && x_1 = 36.268\\
& \ab_2 = (4, 2, 3, 2, 1, 0, 4, 0) && x_2 = 9.0103\\
& \ab_3 = (2, 6, 0, 0, 0, 2, 4, 0) && x_3 = 81.768\\
& \ab_4 = (3, 3, 2, 4, 0, 1, 0, 2) && x_4 = 46.7938\\
& \ab_5 = (7, 3, 3, 1, 3, 1, 2, 0) && x_5 = 77.7732\\
& \ab_6 = (3, 1, 9, 0, 3, 1, 0, 1) && x_6 = 67.1134\\
& \ab_7 = (2, 3, 1, 3, 4, 2, 1, 0) && x_7 = 2.7423\\
& \ab_8 = (0, 2, 1, 0, 3, 1, 2, 1) && x_8 = 3.0309\\[0.5em]
\text{Solution 3}: \quad & \ab_1 = (5, 4, 0, 2, 3, 2, 0, 1) && x_1 = 21.8947\\
& \ab_2 = (3, 8, 4, 4, 1, 3, 0, 0) && x_2 = 22.75\\
& \ab_3 = (4, 2, 0, 3, 4, 0, 1, 1) && x_3 = 17.4934\\
& \ab_4 = (6, 2, 5, 3, 1, 0, 3, 0) && x_4 = 7.0789\\
& \ab_5 = (3, 2, 3, 1, 0, 0, 0, 3) && x_5 = 11.2895\\
& \ab_6 = (4, 4, 2, 3, 0, 2, 1, 1) && x_6 = 26.7303\\
& \ab_7 = (4, 3, 5, 0, 2, 2, 2, 0) && x_7 = 117.25\\
& \ab_8 = (3, 2, 2, 1, 1, 0, 3, 1) && x_8 = 100.0132\\[0.5em]
\text{Solution 4}:\quad & \ab_1 = (1, 1, 2, 2, 2, 2, 1, 1) && x_1 = 4.1\\
& \ab_2 = (2, 2, 4, 1, 2, 1, 3, 0) && x_2 = 53.1\\
& \ab_3 = (1, 3, 3, 1, 3, 0, 0, 2) && x_3 = 10.1\\
& \ab_4 = (6, 4, 2, 1, 2, 2, 2, 0) && x_4 = 130.45\\
& \ab_5 = (2, 2, 2, 1, 0, 0, 2, 2) && x_5 = 48.95\\
& \ab_6 = (4, 0, 6, 1, 2, 1, 1, 1) && x_6 = 16.6\\
& \ab_7 = (4, 3, 2, 4, 1, 1, 1, 1) && x_7 = 5.1\\
& \ab_8 = (2, 4, 5, 2, 1, 1, 1, 1) && x_8 = 56.1\\[0.5em]
\text{Solution 5}:\quad & \ab_1 = (7, 5, 5, 0, 1, 1, 1, 1) && x_1 = 106.8125\\
& \ab_2 = (2, 2, 2, 1, 0, 0, 2, 2) && x_2 = 19.375\\
& \ab_3 = (2, 2, 4, 4, 2, 1, 2, 0) && x_3 = 44.3438\\
& \ab_4 = (3, 0, 0, 2, 1, 2, 2, 1) && x_4 = 3.6875\\
& \ab_5 = (4, 1, 2, 1, 1, 1, 4, 0) && x_5 = 54.2812\\
& \ab_6 = (1, 2, 1, 1, 4, 1, 1, 1) && x_6 = 50.0625\\
& \ab_7 = (3, 5, 0, 1, 2, 1, 0, 2) && x_7 = 0.3437\\
& \ab_8 = (1, 4, 2, 2, 1, 3, 2, 0) && x_8 = 45.5937
\end{alignat*}
There are two important points to note about these solutions. First, notice that all of the columns being used here are very different from the ones used in the column generation solution. In particular, the columns used in the column generation solution are sparser and the magnitudes of the $a_i$'s in those columns are larger. For example, in column $\ab_4$ of the CG solution, one cuts 20 units of demand type 4, whereas the most we cut of demand type 4 in any column of the above five alternate solutions is 4 (e.g., column $\ab_6$ in solution \#5). 

Second, the only overlap in the columns used in these five solutions comes from solution \#4 and solution \#5 (column $\ab_5$ of solution \#4 is the same as column $\ab_2$ of solution \#5). Apart from this one column that appears in two of the solutions, every other column only appears once. When we analyze the 5946 optimal solutions that we found, the corresponding columns sets together contain 5123 unique columns (i.e., letting $S^{(i)}$ denote the $i$th optimal column set, where $i$ ranges from 1 to 5946, we find $| S^{(1)} \cup S^{(2)} \cup \dots \cup S^{(5946)}| = 5123$). For each column, we calculate its incidence, which is the number of column sets in which the column appears. The maximum incidence of any column is 190, with the average incidence over all of the columns being 9.28 (i.e., on average each column appears in roughly 9 column sets). 

Building on the previous statement about the abundance of exactly optimal solutions, there exists an even greater number of \emph{near} optimal solutions. In particular, we can consider the number of solutions that are within $\epsilon = 2.0$ of the optimal objective; note that this absolute gap value translates to a relative gap of $2.0 / 324.5 = 0.62\%$. We find that that there are 18331 distinct column sets out of the 20,000 that are within $\epsilon = 2.0$ of the optimal objective of 324.5. These distinct column sets span 12294 unique columns, with the maximum incidence of any column being 820 and the average incidence being 11.88. 

These last two points, regarding the number of optimal and near-optimal solutions, are important because they directly relate to our analysis of the distributional counterpart in Section~\ref{sec:P_distr}. In particular, Theorem~\ref{theorem:many_BFS_P_distr_bound} of Section~\ref{subsec:P_distr_analysis} asserts that when there exist many $\epsilon$-optimal BFSs where the incidence of any column is low, then the distributional counterpart gap $\Delta v(P_{\distr})$ will be small. Note that although this result is formulated in terms of BFSs, the same proof technique goes through if one replaces these BFSs with solutions that are supported on a subset of the columns and for which any column appears in at most a certain number of supports. Thus, in the context of the cutting stock problem, it makes sense that column randomization does well, because optimal and nearly-optimal solutions that are diverse in terms of their columns exist in great profusion.

\section{Nonparametric Choice Model Estimation Experiments (continued)}
\label{sec:numerics_NCME_appendix}

This section continues the numerical experiments with the nonparametric choice model estimation in Section~\ref{sec:numerics_NCME}.

\subsection{Experiment \#2: Comparison of uniform randomization vs. MNL randomization}
\label{subsec:numerics_NCME_uniform_vs_MNL}

In our second experiment, we compare column randomization with the uniform randomization scheme $\rho_{\Uniform}$ against column randomization with an alternate randomization scheme that we refer to as \emph{MNL randomization}, and denote by $\rho_{\MNL}$. This randomization scheme involves first fitting an MNL model to the observed choice probabilities and then sampling rankings using the random utility model that underlies MNL. The procedure is formally defined below as Algorithm~\ref{alg:sampling_a_ranking_from_MNL}. %

\begin{algorithm}
	\caption{MNL randomization scheme $\rho_{\MNL}$ for the nonparametric choice estimation problem.}
	\label{alg:sampling_a_ranking_from_MNL}
	\begin{algorithmic}[1]
		\REQUIRE Estimated utilities $\hat{u}_1,\dots, \hat{u}_N$ of each product (via maximum likelihood estimation).
		\STATE Initialize $\alpha_{(i,m)} \leftarrow 0$ for $i \in [N]^+$ and $m \in [M]$.
		\STATE Generate $N+1$ independent random variable $\epsilon_i \leftarrow \text{Gumbel}(0,1)$, for $i \in [N]^+$.
		\STATE Set $v_i \leftarrow \hat{u}_i + \epsilon_i$ for $i \in [N]$, $v_0 \leftarrow 0 + \epsilon_0$.
		\STATE Set $\sigma$ to be the ranking such that $v_{\sigma(0)} > v_{\sigma(1)} > v_{\sigma(2)} > \ldots > v_{\sigma(N)} $.
		\FOR {$m \in [M]$}
		\STATE Set $i^* \leftarrow \arg \min_{ i \in S_m \cup \{ 0 \} } \sigma(i)$.
		\STATE Set $\alpha_{(i^*,m)} \leftarrow 1$
		\ENDFOR
		\RETURN Column $\alphab = (\alpha_{(i,m)})_{i \in [N]^+, m \in [M]}$.
	\end{algorithmic}
\end{algorithm}

In this experiment, we again vary $N$, $M$ and we consider two different setups. In Setup 1, we again sample the utility $u_i$ of each product $i$ as $u_i \sim \Uniform([0,1])$. In Setup 2, we instead sample the utility $u_i$ as $u_i \sim \Uniform([0,20])$. 

The rationale for Setup 2 is that when the magnitudes of the utilities $u_1,\dots,u_N$ are large, then the MNL model begins to behave more and more like a ranking based model. (To see this, suppose that $u_1,\dots, u_N$ is a collection of distinct non-zero real numbers; observe that for any set $S \subseteq [N]$ and $i \in S$, $\exp(\alpha u_i) / (1 + \sum_{i' \in S} \exp(\alpha u_i') ) \to \Ibb\{ i = \arg \max_{i' \in S \cup \{ 0\} } u_{i'} \}$ as $\alpha \to \infty$.) Thus, when the magnitudes of $u_1,\dots, u_N$ are large, the corresponding distribution over rankings will be concentrated around the ranking that corresponds to $u_1,\dots, u_N$, i.e., the $\sigma^*$ such that $u_{\sigma^*(0)} > u_{\sigma^*(1)} > u_{\sigma^*(2)} > \dots > u_{\sigma^*(N)}$. For this case, we should expect that $\rho_{\Uniform}$ will perform poorly, as it is unlikely that we will sample a large number of rankings around $\sigma^*$. On the other hand, we should expect $\rho_{\MNL}$ to perform better, as it samples rankings from a fitted MNL model whose parameters should be close to the parameters of the true underlying MNL model; thus, $\rho_{\MNL}$ should generate rankings that are close to $\sigma^*$. Conversely, in Setup 1 (which is identical to our prior setup), we should expect that $\rho_{\MNL}$ should improve over $\rho_{\Uniform}$, but the degree of improvement should be smaller. This is because when $u_1,\dots, u_N \sim \Uniform([0,1])$, the ranking distribution that corresponds to this MNL model will be more diffuse in the space of rankings.

Table~\ref{tb:two_CR_in_ranking_estimation_setup12} below shows the performance of column randomization equipped with the two randomization schemes in both Setup 1 and Setup 2. In the table, we use $Z_{\rho_{\Uniform}}$ and $Z_{\rho_{\MNL}}$ to denote the objective value of column randomization equipped with $\rho_{\Uniform}$ and $\rho_{\MNL}$ respectively. For Setup 1, we can see that in general, $\rho_{\MNL}$ does perform better than $\rho_{\Uniform}$; for example, for $N = 8$, $M = 100$, $K = 500$, the objective value of $\rho_{\MNL}$ is roughly half of that of $\rho_{\Uniform}$. For Setup 2, the edge of $\rho_{\MNL}$ over $\rho_{\Uniform}$ is more stark, with $\rho_{\MNL}$ resulting in objective values that are 2-3 orders of magnitude smaller than those of $\rho_{\Uniform}$. Note that the two forms of column randomization both have minimal computation time requirements: for $\rho_{\MNL}$, the time to carry out the maximum likelihood estimation, sample the columns, and solve $P^{\text{EST}}$ restricted to those columns is no more than 3 seconds across all $(N, M, K)$ combinations. Similarly, the end-to-end computation time for $\rho_{\Uniform}$ is also no more than 3 seconds across all $(N,M,K)$ combinations. 

\begin{table}[h]
	\centering
	\begin{tabular}{rrrrrrrrr} \toprule
		& & & \multicolumn{2}{c}{Setup 1} & \multicolumn{2}{c}{Setup 2} \\
		$N$ & $M$ & $K$ & $Z_{\rho_{\Uniform}}$ & $Z_{\rho_{\MNL}}$ & $Z_{\rho_{\Uniform}}$ & $Z_{\rho_{\MNL}}$ \\ \midrule
	6 &     50 &    500 & 0.06710 & 0.03665 & 4.66656 & 0.06257 \\ 
           &          &   1000 & 0.00013 & 0.00244 & 2.66596 & 0.02707 \\ \midrule
       8 &     50 &    500 & 0.12338 & 0.01776 & 6.89263 & 0.07756 \\ 
          &          &   1000 & 0.00000 & 0.00015 & 5.24286 & 0.03840 \\ \midrule
          &    100 &    500 & 1.04713 & 0.48256 & 17.48573 & 0.14387 \\ 
         &         &   1000 & 0.20209 & 0.02196 & 13.14833 & 0.09023 \\ 
         &        &   1500 & 0.00123 & 0.00103 & 10.09122 & 0.05453 \\ \midrule
      10 &     50 &    500 & 0.27030 & 0.07048 & 10.14161 & 0.08363 \\ 
         &        &   1000 & 0.00003 & 0.00001 & 7.85208 & 0.04030 \\ \midrule
         &    100 &    500 & 1.51977 & 0.91140 & 26.23535 & 0.10948 \\ 
        &       &   1000 & 0.35918 & 0.09530 & 18.81714 & 0.07296 \\ 
         &        &   1500 & 0.03308 & 0.00034 & 15.44066 & 0.05044 \\ 
         &        &   2000 & 0.00000 & 0.00027 & 13.44396 & 0.04081 \\ \bottomrule
	\end{tabular}
	\caption{Performance of $\rho_{\Uniform}$ and $\rho_{\MNL}$ under Setups 1 and 2. \label{tb:two_CR_in_ranking_estimation_setup12} }
	
\end{table}

As with our experiments with the cutting stock problem in Sections~\ref{subsec:numerics_CS_uniform_vs_incremental} and \ref{subsec:numerics_CS_biased_vs_incremental}, this experiment illustrates how the structure of the problem data can affect the performance of column randomization: in Setup 2, column randomization equipped with the basic randomization scheme $\rho_{\Uniform}$ performs poorly. Simultaneously, this experiment again offers an example of how one can use problem-specific knowledge to design the randomization scheme (in this case, fitting an MNL model, and then sampling from the ranking distribution corresponding to that fitted MNL model). We do acknowledge here that $\rho_{\MNL}$ is successful in this experiment because the ground truth model is an MNL model. For other ground truth models (e.g., the nested logit model or the latent-class MNL model), we should no longer expect $\rho_{\MNL}$ to do as well. However, for a different type of ground truth model, one can take the same strategy as in Algorithm~\ref{alg:sampling_a_ranking_from_MNL} where one estimates a different random utility maximization model. (For example, one could fit a latent-class MNL model using expectation-maximization, and then sample rankings from the resulting model.)

Lastly, we also note here that this experiment is congruent with our theoretical results on the distributional counterpart gap under generative model \GMDirichletNum (namely Theorem~\ref{theorem:GMDirichlet_gap_bound_whp}). Recall that in that generative model, the right-hand side is generated as a scaled random convex combination of the set of columns, where the vector of convex combination weights $\thetab$ is drawn uniformly from the $(n-1)$-dimensional unit simplex. Although problem $P^{\text{EST}}$ is not a standard form LP, there is a similarity here as the right hand side vector $\vb$ can also be thought of as being generated by a certain random combination of the columns in $\Ab = [ \alphab_1 \cdots \alphab_{(N+1)!} ]$, with the scale factor $\eta$ of generative model \GMDirichletNum being equal to 1 (see also the discussion in Section~\ref{subsec:P_distr_GM} around the moment estimation problem~\ref{prob:GMDirichlet_moment}). Thus, $\thetab$ can be thought of as the true underlying distribution over rankings. Theorem~\ref{theorem:GMDirichlet_gap_bound_whp} tells us that when $\thetab$ is drawn uniformly, and the randomization scheme is such that one samples columns uniformly over $[n]$, then most of the time the gap should be $O(\log n / \sqrt{n})$. The underlying distribution over rankings that one obtains under Setup 1 is closer to looking like a $\thetab$ drawn under generative model \GMDirichletNum than the same distribution obtained under Setup 2. This, in turn, explains why the performance of $\rho_{\Uniform}$ deteriorates so much from Setup 1 to Setup 2.

\subsection{Experiment \#3: combining column randomization and column generation}
\label{subsec:numerics_NCME_CRthenCG}

In this final experiment, analogously to the experiment in Section~\ref{subsec:numerics_CS_CRthenCG} for the cutting stock problem, we examine the value of using column randomization as a way of warm-starting column generation. We test the same values of $(N, M, K)$ as in Section~\ref{sec:numerics_NCME}. For each $N$ and $M$, we generate 100 problem instances in the same manner as in Section~\ref{sec:numerics_NCME}. Then, for each $K$, we run the column randomization method equipped with $\rho_{\Uniform}$ ten times, and we then use each solution as the initial solution for column generation, which we run until we reach an optimality gap of zero. 

Table~\ref{tb:CR_then_CG_on_ranking_estimation} shows the results of this experiment. The columns labeled $Z_{\rho_{\Uniform}}$ and $T_{\CR}$ indicate the objective value of the column randomization solution and the time required by column randomization, respectively. The next column, $T_{\text{CR-then-CG}}$, indicates the overall time required for the combined method (column randomization followed by column generation) to reach a zero optimality gap. The last column, $T_{\text{CG-only}}$, shows the time required for ordinary column generation (i.e., without any warm starting) to reach an optimality gap of zero. From this table, we again see that using column randomization to warm start column generation can dramatically reduce the time required to reach an optimality gap of zero. 

	\begin{table}
	\centering
	\begin{tabular}{rrrrrrr} \toprule
		$N$   & $M$  & $K$  & $Z_{\rho_{\Uniform}}$  & $T_{\text{CR}}$ (s) & $T_{\text{CR-then-CG}}$ (s)& $T_{\text{CG-only}}$ (s) \\ \midrule
 	      6 & 50 & 500 & 0.02990 & 0.06 & 3.94 & 27.55 \\ 
		    &    & 1000 & 0.00093 & 0.08 & 0.52 & 27.55 \\ \midrule
		  8 & 50 & 500 & 0.17719 & 0.17 & 21.77 & 107.27 \\ 
		    &     & 1000 & 0.00000 & 0.15 & 0.33 & 107.27 \\ \midrule
		    & 100 & 500 & 0.96841 & 0.32 & 397.90 & 705.70 \\ 
		    &     & 1000 & 0.19307 & 0.44 & 222.76 & 705.70 \\ 
		    &     & 1500 & 0.00000 & 0.92 & 1.63 & 705.70 \\ \midrule
		  10 & 50 & 500 & 0.32513 & 0.26 & 86.51 & 284.65 \\ 
		      &    & 1000 & 0.00008 & 0.24 & 2.09 & 284.65 \\ \midrule
		     & 100 & 500 & 1.48889 & 0.36 & 1441.17 & 2311.02 \\ 
		     &     & 1000 & 0.30858 & 0.53 & 737.89 & 2311.02 \\ 
		     &     & 1500 & 0.01132 & 1.06 & 133.56 & 2311.02 \\ 
		     &     & 2000 & 0.00000 & 2.60 & 3.87 & 2311.02 \\ \midrule
		     & 150 & 500 & 2.92142 & 0.82 & 6849.05 & 9831.54 \\ 
		     &     & 1000 & 1.16588 & 0.98 & 5358.42 & 9831.54 \\ 
		     &     & 1500 & 0.47822 & 1.44 & 3897.58 & 9831.54 \\ 
		     &     & 2000 & 0.15500 & 2.25 & 1987.44 & 9831.54 \\ 
		     &     & 2500 & 0.00000 & 5.82 & 10.41 & 9831.54 \\ \bottomrule
      \end{tabular}
      \caption{Performance of the CG method on the nonparametric choice model estimation problem with and without the CR-based warm start. \label{tb:CR_then_CG_on_ranking_estimation} }
      \end{table}

\section{Comparisons to Other Approaches}
\label{sec:comparison}

We complement Section~\ref{sec:literature_review} and make an additional comparison between our work and the two works in the literature.

\subsection{Comparison with \cite{agrawal2014dynamic}}
\label{subsec:comparison_agrawal2014dynamic}

The starting point of \cite{agrawal2014dynamic} is an \emph{online} linear program, which is an online version of the following problem:
\begin{subequations}
	\begin{alignat}{2}
	\label{problem:agrawal2014_offline_LP}
	P_{\OLP}: \quad & \underset{\xb}{\text{maximize}} & & \sum_{j=1}^n \pi_j x_j \\
	& \text{subject to} & \quad & \sum_{j=1}^n a_{ij} x_j \leq b_i, \quad \forall i \in [m], \\
	& & & 0 \leq x_j \leq 1, \quad \forall j \in [n].
	\end{alignat}
\end{subequations}
At the beginning, the decision maker has no information about the ground truth model, except knowing the total number of columns $n$. Time progresses in discrete periods, and at each period $t$, nature randomly reveals a column $\ab_t$ and a coefficient $\pi_t$ to the decision maker. The decision maker then makes a decision $x_t$ based on the history $(\pi_1,\ab_1,x_1,\pi_2,\ab_2,x_2,\ldots,\pi_t,\ab_t)$. Note that this is an \emph{irrevocable} decision: the decision maker cannot change the decisions made in earlier periods later. The decision maker's goal is to maximize the cumulative reward $\sum_{t=1}^n \pi_t x_t$. A critical assumption in \cite{agrawal2014dynamic} is that the arrival order of columns $(\ab_1,\ab_2,\ldots,\ab_n)$ is uniformly distributed over all the permutations (Assumption 1.1 of that paper) and $n$ is known (Assumption 1.2 of that paper).

To solve this online linear program, \cite{agrawal2014dynamic} propose an algorithm called the one-time-learning algorithm (OLA). This algorithm observes the first $s = \lceil \epsilon n \rceil$ periods of time, where $\epsilon \in (0,1)$, and using the $s$ columns observed, it defines a policy for making decisions for the remaining $n - s$ periods/columns. In particular, one considers following primal and dual problem pair defined on the first $s$ columns:
\begin{subequations}
	\begin{alignat}{2}
	P_{\OLA}:\quad & \underset{\xb}{\text{maximize}} & & \sum_{t=1}^s \pi_t x_t \\
	& \text{subject to} & \quad & \sum_{t=1}^s a_{it} x_t \leq (1 - \epsilon) \frac{s}{n} b_i, \quad \forall i \in [m], \\
	& & & 0 \leq x_t \leq 1, \quad \forall t \in [s]. 
	\end{alignat}
\end{subequations}
\begin{subequations}
	\begin{alignat}{2}
	D_{\OLA}: \quad & \underset{\pb, \yb}{\text{minimize}} & & \sum_{i=1}^m (1 - \epsilon) \frac{s}{n} b_i p_i + \sum_{t=1}^s y_t \\
	& \text{subject to} & \quad & \sum_{i=1}^m a_{it} p_i + y_t \geq \pi_t, \quad \forall t \in [s], \\
	& & & p_i \geq 0, \quad \forall i \in [m], \\
	& & & y_t \geq 0, \quad \forall t \in [s]. 
	\end{alignat}
\end{subequations}

Given a dual vector $\pb$ for the dual problem $D_{\OLA}$, define a policy $\xb(\cdot)$ as
\begin{equation}
x_t(\pb) = 
\begin{cases}
0, & \text{if}\ \pi_t \leq \pb^T \ab_t, \\
1, & \text{if}\ \pi_t > \pb^T \ab_t.
\end{cases}
\end{equation}
The OLA algorithm then operates as follows:
\begin{enumerate}
	\item[1)] Initialize $x_t = 0$ for all $t \leq s$. Let $\hat{\pb}$ be the optimal solution to the dual problem $D_{\OLA}$. 
	\item[2)] For $t = s+1,s+2,\ldots,n$, if $a_{it} x_t(\hat{\pb}) \leq b_i - \sum_{j=1}^{t-1} a_{ij} x_j$ for all $i \in [m]$, set $x_t = x_t(\hat{\pb})$; otherwise, set $x_t = 0$. Output $x_t$.
\end{enumerate}

Having provided this overview of the problem setup and the method of \cite{agrawal2014dynamic}, a number of critical differences become apparent. 

First, the problem setup in \cite{agrawal2014dynamic} is different from that of our paper. \cite{agrawal2014dynamic} considers an online problem: their problem is a sequential decision making problem such that at each time $t$, the decision maker receives information (a column) $\ab_t$, and then makes a decision $x_t$. Decisions made in the past cannot be changed. The decision maker cannot know more about the problem until nature reveals more information. In contrast, our problem is not an online problem: the problem is a static problem. %

Second, the source of randomness is different. In the model of \cite{agrawal2014dynamic}, at each time step $t$, nature reveals a column $\ab_t$ uniformly at random from the remaining columns; thus, the randomness is an inherent part of the problem. The OLA method is deterministic, that is to say, it does not introduce (additional) randomness to solve the problem. In other words, \cite{agrawal2014dynamic} uses a deterministic method to solve a stochastic problem, which is the online linear program. In contrast, in our paper, the ground truth model/problem is a large-scale linear program. This problem is a deterministic, one-shot problem -- there is no randomness in how information is revealed to the decision maker, and the decision maker does not need to set decision variables sequentially/in real time -- but it is {\it very large}. Therefore, we introduce randomness in the solution method, i.e., we propose a randomized algorithm to solve this large-scale deterministic problem.

Third, notwithstanding the difference in problem setups, \cite{agrawal2014dynamic} comment on the possibility of using OLA as an offline method to solve large-scale linear programs (see Section 5.3 of that paper). In particular, one first creates a random order of the $n$ columns, samples $s = \lceil \epsilon n \rceil$ columns, solves $D_{\OLA}$, and sets the variables according to the procedure given above. However, this approach is difficult to apply in the setting that we study. First, our linear program is a standard form LP of the form $\min \{ \cb^T \xb \mid \Ab \xb = \bb, \xb \geq \zerob \}$; in particular, the constraints are not inequalities, and the variables do not have a priori upper bounds, so it is not straightforward to adapt the variable-setting procedure of OLA to this more general problem. Second, even if one can overcome this difficulty, OLA fundamentally requires one to iterate through \emph{all} $n$ columns. This is impossible when $n$ is astronomically large. As an example, in the cutting stock problem that we study in Section~\ref{sec:numerics_CS} of our paper, OLA would involve sampling a small set of $s$ patterns, solving a problem to obtain dual variables, and iterating through \emph{every remaining pattern} to set $x_j$ for those patterns according to the dual variables. Although OLA could be useful for solving offline LPs where $n$ is moderately large -- i.e., the full LP is tedious to solve, but solvable -- we do not believe that it is computationally feasible for the case where $n$ is so large that the full LP itself cannot be formed and solved directly. This latter setting is precisely the setting that our method is intended for.

\subsection{Comparison with \cite{vu2018random}}
\label{subsec:comparison_vu2018random}

The approach of \cite{vu2018random} involves reducing the number of constraints in a linear program. 
In particular, instead of solving the problem $\min\{ \cb^T \xb \mid \Ab \xb = \bb, \xb \geq \zerob\}$, one forms a random $k$-by-$m$ matrix $\Tb$ and left-multiplies both sides of the constraint $\Ab \xb = \bb$ by this matrix to obtain the following simplified problem:
\begin{equation}
\min\{ \cb^T \xb \mid \Tb \Ab \xb = \Tb \bb, \xb \geq \zerob \}. \label{prob:vu_sketched_LP}
\end{equation}
This problem has fewer constraints ($k$ constraints, compared to $m$ constraints in the original problem). %

Having given an overview of the random projection method, it is clear that there are a number of important differences. First, while our method involves reducing the number of columns by drawing a random sample of columns, the random projection method of \cite{vu2018random} involves reducing the number of constraints by taking a linear combination of the constraints. This is important because in our problem setting, the number of columns $n$ is assumed to be much larger than the number of rows $m$; thus, random projection does not make the problem simpler to solve.

Second, by replacing the constraint $\Tb \Ab \xb = \Tb \bb$, the solution $\xb$ may not be feasible for the original equality constraint $\Ab \xb = \bb$. In fact, a result of \cite{vu2018random}, Proposition 3, asserts that a solution to problem~\eqref{prob:vu_sketched_LP} is infeasible for the original problem with probability 1. This contrasts with our setup, where if the sampled problem $P_J$ is feasible, the resulting solution is feasible for the complete problem $P$; and additionally, one can augment the sampled set of columns $J$ with a set of columns $J_F$ to guarantee feasibility of the sampled set (see Algorithm~\ref{alg:main_FG} in Section~\ref{subsec:discussion_theorems}).  Although \cite{vu2018random} provide a procedure (Algorithm 1) for retrieving an optimal basic feasible solution under certain conditions with high probability, the probability bound scales like $1 - O(n)$, which for our setting where $n$ is extremely large would yield a low probability. Indeed, the authors of \cite{vu2018random} acknowledge finding ``very high errors'' in applying this retrieval procedure in their numerical experiments (see the discussion in Section 7.2 of \citealt{vu2018random}), and for this reason consider a heuristic modification of their retrieval algorithm. Thus, guaranteeing a feasible solution to the original problem when applying the random projection method is not a triviality. 

Third, a tacit assumption in \cite{vu2018random} is that one can form the matrix $\Ab$ explicitly, and can carry out the multiplication $\Tb \Ab$ exactly. This will in general be impossible for the regime that we are interested, where $n$ can be astronomically large. For example, in the cutting stock example we consider, one would need to form the matrix $\Ab$ containing columns for \emph{all} possible patterns, and then compute $\Tb \Ab$; without even getting to the question of how one solves the sketched problem~\eqref{prob:vu_sketched_LP}, forming $\Ab$ and then $\Tb \Ab$ is clearly computationally infeasible. For this reason, the numerical examples that are considered in \cite{vu2018random} are of a much smaller scale than the ones we consider: $n$ is at most 2400, and the largest computation time reported for solving the original LP $\min\{ \cb^T \xb \mid \Ab \xb = \bb, \xb \geq \zerob\}$ is no more than two minutes.

Lastly, we comment that while the random projection method as originally described in \cite{vu2018random} reduces the number of rows in the primal LP, it is tempting to consider an alternate application of this method where one reduces the number of rows in the dual LP. The dual LP is 
\begin{equation*}
D : \max\{ \pb^T \bb \mid \pb^T \Ab \leq \cb^T \}
\end{equation*}
We can transform this into an equality constrained problem by introducing the slack vector $\sb \in \Rbb^n$: 
\begin{equation*}
D' : \max\{ \pb^T \bb \mid \pb^T \Ab + \sb^T = \cb^T, \sb \geq \zerob \}.
\end{equation*}
Now, we can right-multiply each side of the inequality constraint by a $n$-by-$k$ matrix $\Tb$, where $k \ll n$, resulting in the projected dual problem: 
\begin{equation*}
D'_{\text{RP}} : \max \{ \pb^T \bb \mid \pb^T \Ab \Tb + \sb^T \Tb = \cb^T \Tb \}.
\end{equation*}
The dual of this projected dual problem is 
\begin{equation*}
P_{\text{RP}} : \min\{ \cb^T \Tb \tilde{\xb} \mid \Ab \Tb \tilde{\xb} = \bb, \Tb \tilde{\xb} \geq \zerob \},
\end{equation*}
where $\tilde{\xb}$ is now a $k$-dimensional vector of decision variables (versus an $n$-dimensional vector in the original problem). Although this approach seems promising, again one runs into computation issues. %
In terms of computation, the matrix $\Tb$ is enormous as it has $n$ rows, and one needs to carry out the matrix multiplication $\Ab \Tb$, which for large-scale applications like cutting stock will be impossible. Additionally, although $P_{\text{RP}}$ achieves a reduction in the number of decision variables from $n$ to $k$, there are still $O(n)$ constraints due to the constraint $\Tb \tilde{\xb} \geq \zerob$. 
Lastly, a serious limitation of problem $P_{\text{RP}}$ is that it may be infeasible. (Comparing problem $P_{\text{RP}}$ and the original problem $P$, $P_{\text{RP}}$ is the same as $P$ with the constraint that $\xb$ lies in the lower dimensional subspace $\{ \Tb \tilde{\xb} \mid \tilde{\xb} \in \Rbb^k\}$. With this additional constraint, it is not guaranteed that we can satisfy the equality constraint $\Ab \xb = \bb$ and the nonnegativity constraint $\xb \geq \zerob$.)

We can therefore see that even applying random projection in an alternate fashion is problematic for the large-scale LP setting that we study.

\end{appendices}

\end{document}